\documentclass[reqno,11pt]{amsart}

\usepackage[english]{babel}
\usepackage[utf8]{inputenc}
\usepackage[T1]{fontenc}

\usepackage[margin=1in, letterpaper]{geometry}
\usepackage{amssymb,amsmath,amsfonts,amsthm}
\usepackage{bm,bbm,mathrsfs}
\usepackage[dvipsnames]{xcolor}
\usepackage{url}
\usepackage{enumitem}
\usepackage{mathtools}
\usepackage{placeins}
\usepackage{extarrows}
\usepackage{subfig}
\usepackage[font=small,labelfont=bf]{caption}
\usepackage{floatrow}
\usepackage{verbatim}
\usepackage{wrapfig}
\usepackage{tikz}
\usetikzlibrary{arrows,arrows.meta,cd,decorations.pathmorphing,backgrounds,positioning,fit,matrix}
\usepackage{xcolor,colortbl,graphics,graphicx}
\usepackage[colorlinks=true, allcolors=Blue, pdfencoding=auto]{hyperref}
\usepackage{aliascnt}
\usepackage[noabbrev,capitalize]{cleveref}

\newtheorem{theorem}{Theorem}[section]
\newaliascnt{lemma}{theorem}
\newtheorem{lemma}[lemma]{Lemma}
\aliascntresetthe{lemma}
\crefname{lemma}{Lemma}{Lemmas} 
\newaliascnt{proposition}{theorem}
\newtheorem{proposition}[proposition]{Proposition}
\aliascntresetthe{proposition}
\crefname{proposition}{Proposition}{Propositions}
\newaliascnt{corollary}{theorem}
\newtheorem{corollary}[corollary]{Corollary}
\aliascntresetthe{corollary}
\crefname{corollary}{Corollary}{Corollaries}
\newaliascnt{example}{theorem}
\newtheorem{example}[example]{Example}
\aliascntresetthe{example}
\crefname{example}{Example}{Examples}
\newaliascnt{question}{theorem}

\aliascntresetthe{question}
\crefname{question}{Question}{Questions}
\newaliascnt{conjecture}{theorem}

\aliascntresetthe{conjecture}
\crefname{conjecture}{Conjecture}{Conjectures}
\newaliascnt{assumption}{theorem}

\aliascntresetthe{assumption}
\crefname{assumption}{Assumption}{Assumptions}
\newaliascnt{definition}{theorem}

\aliascntresetthe{definition}
\crefname{definition}{Definition}{Definitions}
\newaliascnt{notation}{theorem}

\aliascntresetthe{notation}
\crefname{notation}{Notation}{Notations}
\theoremstyle{remark}
\newtheorem*{remark}{Remark}
\crefname{equation}{}{}
\numberwithin{equation}{section}
\allowdisplaybreaks

\newcommand{\rad}{\mathrm{rad}}
\newcommand{\Z}{\mathbb{Z}}
\newcommand{\R}{\mathbb{R}}
\newcommand{\Q}{\mathbb{Q}}
\newcommand{\C}{\mathbb{C}}

\renewcommand{\H}{\mathbb{H}}

\newcommand{\msD}{\mathscr{D}}
\newcommand{\msE}{\mathscr{E}}
\newcommand{\msF}{\mathscr{F}}

\newcommand{\msS}{\mathscr{S}}
\newcommand{\mE}{\mathcal{E}}
\newcommand{\mL}{\mathcal{L}}
\newcommand{\ma}{\mathfrak{a}}
\newcommand{\one}{\mathbbm{1}}

\renewcommand{\Re}{\textnormal{Re}}

\newcommand{\GL}{\textnormal{GL}}

\newcommand{\PSL}{\textnormal{PSL}}

\newcommand{\eps}{\varepsilon}

\renewcommand{\emptyset}{\varnothing}

\renewcommand{\bar}{\overline}
\renewcommand{\hat}{\widehat}
\renewcommand{\tilde}{\widetilde}

\renewcommand{\pmod}[1]{\ (\textnormal{mod } #1)}

\newenvironment{psmall}
  {\left(\begin{smallmatrix}}
  {\end{smallmatrix}\right)}

\definecolor{myBlue}{rgb}{0, 0, 0.6}
\newcommand{\fix}[1]{\textcolor{myBlue}{\large (#1)\normalsize}}

\title{Sums of Kloosterman fractions and applications}
\author[Alexandru Pascadi]{Alexandru Pascadi}
\address{Mathematisches Institut, Endenicher Allee 60, 53115 Bonn, Germany}
\email{alexpascadi@gmail.com}

\begin{document}

\begin{abstract}
    We obtain improved bounds for bilinear and trilinear forms with Kloosterman fractions, building on the works of Duke--Friedlander--Iwaniec and Bettin--Chandee. Our results use the spectral theory of automorphic forms, and save $N^{1/12} A^{1/3}$ over the trivial bound when $M = N$ and $A$ is small. In the exceptional Maass spectrum, we rely on a new large sieve inequality with averaging over the levels, which may be of independent interest.

    Applications include an asymptotic for the second moment of the zeta function twisted by a Dirichlet polynomial of length $T^{0.532}$ (improving Bettin--Chandee--Radziwi\l\l's $T^{0.515}$ and Dong--Robles--Zeindler's $T^{0.518}$), as well as an improved level of distribution for unbalanced convolutions (improving results of Fouvry--Radziwi\l\l, Wright, and Dong--Robles--Zeindler).
\end{abstract}

\maketitle

\vspace{-0.5cm}
{
\setlength{\parskip}{0em}
\setcounter{tocdepth}{1}
\tableofcontents
}
\vspace{-1cm}

\section{Introduction}

\subsection{Main results}

Let $M, N \in \Z_+$, $a \in \Z \setminus \{0\}$, and consider bilinear forms of the shape
\begin{equation} \label{eq:bilinear-forms}
    \mathop{\sum_{m=1}^M \sum_{n=1}^N}_{(m,n)=1}
    \alpha_m \beta_n\, e\left(a\frac{\bar m}{n}\right),
\end{equation}
where $e(t)=e^{2\pi i t}$, $\bar{m}m \equiv 1 \pmod{n}$, and $(\alpha_m), (\beta_n)$ are complex sequences. Sums of this kind occur in many arithmetic problems concerning moments of $L$-functions \cite{duke1997representations,bettin2017mean}, the distribution of convolutions in arithmetic progressions \cite{fouvry2022level}, and also level-aspect subconvexity\footnote{For recent bounds towards subconvexity for $\GL_2\times\GL_2$ by spectral reciprocity, see Humphries--Yang \cite{humphries2026rankin}.} \cite{aggarwal2026level}.

A bound for \cref{eq:bilinear-forms} which only depends on the $\ell^2$ norms $\|\alpha\|$, $\|\beta\|$ is equivalent to a bound for the largest singular value\footnote{When $M = N$ and $a$ is small, this matrix is essentially self-adjoint, so one can alternatively work with eigenvalues.} of the $M \times N$ matrix with entries $e(a \bar{m}/n) \one_{(m, n) = 1}$, giving one way to quantify randomness in the distribution of modular inverses. Indeed, this matrix has at most $\min(M, N)$ nonzero singular values and their sum of squares is $\approx MN$, so the largest singular value must be at least $\approx (M+N)^{1/2}$; this is essentially sharp for a matrix of independent random phases. One may therefore conjecture (see, e.g., \cite[(2.20)]{duke1997bilinear}) an upper bound of $(M+N)^{1/2+o(1)} |a|^{o(1)}$ for the largest singular value, which corresponds to square-root cancellation in either $M$ or $N$ in \cref{eq:bilinear-forms}.

However, even beating the trivial Cauchy--Schwarz bound of $\|\alpha\|\|\beta\|\sqrt{MN}$ for \cref{eq:bilinear-forms} has proven to be very challenging in the balanced regime $M \approx N$.
The first power saving in the balanced case was due to Duke--Friedlander--Iwaniec \cite[Theorem~2]{duke1997bilinear}, who saved a factor of $N^{1/48-o(1)}$ for fixed $a$; their proof used a clever amplification argument along with Cauchy--Schwarz, Fourier analysis, and the Weil bound for Kloosterman sums.
Bettin--Chandee \cite[Theorem~1]{bettin2018trilinear} improved this saving to $N^{1/20-o(1)}$, using a more efficient second Cauchy--Schwarz step. 

In this paper, we improve the results of Bettin--Chandee \cite{bettin2018trilinear} using the spectral theory of automorphic forms, including new ingredients in the exceptional spectrum which build on previous work of the author \cite{pascadi2026large}. More specifically, we prove a large sieve inequality for exceptional Maass forms with averaging over the levels and sequences that may depend on the level in a specific way; this is given in \cref{thm:exc-ls-level-avg}, which may be of independent interest. We discuss its use in \cref{subsec:outline}.

Duke--Friedlander--Iwaniec \cite[p.~35, after (3.27)]{duke1997bilinear} had raised the possibility of applying the spectral estimates of Deshouillers--Iwaniec \cite{deshouillers1982kloosterman} to the sums of Kloosterman sums obtained in their argument; they pointed out that the level of the resulting automorphic forms grows rapidly with the length of the amplifier (as well as with a square-full part parameter), making the size of the resulting bounds unclear. 
In the revised version of a preprint, which appeared in the same week as the present manuscript, Dong--Robles--Zeindler \cite{dong2026bilinear} independently pursued a spectral approach and improved the saving in the balanced-sum fixed-$a$ regime to $N^{1/16-o(1)}$. Here we obtain a saving of $N^{1/12-o(1)}$ uniformly in $|a|\ll N^{4/3}$ for the bilinear forms.

\begin{theorem} \label{thm:main-bilinear}
Let $N \in \Z_+$, $a \in \Z \setminus \{0\}$, and $(\alpha_m)_{m=1}^N$, $(\beta_n)_{n=1}^N$
be complex sequences. Then, uniformly in the range
$|a| \ll N^{4/3}$, one has
\[
    \mathop{\sum_{m=1}^N \sum_{n=1}^N}_{(m,n)=1}
    \alpha_m \beta_n\, e\left(a\frac{\bar m}{n}\right)
    \ll \|\alpha\|\|\beta\| N^{11/12+o(1)}.
\]
\end{theorem}

In applications, there is often an additional average over $a \le A$ with a third sequence $(\nu_a)$, which was first exploited by Bettin--Chandee \cite{bettin2018trilinear}. Our second result improves their savings (as well as the recent savings of Dong--Robles--Zeindler \cite{dong2026bilinear}) in both the $N$ and $A$ aspects.

\begin{theorem} \label{thm:main-trilinear}
Let $N, A \in \Z_+$ and $(\alpha_m)_{m=1}^N$,
$(\beta_n)_{n=1}^N$, $(\nu_a)_{a = 1}^A$ be complex sequences. Then, uniformly in the range $A \ll N^{1/14}$, one has
\[
\begin{aligned}
    \mathop{\sum_{m=1}^N \sum_{n=1}^N}_{(m,n)=1} \sum_{a=1}^A 
    \alpha_m \beta_n \nu_a\,
    e\left(a\frac{\bar m}{n}\right)
    &\ll \|\alpha\|\|\beta\|\|\nu\| N^{11/12+o(1)} A^{1/6}.
\end{aligned}
\]
\end{theorem}

The trivial Cauchy--Schwarz bound here is $\|\alpha\|\|\beta\|\|\nu\|N\sqrt A$. For balanced dyadic ranges $m,n\sim N$ and $A\ll N^{1/14}$, up to factors of $N^{o(1)}$, the savings over this bound given by Bettin--Chandee \cite[Theorem~1]{bettin2018trilinear}, Dong--Robles--Zeindler \cite[Theorem~1.1]{dong2026bilinear}, and \cref{thm:main-trilinear} are, respectively,
\[
     N^{1/20}A^{3/20},\qquad N^{1/16}A^{1/8},
     \qquad N^{1/12}A^{1/3}.
\]
For the dependence on larger $|a|$ and $A$, see \cref{cor:balanced-bilinear,cor:balanced-trilinear} respectively. In particular, for $N^{1/14} \ll A \ll N^{1/10}$, \cref{cor:balanced-trilinear} saves $N^{25/286} A^{79/286}$ over the trivial bound up to $N^{o(1)}$ factors, which continues to improve the saving of $N^{1/16} A^{1/8}$ from \cite{dong2026bilinear} in both aspects.

\cref{cor:balanced-spectral-methods} gives a more general result interpolating between the settings of \cref{thm:main-bilinear,thm:main-trilinear}, by considering a phase $e(\vartheta a\bar m/n)$ with an additional integer multiplier $\vartheta$ (with a sum over $a \le A$ but no sum over $\vartheta$); this follows the formulation of Bettin--Chandee \cite{bettin2018trilinear}.
We give more general estimates for potentially-distinct $M, N$ in \cref{cor:main-weil-bound-unbalanced,thm:unbalanced-spectral-methods}; in particular, \cref{cor:main-weil-bound-unbalanced} records an improvement in the $A$-aspect of the bounds of Bettin--Chandee \cite{bettin2018trilinear} using only the Weil bound. For other ranges of the three lengths, see also the work of Shparlinski \cite{shparlinski2021trilinear}.

\subsection{Applications} \label{subsec:applications}
Our first application concerns the second moment of the Riemann zeta function twisted by a Dirichlet polynomial with arbitrary coefficients. Let $T$ be large, let $w$ be a smooth function supported in $[1,2]$ with $w^{(j)}\ll_j(\log T)^j$ for $j\ge0$, and put
\[
    P_Y(s):=\sum_{n\le Y}\frac{c_n}{n^s},
    \qquad c_n\ll n^{o(1)}.
\]
%The implied constants may depend on the bounds for the coefficients and the derivatives of $w$.

\begin{corollary} \label{cor:appl-twisted-zeta}
Let $Y=T^\sigma$ with $0\le\sigma<207/389$ fixed. Then for some $\delta=\delta(\sigma)>0$ one has
\[
\begin{aligned}
    &\int_{\R}\left|\zeta(\tfrac12+it)P_Y(\tfrac12+it)\right|^2
           w(t/T)\,dt\\
    &\quad=\sum_{m,n\le Y}\frac{c_m\overline{c_n}}{[m,n]}
        \int_{\R}\left(\log\frac{t(m,n)^2}{2\pi mn}+2\gamma\right)
                 w(t/T)\,dt
        +O_\sigma(T^{1-\delta}),
\end{aligned}
\]
where $\gamma$ is Euler's constant.
\end{corollary}

Balasubramanian--Conrey--Heath-Brown \cite{balasubramanian1985asymptotic} obtained this asymptotic formula for $\sigma<1/2$, and Bettin--Chandee--Radziwi{\l}{\l} \cite[Theorem~1]{bettin2017mean} extended it to\footnote{Throughout this subsection, $\approx$ means truncation to 3 digits after the decimal point (so we round from below).} $\sigma < 17/33 \approx 0.515$. In the same range, Pratt--Robles \cite[Theorem~1.1]{pratt2018perturbed} allowed shifts of size $O(1/\log T)$ in the two zeta factors. More recently, Dong--Robles--Zeindler \cite[Theorem~1.2]{dong2026bilinear} obtained the shifted formula for $\sigma < 14/27 \approx 0.518$. Our corresponding range is $\sigma < 207/389 \approx 0.532$; we give the shifted formula with an error term uniform in $Y$ in \cref{thm:appl-twisted-zeta}.

Our second application is to the distribution of unbalanced convolutions in arithmetic progressions. 

\begin{corollary} \label{cor:appl-convolutions}
Let $k \in \Z_+$, $\eps>0$, $M,N\ge2$, $Q\ge1$ and $MN/2\le x\le4MN$. Let $(\alpha_m)_{m\sim M}$ and $(\beta_n)_{n\sim N}$ be complex sequences with $|\alpha_m|\le\tau_k(m)$ and $|\beta_n| \le \tau_k(n)$, where $\tau_k$ is the $k$-fold divisor function. Assume that $(\beta_n)$ satisfies the Siegel--Walfisz condition (recalled in \cref{subsec:appl-convolutions}). Suppose that $\exp((\log x)^\eps)\le N$ and $1\le |a|\le x^{175/389}$ with $a\in\Z$. Then, uniformly for
\[
    Q\le x^{207/389-\eps}N^{-468/389},
\]
one has, for every $B>0$,
\begin{equation} \label{eq:appl-convolutions}
    \sum_{\substack{q\sim Q\\(q,a)=1}}
    \Bigg|
        \mathop{\sum\sum}_{\substack{m\sim M, n\sim N\\x<mn\le2x\\mn\equiv a\pmod q}}
            \alpha_m\beta_n
        -\frac{1}{\varphi(q)}
        \mathop{\sum\sum}_{\substack{m\sim M, n\sim N\\x<mn\le2x\\(mn,q)=1}}
            \alpha_m\beta_n
    \Bigg|
    \ll_{k,\eps,B} \frac{x}{(\log x)^B}.
\end{equation}
%The implied constant may also depend on the constants in the Siegel--Walfisz condition.
\end{corollary}

Fouvry--Radziwi{\l}{\l} \cite[Corollary~1.1(i)]{fouvry2022level} obtained $Q\le x^{17/33-\eps}N^{-12/11}$, and Wright \cite[Corollary~2.2(i)]{wright2026trilinear} reduced the exponent $12/11 \approx 1.090$ in the denominator to $28/33 \approx 0.848$, by exploiting a fixed factor in the modulus of the Kloosterman fraction. Dong--Robles--Zeindler \cite[Theorem~1.4]{dong2026bilinear} obtained $|a|^6 Q^{27}N^{35}\le x^{14-\eps}$, giving an upper bound for $Q$ of roughly $x^{14/27}|a|^{-2/9}N^{-35/27}$. For fixed $a$ and $N=x^{o(1)}$, the limiting exponent of $Q$ thus increases from $17/33 \approx 0.515$ in the first two results and $14/27 \approx 0.518$ in the third to $207/389 \approx 0.532$ in \cref{cor:appl-convolutions}. The exponent of $N$ in our denominator is $468/389 \approx 1.203$ (compared to $35/27 \approx 1.296$), and the uniformity in the residue extends to $|a|\le x^{175/389}$; we treat larger residues in \cref{thm:appl-convolutions}.

Finally, consider the determinant equation $m_1n_2 - m_2n_1=\Delta$, with two smooth weights in $m_1, m_2$ and two arbitrary weights in $n_1, n_2$, as in the work of Duke--Friedlander--Iwaniec \cite{duke1997representations}. For fixed $f,g\in C_c^\infty((1,2))$ and complex sequences $(\alpha_n)_{n\sim N},(\beta_n)_{n\sim N}$, put
\[
    S_\Delta(N)=
    \sum_{\substack{n_1,n_2\sim N,\ m_1,m_2\in\Z\\m_1n_2-m_2n_1=\Delta}}
        \alpha_{n_1}\beta_{n_2} f(m_1/N)g(m_2/N).
\]

\begin{corollary} \label{cor:appl-determinant}
Let $N\ge2$. Uniformly for nonzero integers $\Delta$ with $|\Delta|\ll N^{4/3}$, one has
\[
\begin{aligned}
    \Bigg|S_\Delta(N)-
    \sum_{\substack{n_1,n_2\sim N\\(n_1,n_2)\mid\Delta}}
       \frac{(n_1,n_2)}{n_1n_2}\alpha_{n_1}\beta_{n_2}
       \int_{\R} f\left(\frac{x+\Delta}{Nn_2}\right)
                 g\left(\frac{x}{Nn_1}\right)\,dx\Bigg|
    \ll_{f,g} \|\alpha\|\|\beta\|N^{11/12+o(1)}.
\end{aligned}
\]
\end{corollary}

With the same factor of $\|\alpha\|\|\beta\|N^{o(1)}$, Bettin--Chandee \cite[Corollary~1]{bettin2018trilinear} obtained an error of $N^{19/20}$, while Dong--Robles--Zeindler \cite[Theorem~1.3]{dong2026bilinear} obtained $\min(N^{19/20},|\Delta|^{3/8}N^{15/16})$. %Thus for $\Delta=1$, the exponent decreases from $19/20 = 0.95$ and $15/16 \approx 0.937$, respectively, to $11/12 \approx 0.916$. 
We treat larger values of $\Delta$ in \cref{thm:appl-determinant}.

\subsection{Informal outline} \label{subsec:outline}
Let us briefly sketch the new ideas behind \cref{thm:main-trilinear}, ignoring many technical points; we use $\lesssim$ to indicate imprecise inequalities which can be made precise by adding smooth weights, GCD constraints, $N^{o(1)}$ factors, and similar details. In the critical ranges, applying Cauchy--Schwarz in $m$ eliminates the $\alpha_m$'s and leaves us to bound roughly 
\[
    \sum_{m \sim N} \bigg\vert \sum_{n \sim N} \sum_{a \sim A} \beta_n \nu_a\, e\left(a\frac{\bar{m}}{n}\right) \bigg\vert^2
    \stackrel{?}{\lesssim} 
    \|\beta\|^2 \|\nu\|^2 \frac{N^2 A}{S^2},
\]
where $S$ is the desired saving.
Following Duke--Friedlander--Iwaniec \cite{duke1997bilinear}, we introduce an amplifier supported on primes $\ell\sim L$. For technical reasons, we use an amplification setup closer to that of Maynard \cite{maynard2025primes}: we group the products $n\ell$ by their residue modulo $m$ (in fact, modulo a suitable multiple of $m$), and keep this out of Cauchy--Schwarz. On expanding the square, we obtain 
\[
    \sum_{m \sim N} \sum_{\substack{n_1, n_2 \sim N \\ \ell_1, \ell_2 \sim L \\ n_1\ell_1 \equiv n_2\ell_2 \pmod{m}}} 
    \beta_{n_1} \bar\beta_{n_2}
    \sum_{a_1, a_2 \sim A} \nu_{a_1} \bar\nu_{a_2}
    e\left(a_1\frac{\bar m}{n_1} - a_2\frac{\bar m}{n_2} \right)
    \stackrel{?}{\lesssim} 
    \|\beta\|^2 \|\nu\|^2 \frac{N A L^2}{S^2}.
\]
The diagonal terms cost roughly $\|\beta\|^2\|\nu\|^2 NL$, allowing a saving of $S \le (AL)^{1/2}$. Once we have dealt with the terms for which $n_1\ell_1=n_2\ell_2$ or $a_1\ell_1=a_2\ell_2$, we can replace $m$ by the nonzero complementary divisor
$d=(n_1\ell_1-n_2\ell_2)/m$, and focus on the critical range $d \sim L$. Up to a small Archimedean phase, reciprocity rewrites the exponential as
\[
    e\left(d(a_1\ell_1-a_2\ell_2)\frac{\overline{n_1\ell_1}}{n_2\ell_2}\right).
\]
We then apply Cauchy--Schwarz in $n_1,n_2$ to eliminate the $\beta_n$'s. The fact that we keep both $a_1$ and $a_2$ inside allows us to retain the saving $S \le (AL)^{1/2}$ from the diagonal terms, unlike in the work of Bettin--Chandee \cite{bettin2018trilinear}.
In the off-diagonal terms, writing $d',\ell_1',\ell_2',a_1',a_2'$ for the second copy of the inner variables, Poisson summation in $n_1$ brings us to sums of Kloosterman sums like
\begin{equation} \label{eq:sum-kl}
    \sum_{h \sim L^4} e\left(\frac{ht_q}{q}\right)
    \sum_{n_2\sim N}
        S\left(\Delta\bar r,h;sn_2\right),
\end{equation}
where
\begin{equation} \label{eq:sketch-phase}
    r=\ell_1\ell_1'dd',\quad s=\ell_2\ell_2', \quad q=rs, \quad 
    \frac{t_q}{q} \equiv \frac{\bar{\ell_1\ell_2'd'}}{d} + \frac{\bar{\ell_1'\ell_2d}}{d'} \pmod{1},
\end{equation}
\[
    \Delta=d\ell_1'\ell_2'(a_1\ell_1-a_2\ell_2)
    - d'\ell_1\ell_2(a_1'\ell_1'-a_2'\ell_2') \ll AL^4.
\]
The way that the phase $t_q/q$ depends on $d, d', \ell_1, \ell_1', \ell_2, \ell_2'$ will be relevant. Keeping track of the ranges in the Cauchy--Schwarz and Poisson steps, the remaining goal is to bound \cref{eq:sum-kl}, summed over $d, d', \ell_1, \ell_1', \ell_2, \ell_2', a_1, a_1', a_2, a_2'$, by roughly $\|\nu\|^4 N^2A^2L^8/S^4$.

The Kuznetsov trace formula relates the inner sums in \cref{eq:sum-kl} to the Fourier coefficients of automorphic forms $f$ of level $q$ and spectral parameters up to $O(1 + X^{-1})$, where in the critical range $\Delta \asymp AL^4$,
\[
    X := \frac{Ns\sqrt{r}}{\sqrt{|\Delta| L^4}} \asymp \frac{N}{\sqrt{A}} \gg 1.
\]
Let us first consider the regular spectrum where $f$ is a holomorphic form, Eisenstein series, or a Maass form obeying Selberg's eigenvalue conjecture. Applying Cauchy--Schwarz\footnote{One could, in principle, try to keep $a_1, a_1', a_2, a_2'$ inside the Cauchy--Schwarz step and use a second large sieve inequality instead of the density theorem, to pursue better savings in the $A$-aspect (up to $A^{1/2}$). However, the technical complications that arise seem extremely tedious, so we have decided not to pursue this improvement.} in $d, d', \ell_1, \ell_1', \ell_2, \ell_2'$, $a_1, a_1', a_2, a_2', f$, and considering each choice of $a_1, a_1', a_2, a_2'$ separately leaves us with bounding roughly
\[
    \bigg(\sum_{q \asymp L^6} \sum_f^{\Gamma_0(q)} |\lambda_f(\Delta)|^2 \bigg) \bigg(\sum_{q \asymp L^6} \sum_f^{\Gamma_0(q)} \bigg\vert \sum_{h \sim L^4} e\bigg(\frac{ht_q}{q}\bigg) \lambda_f(h) \bigg\vert^2 \bigg) \lesssim \frac{N^2 L^{20}}{S^8},
\]
where the sums are over orthonormal bases, and $\lambda_f$ are Fourier coefficients (or Hecke eigenvalues) normalized to be $1$-bounded on average. A standard density theorem (see \cite[Ch.~11]{iwaniec2021spectral}) bounds the sum in the first parenthesis by $L^{6+6}$ provided that $\Delta \lesssim L^{12}$, so $A \lesssim L^8$ suffices. The large sieve inequalities of Deshouillers--Iwaniec \cite[Theorem 2]{deshouillers1982kloosterman} bound the sum in the second parenthesis by $L^{6+6+4}$, so we need the condition
\[
    L^{12+16} \lesssim \frac{N^2 L^{20}}{S^8} \qquad \iff \qquad S^8 \lesssim \frac{N^2}{L^{8}}. 
\]
Combining this with the constraint $S \le (AL)^{1/2}$ from before, the optimal choice of the amplifier length is $L = N^{1/6} A^{-1/3}$. This gives a final saving of $S \approx N^{1/12} A^{1/3}$, as in \cref{thm:main-trilinear}.

We are not done:
in the \emph{exceptional spectrum} consisting of Maass forms which fail Selberg's eigenvalue conjecture, the spectral side of the Kuznetsov formula contains factors of roughly $X^{2\theta_f}$, where $\tfrac{1}{4}-\theta_f^2$ is the Laplacian eigenvalue of $f$ and $\theta_f > 0$. 
One has $\theta_f \le \tfrac{7}{64}$ due to Kim--Sarnak \cite[Appendix 2]{kim2003functoriality}, but the potential loss of $X^{2\theta_{\max}}$ (see \cref{eq:theta-max}) needs to be considered. Following \cite{pascadi2026large}, we split $X = X_0 \sqrt{X_1 X_2}$ and apply Cauchy--Schwarz, keeping $X_0^{2\theta_{\max}}$ outside and $X_1^{\theta_f}, X_2^{\theta_f}$ inside, and aiming to temper the contributions of $X_1$ and $X_2$ on average over $f$. The aforementioned density theorem (\cref{lem:density-thm}) handles $X_1 \Delta \lesssim L^{12}$, so we may take $X_1 \approx L^8/A$.

To temper the contribution of $X_2$, we need a bound of the shape
\begin{equation} \label{eq:sketch-large-sieve}
    \sum_{q \asymp L^6} \sum_{\substack{f \\ \lambda_f < \frac{1}{4}}}^{\Gamma_0(q)} X_2^{2\theta_f} \bigg\vert \sum_{h \sim L^4} e\bigg(\frac{ht_q}{q}\bigg) \lambda_f(h) \bigg\vert^2 
    \stackrel{?}{\lesssim} 
    L^{16},
\end{equation}
for $X_2$ as large as possible. Such estimates are known as exceptional-spectrum large sieve inequalities \cite{deshouillers1982kloosterman,pascadi2026large}; note that we have averaging over the levels $q$, but the sequence $e(ht_q/q)$ inside the square depends on the level. 
The results of Deshouillers--Iwaniec \cite[Theorem 5]{deshouillers1982kloosterman} would allow taking $X_2 = L^2$, while the author's previous work \cite[Theorem 1.5]{pascadi2026large} would allow $X_2 = L^4$, but neither of these previous results uses the average over $q$. In \cref{thm:exc-ls-level-avg}, we prove a large sieve inequality which exploits this averaging, allowing $X_2 = L^6$ in the setting above. This relies on the Kuznetsov formula and Poisson summation to reach a counting problem (similarly as in \cite{pascadi2026large}, but with a sum over $q$). The sum over $q$ is then absorbed by a divisor bound; see \cref{ex:tq}. 

By combining these exceptional-spectrum estimates, we can completely remove the dependency on $\theta_{\max}$ (i.e., we can take $X_0 = 1$) as long as
\[
    \sqrt{\frac{L^8}{A} \cdot L^6} = \sqrt{X_1 X_2} \gtrsim X = \frac{N}{\sqrt{A}}
    \qquad 
    \iff 
    \qquad 
    N \lesssim L^7.
\]
With the previous optimal choice $L = N^{1/6} A^{-1/3}$, this happens when $A \lesssim N^{1/14}$, which explains the range in \cref{thm:main-trilinear}. For larger values of $A$, \cref{cor:balanced-trilinear} gives a bound which depends on $\theta_{\max}$, but the large sieve inequality from \cref{thm:exc-ls-level-avg} is still useful in reducing this dependency.

\subsection{Some technical highlights} \label{subsec:technical}
We need to employ several technical tricks (not present in previous works on this topic) to formally justify the outline from \cref{subsec:outline}. Below are some examples:
\begin{itemize} 
\item[$(i)$.] We amplify using residue classes modulo $em$, with $e \sim E$ large enough that the complementary quotient $|n_1\ell_1-n_2\ell_2|/(em)$ is less than $L$. This is helpful at various points; e.g., it ensures that neither of $\ell_1, \ell_2$ divides the complementary quotient if it is nonzero.
\item[$(ii)$.] When $g=(n_1,n_2)>1$, we write $n_1=gu_1$, $n_2=gg_2u_2$, where $g_2\mid g^\infty$ and $(u_2,g)=1$. The reciprocity calculation then produces an extra phase modulo $gg_2$. We group terms according to a residue modulo $gg_2$ before the second Cauchy--Schwarz, so that equality of these residues supplies an additional congruence on expanding the square. 
\item[$(iii)$.] A box decomposition in the spirit of Fouvry \cite{fouvry1984autour} (see \cref{lem:box-decomposition}) rules out the `cross-diagonal' terms $\ell_1 = \ell_2'$ and $\ell_2 = \ell_1'$, ensuring some important coprimality relations.
\end{itemize}

\subsection{General notation}
We write $m\sim M$ for $M<m\le2M$, and $m\asymp M$ for $m$ bounded above and below by positive constant multiples of $M$. We use the standard asymptotic notation in analytic number theory, so $\ll$ and $O(\cdot)$ are interchangeable, and subscripts indicate the parameters on which the implied constants may depend. In expressions like $x^{o(1)}$, the exponent denotes a function that tends to $0$ as $x \to \infty$; in an upper bound, factors of $x^{o(1)}$ may be replaced by $O_\eps(x^\eps)$ for every $\eps>0$. In particular, we may absorb factors of the divisor function $\tau(n)$ into $x^{o(1)}$ whenever $1\le n \ll x^{O(1)}$.

We write $\tau_k$ for the $k$-fold divisor function.
We use $(m,n)$ for the greatest common divisor and $[m,n]$ for the least common multiple of $m,n$, and occasionally write $(m;n)$ for an ordered pair. We write $d\mid g^\infty$ if every prime divisor of $d$ divides $g$, and $\rad(d)$ for the product of the distinct primes dividing $d$ (so $d \mid g^\infty \iff \rad(d) \mid g$). Congruences modulo a negative integer are taken modulo its absolute value. The notation $r\pmod g^\times$ restricts $r$ to the invertible residue classes modulo $g$. If $r$ is given modulo $g$, then $(r,d)$ is well-defined for $d\mid g$, and the condition $(r,d)=1$ is well-defined more generally for $d\mid g^\infty$. For an integrable function $F$ on $\R^k$, we normalize its Fourier transform by
$\widehat F(\xi) := \int_{\R^k} F(x) e(-x\cdot\xi)\, dx$. 

Throughout this paper, we let $\Phi : \R \to [0, \infty)$ be a fixed smooth majorant of the indicator function of the interval $[\tfrac{1}{8},8]$, with the properties
\begin{equation} \label{eq:smooth-maj}
    \Phi \one_{[\frac{1}{8},8]} = \one_{[\frac{1}{8},8]}, \qquad \text{supp}(\Phi) \subset [0.1, 10], \qquad  \Phi^{(j)} \ll_j 1.
\end{equation} 
Further notation for automorphic forms is introduced in \cref{subsec:aut-forms-kuznetsov}.

\subsection*{Acknowledgements}
The author would like to thank James Maynard, Valentin Blomer, Jori Merikoski, and Edgar Assing for helpful discussions. The author was supported by the ERC Advanced Grant 101054336 and through Germany’s Excellence Strategy grant EXC-2047/1-390685813. 

\subsection*{AI Statement}
All of the ideas in this manuscript are due to the author, with the key results dating back to mid-2025. ChatGPT Astra was used to check for errors in an earlier version of this manuscript, and to assist with writing during the final stages of this project.

\section{Preliminaries} \label{sec:preliminaries}

\subsection{Standard arithmetic tools}
The following results will be used throughout our computations in \cref{sec:kloosterman-fractions-to-sums}.

\begin{lemma}[Swapping moduli] \label{lem:bezout}
For $a, b \in \Z_+$ and $x, y \in \Z$ such that $(a, b) = (x, ab) = 1$, one has
\[
    \frac{\bar{ax}y}{b} + \frac{\bar{bx}y}{a} \equiv \frac{\bar{x}y}{ab} \pmod{1}.
\]
\end{lemma}

\begin{proof}
This is often referred to as B\'ezout's identity, and follows from the Chinese Remainder Theorem.
\end{proof}

\begin{lemma}[Integers with specified prime factors] \label{lem:spec-prime-factors}
For any $g \in \Z_+$ and $D \ge 1$, one has
\[
    \#\{d \le D : d \mid g^\infty\} \ll (Dg)^{o(1)}.
\]
\end{lemma}
\begin{proof}
For any $\eps > 0$, we have
\[
    \#\{d \le D : d \mid g^\infty\}
    \le 
    \sum_{d \mid g^\infty} \left(\frac{D}{d}\right)^\eps
    =
    D^\eps
    \prod_{\text{prime }p \mid g} \frac{p^\eps}{p^\eps - 1}.
\]
Since we have $p^{\eps} \ge 2$ except for $O_\eps(1)$ primes $p$, we conclude that
\[
    \#\{d \le D : d \mid g^\infty\} \ll_\eps D^\eps \prod_{\text{prime }p \mid g} p^\eps \le (Dg)^\eps,
\]
which gives the desired bound.
\end{proof}

\begin{lemma}[GCD factors on average] \label{lem:gcd-avg}
Let $a, b \in \Z \setminus \{0\}$ and $M, N \gg 1$. Then one has %\fix{You might not need the second bound}
\[
    \sum_{n \sim N} (a, n) \ll |a|^{o(1)} N,
\]
\iffalse 
\begin{equation} \label{eq:gcd-avg-2}
    \sum_{\substack{m \sim M \\ n \sim N \\ am + bn \neq 0}} (am + bn, mn) \ll (|ab|MN)^{o(1)} MN.
\end{equation}
\fi
\end{lemma}
\begin{proof}
We write
\[
    \sum_{n \sim N} (a, n) 
    \le 
    \sum_{d \mid a} d \sum_{n \sim N} \one_{d \mid n}
    \ll 
    \sum_{d \mid a} d \frac{N}{d},
\]
and finish using the divisor bound. 
\iffalse
For \cref{eq:gcd-avg-2}, we let $g = (m, n)$ and put $g$ into dyadic ranges to rewrite our sum as
\[
\begin{aligned}
    \sum_{\substack{m \sim M \\ n \sim N \\ am + bn \neq 0}} (am + bn, mn)
    &=
    \sum_{g \le 2\min(M, N)} \sum_{\substack{m \sim M/g \\ n \sim N/g \\ (m, n) = 1 \\ am + bn \neq 0}} g(am + bn, gmn)
    \\
    &\ll 
    \min(M, N)^{o(1)}
    \max_{G \le \min(M, N)}
    G
    \sum_{\substack{\frac{M}{2G}\le m \le \frac{2M}{G} \\ \frac{N}{2G} \le n \le \frac{2N}{G} \\ (m, n) = 1 \\ am + bn \neq 0}} 
    \sum_{g \sim G} 
    (am + bn, gmn)
    .
\end{aligned}
\]
But when $(m, n) = 1$, we have
\[
\begin{aligned}
    (am + bn, gmn) 
    \le (am + bn, g) (am + bn, m) (am + bn, n)&
    \\
    =
    (am + bn, g) (bn, m) (am, n)&
    =
    (am + bn, g) (m, b) (n, a).
\end{aligned}
\]
Using \cref{eq:gcd-avg-1},
we can absorb the $(am + bn, g)$ factor when evaluating the sum over $g$, and then the $(m, b)$ and $(n, a)$ factors when evaluating the sums over $m$ and $n$. This gives
\[
\begin{aligned}
    \sum_{\substack{m \sim M \\ n \sim N \\ am + bn \neq 0}} (am + bn, mn)
    &\ll 
    \min(|ab|MN)^{o(1)}
    \max_{G \le \min(M, N)}
    G
    \frac{M}{G} \frac{N}{G} G
    \\
    &=
    \min(|ab|MN)^{o(1)}
    MN,
\end{aligned}
\]
which completes our proof.
\fi
\end{proof}

We recall the classical Kloosterman sums: for $c \in \Z_+$ and $m, n \in \Z$ (or $\Z/c\Z$), we write
\[
    S(m, n; c) := \sum_{x \pmod{c}^\times} e\left(\frac{mx+n\bar{x}}{c}\right).
\]

\begin{lemma}[Weil--Ramamanujan bound] \label{lem:weil-bound}
For $c \in \Z_+$ and $m, n \in \Z$ (or $\Z/c\Z$), one has
\begin{equation} \label{eq:weil-bound}
    S(m, n; c) \ll \tau(c) \sqrt{(m, n, c) c}
\end{equation}
and
\begin{equation} \label{eq:ram-bound}
    |S(0, n; c)| \le (n, c).
\end{equation}
\end{lemma}
\begin{proof}
The first bound is \cite[Corollary 11.12]{iwaniec2004analytic}, while the second bound (for Ramanujan sums) follows from M\"obius inversion; see \cite[(3.2) and (3.3)]{iwaniec2004analytic}.
\end{proof}

\begin{lemma}[Kloosterman completion] \label{lem:kloosterman-completion}
Let $c, q \in \Z_+$ with $(c, q) = 1$, $a \in \Z/c\Z$, and $b \in \Z/q\Z$. Let $N, x \gg 1$, and suppose that $c, q, N \ll x^{O(1)}$. Then for $\Phi$ as in \cref{eq:smooth-maj}, any $\eps > 0$, and any $H \gg x^\eps \frac{cq}{N}$, one has
\[
    \sum_{\substack{n \in \Z \\ (n, c) = 1 \\ n \equiv b \pmod{q}}}
    \Phi\left(\frac{n}{N}\right) e\left(a\frac{\bar{n}}{c}\right) 
    =
    \frac{N}{cq}
    \sum_{|h| \le H}
    \hat\Phi\left(\frac{hN}{cq}\right)
    e\left(h\frac{b\bar{c}}{q}\right)
    S(a\bar{q}, h; c)
    +
    O_\eps\left(x^{-100}\right).
\]
\end{lemma}

\begin{proof}
This is a standard consequence of truncated Poisson summation, which follows directly from \cite[Lemma 13.5]{maynard2025primes}, for example. 
\end{proof}

Combining \cref{lem:kloosterman-completion,lem:weil-bound} leads to a standard bound for incomplete Kloosterman sums.

\begin{corollary}[Incomplete Kloosterman bound] 
\label{cor:incomplete-kl-bound}
Let $c \in \Z_+$, $a \in \Z/c\Z$, and $N \gg 1$. Then for $\Phi$ as in \cref{eq:smooth-maj}, one has
\[
    \sum_{\substack{n \in \Z \\ (n, c) = 1}}
    \Phi\left(\frac{n}{N}\right) e\left(a\frac{\bar{n}}{c}\right) 
    \ll 
    (cN)^{o(1)} \left(\frac{N}{c} (a, c) + \sqrt{\frac{c}{(a,c)}}\right).
\]
\end{corollary}

\begin{proof}
This is \cite[Lemma 16.1]{maynard2025primes}.
\end{proof}

\begin{lemma}[$\GL_1$ large sieve] \label{lem:gl1-large-sieve}
Let $A \gg 1$ and $(\nu_a)_{a \sim A}$ be any complex sequence. Then for any $q \in \Z_+$ one has
\[
    \sum_{x \pmod{q}} \left\vert \sum_{a \sim A} \nu_a e\left(\frac{ax}{q}\right) \right\vert^2 \ll (q+A)\|\nu\|^2,
\]
and for any $Q \gg 1$ one has
\[
    \sum_{q \sim Q} \sum_{x \pmod{q}^\times} \left\vert \sum_{a \sim A} \nu_a e\left(\frac{ax}{q}\right) \right\vert^2 \ll (Q^2+A)\|\nu\|^2.
\]
\end{lemma}

\begin{proof}
By expanding the square, evaluating the sum over $x$, and Cauchy--Schwarz, the first sum is
\[
    q\sum_{a_1, a_2 \sim A} \nu_{a_1} \bar\nu_{a_2} \one_{a_1 \equiv a_2 \pmod{q}}
    =
    q \sum_{r \pmod{q}} \Bigg\vert \sum_{\substack{a \sim A \\ a \equiv r \pmod{q}}} \nu_a \Bigg\vert^2 
    \ll 
    q \left(1 + \frac{A}{q}\right) \sum_{a \sim A} |\nu_a|^2,
\]
which proves the first claim. The second claim is (a particular case of) the classical additive $\GL_1$ large sieve, see e.g. \cite[Theorem 7.11]{iwaniec2004analytic}.
\end{proof}

\begin{lemma}[Box decomposition] \label{lem:box-decomposition}
Let $X \ge 2$ and $B_X := ((X, 2X] \cap \Z) \times ((X, 2X] \cap \Z)$. There exist a positive integer $J_X \ll \log X$ and subsets $(B_{X,j})_{1 \le j \le J_X}$ of $B_X$ such that for each $j$,
\begin{equation} \label{eq:disjoint-projections}
    \left(\bigcup_{(x; y) \in B_{X,j}} \{x\} \right)
    \cap 
    \left(\bigcup_{(x; y) \in B_{X,j}} \{y\} \right)
    =
    \emptyset,
\end{equation}
and such that
\begin{equation} \label{eq:box-decomposition}
    \one_{B_X}
    =
    \one_{B_X \cap \{(x; x) : x \in \Z\}} +
    \sum_{j=1}^{J_X} \one_{B_{X,j}}.
\end{equation}
In other words, the box $B_X$ can be written as a disjoint union of its `diagonal' and the sets $B_{X,j}$, where each $B_{X,j}$ has disjoint coordinate-projections.
\end{lemma}

\begin{proof}
Our construction is a simpler version of Fouvry's \cite[Lemme 6]{fouvry1984autour}, and follows \cref{fig:Box-Dec}.
\begin{figure}[ht]
    \includegraphics[scale=0.5]{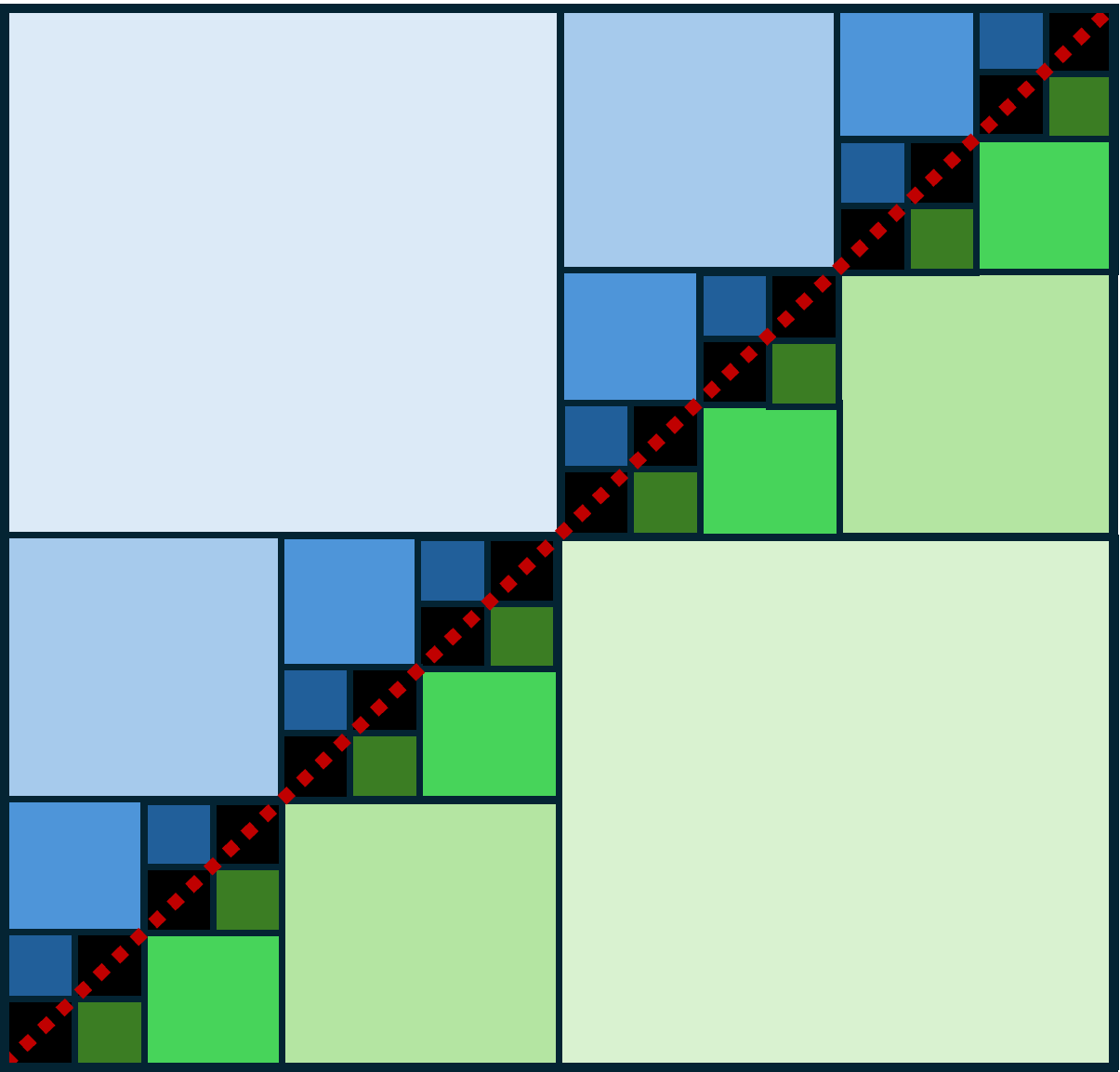}
\caption{Dyadic box decomposition. Each shade of blue and green corresponds to a different set $B_{X,j}$. Diagonal shown in red.}
\label{fig:Box-Dec}
\end{figure}

Let $B = (a, b] \times (a, b]$ be a box in $\R^2$, with $a, b \in \R$ such that $a < b$. For any such $B$, we define the four ``quarter''-boxes
\[
\begin{aligned}
    Q_1(B) &:= \left(a, \tfrac{a+b}{2}\right] \times \left(\tfrac{a+b}{2}, b\right], 
    \qquad\quad 
    Q_2(B) := \left(\tfrac{a+b}{2}, b\right] \times \left(\tfrac{a+b}{2}, b\right],
    \\
    Q_3(B) &:= \left(a, \tfrac{a+b}{2}\right] \times \left(a, \tfrac{a+b}{2}\right], 
    \qquad\quad 
    Q_4(B) := \left(\tfrac{a+b}{2}, b\right] \times \left(a, \tfrac{a+b}{2}\right].
\end{aligned}
\]
Note that 
\begin{equation} \label{eq:disjoint-coord-1}
    \left(\bigcup_{(x; y) \in Q_j(B)} \{x\} \right)
    \cap 
    \left(\bigcup_{(x; y) \in Q_j(B)} \{y\} \right)
    =
    \emptyset,
    \qquad\qquad 
    \text{for } j \in \{1, 4\},
\end{equation}
and that
\begin{equation} \label{eq:disjoint-coord-2}
    \left(\bigcup_{(x; y) \in Q_2(B)} \{x\} \cup \{y\} \right)
    \cap 
    \left(\bigcup_{(x; y) \in Q_3(B)} \{x\} \cup \{y\} \right)
    =
    \emptyset.
\end{equation}
Moreover, as indicator functions of subsets of $\R^2$, one has
\begin{equation} \label{eq:box-decomp-step}
    \one_B = \one_{Q_1(B)} + \one_{Q_2(B)} + \one_{Q_3(B)} + \one_{Q_4(B)}.
\end{equation}

Given $k \in \Z_+$ and indices $i_1, \ldots, i_{k-1} \in \{2, 3\}$ and $i_k \in \{1, 2, 3, 4\}$, we write $Q_{i_k,\ldots,i_1}(B) := Q_{i_k}(Q_{i_{k-1}}(\cdots Q_{i_1}(B) \cdots ))$. By recursively expanding $\one_{Q_2(B)}$ and $\one_{Q_3(B)}$ via \cref{eq:box-decomp-step}, we obtain that for any $J \in \Z_+$, one has
\begin{equation} \label{eq:box-decomp-recursive}
    \one_B = \left(\sum_{j=1}^{J} \sum_{i_1, \ldots, i_{j-1} \in \{2, 3\}} \left(\one_{Q_{1,i_{j-1},\ldots,i_1}(B)} + \one_{Q_{4,i_{j-1},\ldots,i_1}(B)}\right) \right)
    +
    \sum_{i_1, \ldots, i_{J} \in \{2, 3\}}
    \one_{Q_{i_J,\ldots,i_1}(B)}.
\end{equation}
We now take $B := (X, 2X] \times (X, 2X]$, $J := \lceil \log_2 X \rceil$, and for each $j \in \{1, \ldots, J\}$,
\[
    T_j := \bigcup_{i_1,\ldots, i_{j-1} \in \{2, 3\}} Q_{1,i_{j-1},\ldots,i_1}(B),
    \qquad\quad
    U_j := \bigcup_{i_1,\ldots, i_{j-1} \in \{2, 3\}} Q_{4,i_{j-1},\ldots,i_1}(B),
\]
which are disjoint unions. Note that each set $Q_{1,i_{j-1},\ldots,i_1}(B)$ has disjoint coordinate-projections due to \cref{eq:disjoint-coord-1}, and that two distinct sets of this shape have completely-disjoint coordinate-sets due to \cref{eq:disjoint-coord-2}. It follows that \cref{eq:disjoint-projections} holds with $B_{X,j}$ replaced by $T_j$, and similarly by $U_j$.

Now suppose we have $(x; y) \in Q_{i_J,\ldots,i_1}(B) \cap (\Z \times \Z)$ for some $i_1, \ldots, i_J \in \{2, 3\}$. An easy induction shows $|x-y| < 2^{-J} X \le 1$, which forces $x = y$. Since $T_j$ and $U_j$ do not contain any diagonal pairs, it follows that
\[
    (\Z \times \Z) \cap \bigcup_{i_1,\ldots,i_J \in \{2, 3\}} Q_{i_J,\ldots,i_1}(B)
    =
    B_X \cap \{(x; x) : x \in \Z\}.
\]
Multiplying \cref{eq:box-decomp-recursive} by $\one_{\Z \times \Z}$ now yields the desired decomposition in \cref{eq:box-decomposition}, up to relabeling the sets $(T_j \cap \Z^2)$, $(U_j \cap \Z^2)$ as $(B_{X,j})_{1 \le j \le J_X}$ where $J_X := 2J$.
\end{proof}

\subsection{Spectral estimates with automorphic forms}
\label{subsec:aut-forms-kuznetsov}

We now discuss the set-up for our spectral arguments from \cref{sec:kuznetsov-exc-spectrum}. We will mostly follow the notation of (and borrow some preliminary results from) \cite[\S 3.2--3.4]{pascadi2026large}, which built heavily on the work of Deshouillers--Iwaniec \cite{deshouillers1982kloosterman}.

For a positive integer $q$, we recall the congruence subgroup $\Gamma_0(q) = \{\pm \begin{psmall} a & b \\ c & d \end{psmall} \in \PSL_2(\Z) : q \mid c\}$. This group acts on the upper-half plane by M\"obius transformations, giving rise to the fundamental domain $\Gamma_0(q) \backslash \H$, with cusps of the shape $\tfrac{u}{w}$ where $w \mid q$ and $(u, w) = 1$. For a cusp $\ma \in \Q \cup \{\infty\}$ equivalent to such $\tfrac{u}{w}$, we write $\mu(\ma) := (w, \tfrac{q}{w})/q$, and also $\Gamma_\ma := \{\sigma \in \Gamma_0(q) : \sigma \ma = \ma\}$ for its stabilizer. For instance, $\mu(\infty) = q^{-1}$ since $\infty = \begin{psmall} 1 & 0 \\ -q & 1 \end{psmall} \tfrac{1}{q}$, and $\Gamma_\infty = \{\begin{psmall} 1 & n \\ 0 & 1 \end{psmall} : n \in \Z\}$. 
\iffalse 
All cusps with $\mu(\ma) = q^{-1}$ have the shape $\ma = \tau \tfrac{1}{s}$, where $\tau \in \Gamma_0(q)$ and $q = rs$ with $(r, s) = 1$; following \cite{watt1995kloosterman,pascadi2026large}, to such a cusp we associate the scaling matrix
\begin{equation} \label{eq:scaling-matrix}
    \sigma_{\ma}
    :=
    \tau 
    \begin{pmatrix}
        \sqrt{r} & -\frac{\bar{s}}{\sqrt{r}} \\ 
        s\sqrt{r} & \bar{r} \sqrt{r}
    \end{pmatrix},
\end{equation}
for some choice of $\bar{r}, \bar{s} \in \Z$ with $r\bar{r} + s\bar{s} = 1$. In particular, we choose $\sigma_\infty = \begin{psmall} 1 & 0 \\ -q & 1 \end{psmall} \begin{psmall} 1 & 0 \\ q & 1 \end{psmall} = \begin{psmall} 1 & 0 \\ 0 & 1 \end{psmall}$.
\fi 

The Kuznetsov trace formula \cite{deshouillers1982kloosterman,kuznetsov1980petersson} relates smooth sums over $c$ of Kloosterman sums $S(m\bar{r}, n; sc)$ to the data of automorphic forms of level $q = rs$, including their spectral parameters and their Fourier coefficients at $m, n$. Specifically, the spectral side of the Kuznetsov formula includes Maass cusp forms, holomorphic cusp forms, and Eisenstein series. Since our new inputs lie entirely in the exceptional spectrum, consisting of those Maass cusp forms which might fail the Selberg eigenvalue conjecture, we only need to introduce specific notation for the Maass spectrum.

The Maass cusp forms of level $q$ are smooth, $\Gamma_0(q)$-invariant, square-integrable functions $f \in L^2(\Gamma_0(q) \backslash \H)$, vanishing at the cusps, which are eigenforms of the hyperbolic Laplacian with eigenvalues $\lambda_f \in (0, \infty)$. The \emph{Selberg eigenvalue conjecture} \cite{selberg1965estimation,sarnak1995selberg} asserts that one always has $\lambda_f \ge \tfrac{1}{4}$, and the spectral parameter
\[
    \theta_f
    \in 
    i[0,\infty) \cup
    \left(0, \frac{7}{64}\right],
    \qquad 
    \theta_f^2 = \frac{1}{4} - \lambda_f
\]
measures the failure of Selberg's conjecture for a given form $f$; the upper bound of $\tfrac{7}{64}$ is due to Kim--Sarnak \cite[Appendix 2]{kim2003functoriality}, using the functoriality of fourth symmetric powers of $\GL_2$ forms.
We refer to those forms $f$ such that $\lambda_f < \tfrac{1}{4}$ (i.e., $\theta_f > 0$) as \emph{exceptional} -- conjecturally these exceptional forms do not exist, but unconditionally they produce loss factors of the shape $X^{\theta_f}$ in bounds for sums of Kloosterman sums. We will often write
\begin{equation} \label{eq:on-basis}
    \sum_f^{\Gamma_0(q)},
    \qquad 
    \text{respectively}
    \qquad 
    \sum_{\substack{f \\ \lambda_f < \frac{1}{4}}}^{\Gamma_0(q)},
\end{equation}
for a sum over an \emph{orthonormal basis} (arbitrary, unless specified otherwise) of Maass cusp forms on $\Gamma_0(q) \backslash \H$ with respect to the Petersson inner product, respectively for the restriction of such a basis to the exceptional spectrum. We also write
\begin{equation} \label{eq:theta-max}
    \theta_{\max} := \sup_{q \in \Z_+} \sup_f^{\Gamma_0(q)} \sqrt{\max\left(0, \frac{1}{4}-\lambda_f\right)}
    \in \left[0, \frac{7}{64}\right].
\end{equation}

Due to invariance under translation by $\begin{psmall} 1 & 1 \\ 0 & 1 \end{psmall}$, Maass cusp forms have Fourier expansions around $\infty$. More generally, one can Fourier-expand a Maass cusp form $f$ around any cusp $\ma$, using a scaling matrix $\sigma_{\ma}$ such that $\sigma_\ma \infty = \ma$ and $\sigma_{\ma}^{-1} \Gamma_\ma \sigma_{\ma} = \Gamma_\infty$. We write $\rho_{f\ma}(n)$ for the Fourier coefficient of $f$ at $n \in \Z \setminus \{0\}$ around $\ma$ with an implicit choice of $\sigma_\ma$, normalized as in \cite{deshouillers1982kloosterman,pascadi2026large} (in this paper, we will mainly use cusps with $\mu(\ma) = q^{-1}$, and one should think of $\rho_{f\ma}(n)$ as having size $\approx q^{-1/2}$ on average). Specifically, for $z = x+iy \in \H$, we have
\[
    f(\sigma_\ma z) = y^{1/2} \sum_{n \in \Z\setminus \{0\}} \rho_{f\ma}(n)\, K_{\theta_f}(2\pi|n|y)\, e(nx),
\]
with Whittaker functions $K$ normalized as in \cite[p.\,264]{deshouillers1982kloosterman} (so for example, $\rho_{f\ma}(n)$ is larger by a factor of $2\sqrt{|n|}$ here than in \cite{iwaniec2021spectral}).
Throughout this paper, we use the convention that $\sigma_\infty = \begin{psmall} 1 & 0 \\ 0 & 1 \end{psmall}$; with this convention, if $f$ is a Hecke eigenform, then for $(n, q) = 1$ the Fourier coefficients $\rho_{f\infty}(n)$ are proportional to the Hecke eigenvalues of $f$.

%\fix{Do you really need the choice of scaling matrices?? Probably not since you can just use $\infty$ for the level-avg version...}

\iffalse
\begin{proposition}[Kuznetsov formula \cite{kuznetsov1980petersson,deshouillers1982kloosterman}]\label{prop:kuznetsov}
    ...
\end{proposition}
\begin{proof}
    
\end{proof}
\fi 

A key role in the spectral theory of automorphic forms is played by on-average estimates for the Fourier coefficients and spectral parameters of automorphic forms, known as \emph{density theorems} and \emph{large sieve inequalities}. Following Deshouillers--Iwaniec \cite{deshouillers1982kloosterman}, such results can be proven using trace formulae and bounds for sums of Kloosterman sums\footnote{Duality and Rankin--Selberg $L$-functions give rise to an alternative approach to density theorems \cite{lichtman2024density} and large sieve inequalities \cite{duke2000problem,thorner2021unconditional,pascadi2025large}, which has so far been relevant in higher-rank settings.}; this is not circular even if the original goal was to prove bounds for sums of Kloosterman sums via a trace formula, since an application of Cauchy--Schwarz in the middle can prevent circular arguments.

In the exceptional spectrum, the main feature of density theorems and large sieve inequalities is absorbing a factor of $X^{2\theta_f}$ on average over $f$, where $X$ is as large as possible.
The first result that we need is classical (best attributed to Iwaniec \cite[Ch.~11]{iwaniec2021spectral}), and follows by combining a trace formula with the Weil bound for Kloosterman sums in \cref{eq:weil-bound}.

\begin{lemma}[Density theorem] \label{lem:density-thm}
With the notation above, for any $q \in \Z_+$ and $n \in \Z \setminus \{0\}$, one has
\[
    \sum_{\substack{f \\ \lambda_f < \frac{1}{4}}}^{\Gamma_0(q)}
    X^{2\theta_f} |\rho_{f\ma}(n)|^2 
    \ll_\eps 
    (q|n|)^\eps
    \left(1 + (q, n)^{1/2}\mu(\ma)\sqrt{|n|}\right),
\]
for any nonnegative 
\[
    X \ll \max \left(1, \frac{\mu(\ma)^{-2}}{(q,n)|n|} \right).
\]
\end{lemma}

\begin{proof}
This is the same as \cite[(3.28)]{pascadi2026large}, which follows immediately from \cite[Lemma 2.9]{topacogullari2017certain} with $q_0 = 1$ (the Fourier coefficients in loc.~cit.~are normalized to be smaller by a factor $\asymp \sqrt{|n|}$).
The case $\ma = \infty$ also follows from \cite[(11.23)]{iwaniec2021spectral}: assume $X > e^4$ without loss of generality, and take $L = \tfrac{1}{2}\log X$ in loc.~cit.; note that $t_j$ corresponds to $-i\theta_f$ in our notation, so $(\sin(t_j L)/(t_j L))^4 \asymp X^{2\theta_f + o(1)}$.
\end{proof}

\iffalse
\begin{lemma}[Regular-spectrum large sieve \cite{deshouillers1982kloosterman}]
\label{lem:reg-large-sieve}
\end{lemma}

\begin{proof}
    
\end{proof}
\fi

While \cref{lem:density-thm} shows that $X^{2\theta_f} |\rho_{f\ma}(n)|^2$ is small on average over $f$, many applications require exploiting cancellation in an additional sum over the $n$-variable; following Deshouillers--Iwaniec \cite{deshouillers1982kloosterman}, this can be done through a large sieve inequality. For now, we state a preliminary bound towards exceptional-spectrum large sieve inequalities, which will be used in the proof of \cref{thm:exc-ls-level-avg}.

\begin{lemma}[Towards the exceptional-spectrum large sieve] 
\label{lem:to-exc-large-sieve}
Let $X, N \gg 1$, $(\alpha_n)$ be a complex sequence supported on $n \asymp N$, and $\Phi : \R \to \C$ be a nonnegative smooth function with $\Phi(x)$ supported in $x \asymp 1$, $\Phi^{(j)}(x) \ll_j 1$ for all $j \ge 0$, and $\int \Phi(x) dx \gg 1$. Then with the notation above, one has
\[
\begin{aligned}
    \sum_{\substack{f \\ \lambda_f < \frac{1}{4}}}^{\Gamma_0(q)}
    X^{2\theta_f} \left\vert \sum_n \alpha_n\, \rho_{f\infty}(n)\right\vert^2 
    &\ll 
    \left\vert \sum_{c \in \Z_+} 
    \frac{1}{cq}
    \sum_{m, n} \alpha_m \bar\alpha_n S(m, n; cq)\, \Phi\left(\frac{\sqrt{mn}}{cq}X\right)
    \right\vert 
    \\
    &+ 
    \left(1 + \frac{N^{1+o(1)}}{q}\right) \|\alpha\|^2.
\end{aligned}
\]
\end{lemma}

\begin{proof}
This follows identically as in \cite[Corollary 3.10]{pascadi2026large} (except that here $m, n$ may have support $\asymp N$ rather than $\sim N$), using $\ma = \infty$, $c \gets cq$, and the explicit description of Kloosterman sums in \cite[(3.13)]{pascadi2026large}. The result is essentially present in \cite[p.~271]{deshouillers1982kloosterman}; it follows from an application of the Kuznetsov trace formula, bounds for Bessel functions, and the regular-spectrum large sieve inequalities of Deshouillers--Iwaniec \cite[Theorem 2]{deshouillers1982kloosterman}.
\end{proof}

\section{From Kloosterman fractions to Kloosterman sums} \label{sec:kloosterman-fractions-to-sums}

This section requires no knowledge of automorphic forms. On a broad level, the computations here are similar to previous works \cite{duke1997bilinear,bettin2018trilinear}, but on a technical level, there are several key differences (summarized in \cref{subsec:technical}) which make our final results possible.

\subsection{Setup} \label{subsec:setup}
Let $M, N, A \ge \tfrac{1}{2}$, $\vartheta \in \Z \setminus \{0\}$, and $(\alpha_m)_{m \sim M}, (\beta_n)_{n \sim N}, (\nu_a)_{a \sim A}$ be complex sequences; we may extend these with zeroes outside of the ranges $m \sim M$, $n \sim N$, $a \sim A$. 
We aim to bound the sum 
\begin{equation} \label{eq:S1}
    \msS_1 := \mathop{\sum\sum}_{\substack{m \sim M, n \sim N \\ (m, n) = 1 }} \sum_{a \sim A} \alpha_m \beta_n  \nu_a\, 
    e\left(\vartheta a \frac{\bar{m}}{n} \right),
\end{equation}
under the following simplifying assumptions (which will be removed later in the proofs of \cref{thm:unbalanced-spectral-methods,cor:main-weil-bound-unbalanced}):
\begin{itemize}
    \item[$(i)$.] One has $|\vartheta|, A \ll (MN)^{O(1)}$.
    \item[$(ii)$.] $(\beta_n)$ is supported on integers $n$ with $(n, \vartheta) = 1$. 
\end{itemize}

%\fix{You might also get (almost for free) something good if $\beta_n$ are supported on multiples of a given $g_0$; but you can check that at the very end since only the first prop will be affected.}

%\fix{It might help to assume that $\nu_a \ll 1$, but perhaps it's best to include the $\ell^\infty$ norm in your bounds whenever it's necessary! That will be totally fine for applications due to the $\vartheta$ aspect.}

In anticipation of inserting an amplifier, we also let, for $k \in \Z \setminus \{0\}$,
\begin{equation} \label{eq:L-defn}
    C_0+2\log(2|\vartheta|M) \le L \ll (MN)^{O(1)}, \qquad \mL := \{\ell \sim L : \ell \text{ prime}\}, \qquad \mL_k := \{\ell \in \mL : (\ell, k) = 1\}.
\end{equation}
Here $C_0$ is a sufficiently large absolute constant. Note that if $2\log |k| \le L$, one has 
\begin{equation} \label{eq:L-set-bound}
    \frac{L}{\log L} \ll |\mL_k| \le |\mL| \asymp \frac{L}{\log L}.
\end{equation}
We also denote, throughout this section,
\begin{equation} \label{eq:E-defn}
    E := 1000 \max\left(1, \frac{N}{M}\right) \asymp \frac{M+N}{M},
    \qquad
    \mE := \{e \sim E : e \text{ does not have a prime factor $\ell \sim L$}\}.
\end{equation}
Note that
\begin{equation} \label{eq:E-set-bound}
    |\mE| \ge E-1 - \sum_{\text{prime }\ell \sim L} \left\lfloor\frac{2E}{\ell}\right\rfloor \ge 
    E-1-O\left(\frac{E}{\log L}\right) \asymp E.
\end{equation}

We begin with a preparatory lemma, which will help eliminate some bad ranges in our final proofs.

\begin{lemma}\label{lem:large-sieve-large-powers}
Let $\msS_1$ be as in \cref{eq:S1}, let $C > 0$, and suppose that $(\beta_n)$ is supported on integers $n$ with 
\[
    \frac{n}{\rad(n)} \ge C.
\]
Then, one has %\fix{You might not need the first bound}
\[
    \msS_1^2 \ll N^{o(1)} \|\alpha\|^2 \|\beta\|^2 \|\nu\|^2 (M+N) \min\left(N+\frac{A}{N}, \frac{N+A}{C}\right).
\]
\end{lemma}
\begin{proof}
We apply Cauchy--Schwarz in the $m, n$ variables to obtain
\[
    \msS_1^2 \ll \|\alpha\|^2 \|\beta\|^2 \sum_{\substack{n \sim N \\ n \ge C\rad(n) \\ (n, \vartheta) = 1}}\ \sum_{\substack{m \sim M \\ (m, n) = 1}} \left\vert
    \sum_{a \sim A} \nu_a\, 
    e\left(\vartheta a \frac{\bar{m}}{n} \right)
    \right\vert^2.
\]
Since the inner sum only depends on the residue of $m$ modulo $n$, we can complete the sum in $x = \vartheta\bar{m}$ modulo $n$ (with some overlap if $M > N$) to obtain
\[
    \msS_1^2 \ll \|\alpha\|^2 \|\beta\|^2 \left(1 + \frac{M}{N}\right)  \sum_{\substack{n \sim N \\ n \ge C\rad(n)}}
    \sum_{x \pmod{n}^\times} \left\vert
    \sum_{a \sim A} \nu_a\, 
    e\left(a \frac{x}{n} \right)
    \right\vert^2.
\]
By \cref{lem:gl1-large-sieve}, we conclude that
\[
    \msS_1^2 \ll \|\alpha\|^2 \|\beta\|^2 \|\nu\|^2 \left(1 + \frac{M}{N}\right) \min\Big(N^2+A, \sum_{\substack{n \sim N \\ n \ge C\rad(n)}}
    (N+A) \Big)
\]
To complete the proof, we bound
\[
    \sum_{\substack{n \sim N \\ \rad(n) \le \frac{2N}{C}}} 1
    \le 
    \sum_{r \le \frac{2N}{C}} \sum_{\substack{n \sim N \\ n \mid r^\infty}} 1
    \ll 
    \frac{N^{1+o(1)}}{C},
\]
by \cref{lem:spec-prime-factors}.
\end{proof}

\subsection{The first Cauchy--Schwarz} \label{subsec:first-CS}
Here we perform the amplification argument.

\begin{proposition}\label{prop:first-CS}
Let $\msS_1$ be as in \cref{eq:S1}, and let $C \ge \max\{\tfrac{n}{\rad(n)} : \beta_n \neq 0\}$. Then recalling \cref{eq:L-defn,eq:E-defn,eq:smooth-maj}, one has
\[
    \msS_1^2 \ll (MN)^{o(1)} \left(\frac{M(M+N)}{L}\left(1 + \frac{A^2}{N} + \frac{\sqrt{NL}A}{M}\right)\|\alpha\|^2 \|\beta\|^2\|\nu\|^2 + \frac{\|\alpha\|^2 M}{L^2} \sup_{\substack{G \le 30NL/M \\ G_2 \le C}} |\msS_2| \right),
\]
where
\begin{equation} \label{eq:S2}
\begin{aligned}
    \msS_2 := \sum_{e \in \mE} \sum_{\substack{g \sim G \\  g_2 \sim G_2 \\ g_2 \mid g^\infty }} \sum_{\substack{u_1 \sim N/g \\ u_2 \sim N/(gg_2) \\ (u_1, g_2u_2) = (u_2, g) = 1}} \beta_{gu_1} \bar\beta_{gg_2u_2}
    &\sum_{\substack{\ell_1, \ell_2 \in \mL_{\vartheta g} \\ (u_1\ell_1, u_2\ell_2) = 1 \\ a_1, a_2 \sim A \\ a_1\ell_1 \neq a_2\ell_2}} \nu_{a_1} \bar\nu_{a_2}
    \\
    &\times
    \sum_{\substack{m \in \Z,\,(m, g) = 1 \\ em \mid g(u_1\ell_1 - g_2u_2\ell_2)}} \Phi\left(\frac{m}{M}\right)
    e\left(\vartheta a_1 \frac{\bar{m}}{gu_1} - \vartheta a_2 \frac{\bar{m}}{gg_2u_2}\right).
\end{aligned}
\end{equation}
\end{proposition}

\begin{proof}
We amplify in a different way than in \cite{duke1997bilinear,bettin2018trilinear}, using residue classes $r \pmod{em}$ rather than Dirichlet characters modulo $m$, where $e$ is a suitably-chosen integer; this is partly inspired by Maynard \cite{maynard2025primes}. Recall the set $\mE$ from \cref{eq:E-defn,eq:E-set-bound}. 
Before applying Cauchy--Schwarz, we note that for any $e \in \mE$ and any $\ell \in \mL_{\vartheta m}$, one has
\[
\begin{aligned}
    \sum_{\substack{n \sim N \\ (n, m) = 1 }} \beta_n \sum_{a \sim A} \nu_a\,
    e\left(\vartheta a \frac{\bar{m}}{n} \right)
    &=
    \sum_{r \pmod{em}}
    \sum_{\substack{n \sim N \\ (n, m) = 1 \\ n\ell \equiv r \pmod{em}}} \beta_n \sum_{a \sim A} \nu_a\,
    e\left(\vartheta a \frac{\bar{m}}{n} \right)
    \\
    &\ll 
    \frac{1}{E}
    \sum_{e \in \mE}
    \Bigg\vert
    \sum_{r \pmod{em}} \frac{1}{|\mL_{\vartheta m}|} 
    \sum_{\ell \in \mL_{\vartheta m}}
    \sum_{\substack{n \sim N \\ (n, m) = 1 \\ n\ell \equiv r \pmod{em}}} \beta_n \sum_{a \sim A} \nu_a\,
    e\left(\vartheta a \frac{\bar{m}}{n} \right)
    \Bigg\vert.
\end{aligned}
\]
Then, by Cauchy--Schwarz in $m, e, r$, using \cref{eq:L-set-bound} and inserting a smooth majorant in $m$ as in \cref{eq:smooth-maj}, we get
\begin{equation} \label{eq:first-CS-step}
\begin{aligned} 
    \msS_1^2 
    &\ll 
    \left( 
    \sum_{m \sim M} |\alpha_m| \frac{1}{E} \sum_{e \in \mE}
    \sum_{r \pmod{em}} \frac{1}{|\mL_{\vartheta m}|} 
    \Bigg\vert
    \sum_{\ell \in \mL_{\vartheta m}}
    \sum_{\substack{n \sim N \\ (n, m) = 1 \\ n\ell \equiv r \pmod{em}}} \beta_n \sum_{a \sim A} \nu_a\,
    e\left(\vartheta a \frac{\bar{m}}{n} \right)
    \Bigg\vert \right)^2
    \\ 
    &\ll 
    (MN)^{o(1)} \frac{1}{E^2L^2} \|\alpha\|^2 E^2 M
    \sum_{e \in \mE} \sum_{m \in \Z} \Phi\left(\frac{m}{M}\right) \sum_{r \pmod{em}} 
    \Bigg\vert \sum_{\substack{n \sim N \\ \ell \in \mL_\vartheta \\ (m, n\ell) = 1 \\ n\ell \equiv r \pmod{em}}} \beta_n \sum_{a \sim A} \nu_a\,
    e\left(\vartheta a \frac{\bar{m}}{n} \right) \Bigg\vert^2.
\end{aligned}
\end{equation}
By expanding the square and executing the summation over $r$, it follows that
\begin{equation} \label{eq:S1-dec}
    \msS_1^2
    \ll
    (MN)^{o(1)}
    \frac{\|\alpha\|^2 M}{L^2} \msS_1',
\end{equation}
where
\[
\begin{aligned}
    \msS_1' :=
    \sum_{e \in \mE}
    \sum_{n_1, n_2 \sim N} \beta_{n_1} \bar\beta_{n_2}
    \sum_{\substack{\ell_1, \ell_2 \in \mL_\vartheta \\ a_1, a_2 \sim A}} \nu_{a_1} \bar\nu_{a_2}
    \sum_{\substack{m \in \Z \\ (m, n_1n_2\ell_1\ell_2) = 1  \\ em \mid n_1\ell_1 - n_2\ell_2 \\ }} \Phi\left(\frac{m}{M}\right)
    e\left(\vartheta a_1 \frac{\bar{m}}{n_1} - \vartheta a_2 \frac{\bar{m}}{n_2}\right).
\end{aligned}
\]
We now split
\begin{equation} \label{eq:S1-prime-dec}
    \msS_1' = \msD_1 + \msE_1 + \msF_1 + \msS_1'',
\end{equation}
where: 
\begin{itemize} 
\item[$(i)$.] $\msD_1$ contains all the terms with $a_1\ell_1 = a_2\ell_2$ and $n_1\ell_1 = n_2\ell_2$,
\item[$(ii)$.] $\msE_1$ contains all the terms with $a_1\ell_1 = a_2\ell_2$ and $n_1\ell_1 \neq n_2\ell_2$,
\item[$(iii)$.] $\msF_1$ contains all the terms with $a_1\ell_1 \neq a_2\ell_2$, and $n_1\ell_1 = n_2\ell_2$,
\item[$(iv)$.] $\msS_1''$ contains all the terms with $a_1\ell_1 \neq a_2\ell_2$ and $n_1\ell_1 \neq n_2\ell_2$. 
\end{itemize}

First, we bound $\msD_1$ by summing trivially over $m, e$, Cauchy--Schwarz and the divisor bound:
\[
\begin{aligned}
    \msD_1 \ll 
    ME
    \sum_{a^* \le 4AL}
    \sum_{n^* \le 4NL}
    \left( \sum_{\substack{n \sim N \\ \ell \sim L \\ a \sim A}} \one_{n\ell = n^*} \one_{a\ell = a^*} |\beta_n \nu_a| \right)^2
    &\ll 
    ME (NLA)^{o(1)}
    \sum_{\substack{n \sim N \\ \ell \sim L \\ a \sim A}}
    \left\vert \beta_n \nu_a \right\vert^2
    \\
    &\ll
    (MN)^{o(1)} MEL \|\beta\|^2 \|\nu\|^2.
\end{aligned}
\]
Secondly, we bound $\msE_1$ by the divisor bound for $em \mid n_1\ell_1 - n_2\ell_2$, Cauchy--Schwarz, and the divisor bound again:
\[
\begin{aligned}
    \msE_1 \ll (NL)^{o(1)} \sum_{n_1, n_2 \sim N} |\beta_{n_1} \beta_{n_2}| \sum_{a^* \le 4AL} \left( \sum_{\substack{\ell \sim L \\ a \sim A}} \one_{a\ell = a^*} |\nu_a| \right)^2
    &\ll 
    (NL)^{o(1)} 
    N \|\beta\|^2 (LA)^{o(1)} \sum_{\substack{\ell \sim L \\ a \sim A}} |\nu_a|^2
    \\
    &\ll
    (MN)^{o(1)} NL \|\beta\|^2 \|\nu\|^2.
\end{aligned}
\]
Thirdly, we estimate $\msF_1$ via \cref{cor:incomplete-kl-bound}:
\[
\begin{aligned}
    \msF_1 &\ll E
    \sum_{n^* \le 4NL} 
    \sum_{\substack{n_1, n_2 \sim N \\ \ell_1, \ell_2 \in \mL_\vartheta \\ a_1, a_2 \sim A \\ n_1\ell_1 = n_2\ell_2 = n^* \\ a_1\ell_1 \neq a_2\ell_2}} 
    |\beta_{n_1} \beta_{n_2}
    \nu_{a_1} \nu_{a_2}|
    \left\vert 
    \sum_{\substack{m \in \Z \\ (m, n^*) = 1}} \Phi\left(\frac{m}{M}\right)
    e\left(\vartheta (a_1\ell_1 - a_2\ell_2) \frac{\bar{m}}{n^*}\right)
    \right\vert
    \\
    &\ll 
    (MN)^{o(1)} E \sum_{n^* \le 4NL} 
    \sum_{\substack{n_1, n_2 \sim N \\ \ell_1, \ell_2 \in \mL_\vartheta \\ a_1, a_2 \sim A \\ n_1\ell_1 = n_2\ell_2 = n^* \\ a_1\ell_1 - a_2\ell_2 \neq 0}} 
    |\beta_{n_1} \beta_{n_2}
    \nu_{a_1} \nu_{a_2}|
    \left(\frac{M}{NL}(\vartheta(a_1\ell_1-a_2\ell_2),n^*) + \sqrt{NL}\right).
\end{aligned}
\]
Since $\mL_\vartheta$ only contains primes that do not divide $\vartheta$, and since $(\beta_n)$ were assumed (at the beginning of \cref{sec:kloosterman-fractions-to-sums}) to be supported on integers that are coprime with $\vartheta$, we have $(\vartheta(a_1\ell_1-a_2\ell_2),n^*) = (a_1\ell_1 - a_2\ell_2, n^*)$; the latter is $O(AL)$ since $a_1\ell_1 \neq a_2\ell_2$. Using this, Cauchy--Schwarz, and the divisor bound, we get
\[
\begin{aligned}
    \msF_1 &\ll 
    (MN)^{o(1)} E \left(\frac{MA}{N} + \sqrt{NL}\right) \sum_{n^* \le 4NL} A \|\nu\|^2 \left\vert \sum_{\substack{n \sim N \\ \ell \in \mL_\vartheta}} \one_{n\ell = n^*} |\beta_n| \right\vert^2
    \\ 
    &\ll 
    (MN)^{o(1)}
    E A \|\nu\|^2 L \|\beta\|^2  \left(\frac{M}{NL}AL + \sqrt{NL}\right).
\end{aligned}
\]

Since $ME \asymp M + N$ by \cref{eq:E-defn}, we conclude that
\begin{equation} \label{eq:D1-E1-etc-bound}
    \msD_1 + \msE_1 + \msF_1 \ll (MN)^{o(1)}
    (M+N)L\left(1 + \frac{A^2}{N} + \frac{\sqrt{NL}A}{M}\right) \|\beta\|^2 \|\nu\|^2.
\end{equation}

Finally, it remains to develop $\msS_1''$. By writing $g := (n_1, n_2)$ and $u_i := n_i/g$, we reach
\[
    \msS_1'' = 
    \sum_{e \in \mE}
    \sum_{g \le 2N}
    \sum_{\substack{u_1, u_2 \sim N/g \\ (u_1, u_2) = 1}} \beta_{gu_1} \bar\beta_{gu_2}
    \sum_{\substack{\ell_1, \ell_2 \in \mL_\vartheta \\ u_1\ell_1 \neq u_2\ell_2 \\ a_1, a_2 \sim A \\ a_1\ell_1 \neq a_2\ell_2}} \nu_{a_1} \bar\nu_{a_2}
    \sum_{\substack{m \in \Z \\ (m, gu_1u_2\ell_1\ell_2) = 1 \\ em \mid gu_1\ell_1 - gu_2\ell_2}} \Phi\left(\frac{m}{M}\right)
    e\left(\vartheta a_1 \frac{\bar{m}}{gu_1} - \vartheta a_2 \frac{\bar{m}}{gu_2}\right).
\]
The (first) payoff of amplification modulo $em$ with $e \in \mE$ as in \cref{eq:E-defn} is that now, since $m \ge M/10$ on the support of $\Phi(m/M)$, we have $em \ge EM/10 \ge 100N$, and hence 
\[
    \frac{|gu_1\ell_1 - gu_2\ell_2|}{em} \le \frac{3NL}{100N} < L.
\]
In particular, when $u_1\ell_1 \neq u_2\ell_2$, the primes $\ell_1, \ell_2 \sim L$ cannot be among the divisors of $|gu_1\ell_1 - gu_2\ell_2|/(em)$. Recalling \cref{eq:L-defn,eq:E-defn}, we have $(em, \ell_1\ell_2) = 1$, so actually $\ell_1$ and $\ell_2$ cannot divide $gu_1\ell_1 - gu_2\ell_2$; this implies that $\ell_1 \neq \ell_2$, $\ell_1 \nmid gu_2$, and $\ell_2 \nmid gu_1$. Therefore,
\[
    \msS_1'' = 
    \sum_{e \in \mE}
    \sum_{g \le 2N}
    \sum_{u_1, u_2 \sim N/g} \beta_{gu_1} \bar\beta_{gu_2}
    \sum_{\substack{\ell_1, \ell_2 \in \mL_{\vartheta g} \\ (u_1\ell_1, u_2\ell_2) = 1 \\ a_1, a_2 \sim A \\ a_1\ell_1 \neq a_2\ell_2}} \nu_{a_1} \bar\nu_{a_2}
    \sum_{\substack{m \in \Z \\ (m, gu_1u_2\ell_1\ell_2) = 1  \\ em \mid gu_1\ell_1 - gu_2\ell_2}} \Phi\left(\frac{m}{M}\right)
    e\left(\vartheta a_1 \frac{\bar{m}}{gu_1} - \vartheta a_2 \frac{\bar{m}}{gu_2}\right).
\]
We also note that the condition $(m, gu_1u_2\ell_1\ell_2) = 1$ can now be replaced with $(m, g) = 1$. Indeed, given that $(m, g) = 1$ and $em \mid gu_1\ell_1 - gu_2\ell_2$, we must have $m \mid u_1\ell_1 - u_2\ell_2$, and then the constraint $(u_1\ell_1, u_2\ell_2) = 1$ forces $(m, u_1\ell_1) = (m, u_2\ell_2) = 1$. Finally, from $m \mid u_1\ell_1 - u_2\ell_2$, we deduce that
\[
    1 \le \frac{|u_1\ell_1-u_2\ell_2|}{m} \le \frac{30NL}{gM},
\]
so $g \le 30NL/M$.
Finally, we denote 
\[
    g_2 := (u_2, g^\infty),
\]
and note that $g_2 \mid \tfrac{gu_2}{\rad(gu_2)}$ (indeed, if a prime power $p^k \mid g_2$, then $p^k \mid u_2$ and $p \mid g$, so $p^{k+1} \mid gu_2$ and thus $p^k \mid \tfrac{gu_2}{\rad(gu_2)}$). Whenever $\beta_{gu_2} \neq 0$, we thus have $g_2 \le C$ by the assumption of \cref{prop:first-CS}.

By the change of variables $u_2 \gets g_2 u_2$, we then obtain
\[
\begin{aligned}
    \msS_1'' =
    \sum_{e \in \mE}
    \sum_{\substack{g \le 30NL/M \\ g_2 \le C \\ g_2 \mid g^\infty }}
    \sum_{\substack{u_1 \sim N/g \\ u_2 \sim N/(gg_2) \\ (u_1, g_2u_2) = (u_2, g) = 1}}
    &\beta_{gu_1} \bar\beta_{gg_2u_2}
    \sum_{\substack{\ell_1, \ell_2 \in \mL_{\vartheta g} \\ (u_1\ell_1, u_2\ell_2) = 1 \\ a_1, a_2 \sim A \\ a_1\ell_1 \neq a_2\ell_2}} \nu_{a_1} \bar\nu_{a_2}
    \\
    &\times \sum_{\substack{m \in \Z,\,(m, g) = 1 \\ em \mid g(u_1\ell_1 - g_2u_2\ell_2)}} \Phi\left(\frac{m}{M}\right)
    e\left(\vartheta a_1 \frac{\bar{m}}{gu_1} - \vartheta a_2 \frac{\bar{m}}{gg_2u_2}\right).
\end{aligned}
\]
Putting $g, g_2$ into dyadic ranges gives
\begin{equation} \label{eq:S1-second-dec}
    \msS_1'' \ll (NL)^{o(1)} \sup_{\substack{G \le 30NL/M \\ G_2 \le C}} |\msS_2|,
\end{equation}
where $\msS_2$ is as in \cref{eq:S2}. Combining \cref{eq:S1-dec,eq:S1-prime-dec,eq:D1-E1-etc-bound,eq:S1-second-dec} then completes our proof. 
\end{proof}

\subsection{Switching divisors and swapping moduli}
The next step is to replace the summation over $m \asymp M$ with a summation over a complementary divisor $d$. Manipulating the resulting exponential phases becomes a fairly convoluted task, especially due to the presence of the GCD factors $g, g_2$.

\begin{proposition} \label{prop:switch-div}
Let $\msS_2$ be as in \cref{eq:S2}. Then one has
\[
    \msS_2 \ll (MN)^{o(1)} \left(\frac{|\vartheta|L^2 A^2}{M} \|\beta\|^2 \|\nu\|^2 + \sup_{D \le 30 \frac{NL}{GM}} |\msS_3|\right),
\]
where
\begin{equation} \label{eq:S3}
\begin{aligned}
    \msS_3 := \sum_{\substack{g \sim G \\  g_2 \sim G_2 \\ g_2 \mid g^\infty }} \sum_{\substack{u_1 \sim N/g \\ u_2 \sim N/(gg_2) \\ (u_1, g_2u_2) = (u_2, g) = 1}} \beta_{gu_1} \bar\beta_{gg_2u_2}
    \sum_{\substack{\ell_1, \ell_2 \in \mL_{\vartheta g} \\ (u_1\ell_1, u_2\ell_2) = 1 \\ a_1, a_2 \sim A \\ a_1\ell_1 \neq a_2\ell_2}} \nu_{a_1} \bar\nu_{a_2}
    \sum_{\substack{|d| \sim D \\ (u_1\ell_1-g_2u_2\ell_2, gd) = |d| \\ m := (u_1\ell_1-g_2u_2\ell_2)/d}} \ \sum_{\substack{e \in \mE \\ e \mid dg}}
    \\
    \times\
     \Phi\left(\frac{u_1\ell_1 - g_2u_2\ell_2}{dM}\right)
    e\left(\vartheta d(a_1\ell_1 - a_2\ell_2)\frac{\bar{u_1\ell_1gg_2}}{u_2\ell_2} + \vartheta(a_1\ell_1-a_2\ell_2)\frac{\bar{mu_2\ell_2}}{gg_2}\right).
\end{aligned}
\end{equation}
\end{proposition}

\begin{proof} 
\iffalse
We start by putting $m$ into residue classes $r \pmod{gg_2}$ to obtain
\begin{equation} \label{eq:S2-residues}
\begin{aligned}
    \msS_2 = \sum_{e \in \mE} \sum_{\substack{g \sim G \\  g_2 \sim G_2 \\ g_2 \mid g^\infty  \\ r \pmod{gg_2}^\times}} \sum_{\substack{u_1 \sim N/g \\ u_2 \sim N/(gg_2) \\ (u_1, g_2u_2) = (u_2, g) = 1}} \beta_{gu_1} \bar\beta_{gg_2u_2}
    &\sum_{\substack{\ell_1, \ell_2 \in \mL_{\vartheta g} \\ (u_1\ell_1, u_2\ell_2) = 1 \\ a_1, a_2 \sim A}} \nu_{a_1} \bar\nu_{a_2}
    \\
    &\times
    \sum_{\substack{m \in \Z \\ em \mid g(u_1\ell_1 - g_2u_2\ell_2) \\ m \equiv r \pmod{gg_2}}} \Phi\left(\frac{m}{M}\right)
    e\left(a_1 \frac{\bar{m}}{gu_1} - a_2 \frac{\bar{m}}{gg_2u_2}\right).
\end{aligned}
\end{equation}
\fi
In the inner summation of \cref{eq:S2}, since $(m, g) = 1$, we may denote
\[
    d := \frac{u_1\ell_1 - g_2u_2\ell_2}{m} \in \Z \setminus \{0\},
    \qquad\quad 
    |d| \le 30\frac{NL}{GM}.
\]
Note that $(u_2\ell_1\ell_2, gg_2) = (u_1\ell_1, g_2u_2\ell_2) = (m, gg_2u_1u_2\ell_1\ell_2) = 1$ and that
\[
    \bar{u_1} \equiv \ell_1 \bar{g_2u_2\ell_2} \pmod{m},
    \qquad\quad 
    \bar{m} \equiv d \bar{u_1\ell_1} \pmod{u_2\ell_2}.
\]
Using these and \cref{lem:bezout} three times, the argument of the exponential in \cref{eq:S2} can thus be rewritten as $\vartheta$ times
\begin{equation} \label{eq:argument-manipulation}
\begin{aligned}
    &a_1\frac{\bar{m}}{gu_1} - a_2\frac{\bar{m}}{gg_2u_2}  
    \\ 
    &\equiv
    \left(\frac{a_1}{mgu_1} - a_1 \frac{\bar{gu_1}}{m}\right) - \left(a_2\frac{\bar{mgg_2}}{u_2} +  a_2\frac{\bar{mu_2}}{gg_2}\right)
    \\
    &\equiv 
    \frac{a_1}{mgu_1} - a_1 \ell_1 \frac{\bar{gg_2u_2\ell_2}}{m}
    - da_2\frac{\bar{u_1\ell_1gg_2}}{u_2}  - a_2\frac{\bar{mu_2}}{gg_2}
    \\
    &\equiv 
    \frac{a_1}{mgu_1} - \left(\frac{a_1 \ell_1}{mgg_2u_2\ell_2} - a_1\ell_1 \frac{\bar{m}}{gg_2u_2\ell_2}\right)
    - da_2\frac{\bar{u_1\ell_1gg_2}}{u_2}  - a_2\frac{\bar{mu_2}}{gg_2}
    \\
    &\equiv 
    \frac{a_1}{mgu_1} - \frac{a_1 \ell_1}{mgg_2u_2\ell_2} + \left(a_1\ell_1 \frac{\bar{mgg_2}}{u_2\ell_2}
    +
    a_1\ell_1 \frac{\bar{mu_2\ell_2}}{gg_2}
    \right)
    - da_2\frac{\bar{u_1\ell_1gg_2}}{u_2}  - a_2\frac{\bar{mu_2}}{gg_2}
    \\
    &\equiv 
    \frac{a_1 (g_2u_2\ell_2 - u_1\ell_1)}{mgg_2u_1u_2\ell_2} + 
    d a_1\ell_1 \frac{\bar{gg_2u_1\ell_1}}{u_2\ell_2}
    +
    a_1\ell_1 \frac{\bar{mu_2\ell_2}}{gg_2}
    - da_2\ell_2\frac{\bar{u_1\ell_1gg_2}}{u_2\ell_2}  - a_2\ell_2\frac{\bar{mu_2\ell_2}}{gg_2}
    \\
    &\equiv 
    -
    \frac{da_1}{gg_2u_1u_2\ell_2}
    +
    d (a_1\ell_1 - a_2\ell_2) \frac{\bar{gg_2u_1\ell_1}}{u_2\ell_2}
    +
    (a_1\ell_1 - a_2\ell_2) \frac{\bar{mu_2\ell_2}}{gg_2}
    \quad 
    \pmod{1}.
\end{aligned}
\end{equation}
The first fraction is bounded by
\begin{equation} \label{eq:argument-error}
    \frac{da_1}{gg_2u_1u_2\ell_2}
    \ll 
    \frac{\frac{NL}{GM} A}{N^2 L G^{-1}}
    = 
    \frac{A}{MN}.
\end{equation}
It follows from \cref{eq:S2,eq:argument-manipulation,eq:argument-error} that
\begin{equation} \label{eq:S2-dec}
    \msS_2 = \msS_2' + \msE_2,
\end{equation}
where
\[
\begin{aligned}
    \msS_2' := \sum_{e \in \mE} \sum_{\substack{g \sim G \\  g_2 \sim G_2 \\ g_2 \mid g^\infty }} \sum_{\substack{u_1 \sim N/g \\ u_2 \sim N/(gg_2) \\ (u_1, g_2u_2) = (u_2, g) = 1}} \beta_{gu_1} \bar\beta_{gg_2u_2}
    \sum_{\substack{\ell_1, \ell_2 \in \mL_{\vartheta g} \\ (u_1\ell_1, u_2\ell_2) = 1 \\ a_1, a_2 \sim A \\ a_1\ell_1 \neq a_2\ell_2}} \nu_{a_1} \bar\nu_{a_2}
    \sum_{\substack{m \in \Z,\,(m, g) = 1 \\ em \mid g(u_1\ell_1 - g_2u_2\ell_2) \\ d := (u_1\ell_1-g_2u_2\ell_2)/m}}
    \\
    \times\
     \Phi\left(\frac{m}{M}\right)
    e\left(\vartheta d(a_1\ell_1 - a_2\ell_2)\frac{\bar{u_1\ell_1gg_2}}{u_2\ell_2} + \vartheta(a_1\ell_1-a_2\ell_2)\frac{\bar{mu_2\ell_2}}{gg_2}\right),
\end{aligned}
\]
and, absorbing $e, m$ by the divisor bound,
\begin{equation} \label{eq:E2-bound}
    \begin{aligned} 
    \msE_2 &\ll 
    \sum_{e \in \mE} \sum_{\substack{g \sim G \\  g_2 \sim G_2 \\ g_2 \mid g^\infty }} \sum_{\substack{u_1 \sim N/g \\ u_2 \sim N/(gg_2) \\ (u_1, g_2u_2) = (u_2, g) = 1}} |\beta_{gu_1} \beta_{gg_2u_2}|
    \sum_{\substack{\ell_1, \ell_2 \in \mL_{\vartheta g} \\ (u_1\ell_1, u_2\ell_2) = 1 \\ a_1, a_2 \sim A \\ a_1\ell_1 \neq a_2\ell_2}} |\nu_{a_1} \nu_{a_2}|
    \sum_{\substack{m \in \Z,\,(m, g) = 1 \\ em \mid g(u_1\ell_1 - g_2u_2\ell_2)}} \Phi\left(\frac{m}{M}\right) \frac{|\vartheta|A}{MN}
    \\
    &\ll
    (MN)^{o(1)} \sum_{n_1, n_2 \sim N} |\beta_{n_1} \beta_{n_2}|
    L^2 \sum_{a_1, a_2 \sim A} |\nu_{a_1} \nu_{a_2}|
     \frac{|\vartheta|A}{MN}
    \\
    &
    \ll 
    (MN)^{o(1)} \frac{|\vartheta|L^2 A^2}{M} \|\beta\|^2 \|\nu\|^2.
    \end{aligned}
\end{equation}
We may now transition from the conditions on $m$ to conditions on $d$ as follows:
\[
    \begin{cases} 
    m \in \Z,\ (m, g) = 1 \\ 
    em \mid g(u_1\ell_1-g_2u_2\ell_2).
    \end{cases} 
    \qquad \iff \qquad 
    \begin{cases} 
    (u_1\ell_1 - g_2u_2\ell_2, gd) = |d| \\ 
    e \mid dg.
    \end{cases} 
\]
giving
\[
\begin{aligned}
    \msS_2' = \sum_{\substack{g \sim G \\  g_2 \sim G_2 \\ g_2 \mid g^\infty }} \sum_{\substack{u_1 \sim N/g \\ u_2 \sim N/(gg_2) \\ (u_1, g_2u_2) = (u_2, g) = 1}} \beta_{gu_1} \bar\beta_{gg_2u_2}
    \sum_{\substack{\ell_1, \ell_2 \in \mL_{\vartheta g} \\ (u_1\ell_1, u_2\ell_2) = 1 \\ a_1, a_2 \sim A \\ a_1\ell_1 \neq a_2\ell_2}} \nu_{a_1} \bar\nu_{a_2}
    \sum_{\substack{d \in \Z \setminus \{0\} \\ (u_1\ell_1-g_2u_2\ell_2, gd) = |d| \\ m := (u_1\ell_1-g_2u_2\ell_2)/d}} \ \sum_{\substack{e \in \mE \\ e \mid dg}}
    \\
    \times\
     \Phi\left(\frac{u_1\ell_1 - g_2u_2\ell_2}{dM}\right)
    e\left(\vartheta d (a_1\ell_1 - a_2\ell_2)\frac{\bar{u_1\ell_1gg_2}}{u_2\ell_2} + \vartheta(a_1\ell_1-a_2\ell_2)\frac{\bar{mu_2\ell_2}}{gg_2}\right).
\end{aligned}
\]
Finally, we put $|d|$ in dyadic ranges to obtain
\begin{equation} \label{eq:S2-second-dec}
    \msS_2' \ll (MN)^{o(1)} \sup_{D \le 30\frac{NL}{GM}} |\msS_3|,
\end{equation}
where $\msS_3$ is as in \cref{eq:S3}.
Putting \cref{eq:S2-dec,eq:E2-bound,eq:S2-second-dec} together completes our proof.
\end{proof}

\subsection{Separation of variables}
To aid our later computations, we take a moment to simplify our sums by separating variables in the smooth weight depending on $u_1\ell_1 - g_2u_2\ell_2$, and by removing the exponential phase modulo $gg_2$.

\begin{proposition}\label{prop:sep-vars}
Let $\msS_3$ be as in \cref{eq:S3}. Then one has
\[
    \msS_3 \ll (MN)^{o(1)} \left(\sup_{t, \xi \in \R} \msS_4 + (MN)^{-50} \|\beta\|^2 \|\nu\|^2\right),
\]
with
\begin{equation} \label{eq:S4}
\begin{aligned}
    \msS_4 := \sum_{\substack{g \sim G \\  g_2 \sim G_2 \\ g_2 \mid g^\infty  \\ r \pmod{gg_2}}} \sum_{\substack{u_1 \sim N/g \\ u_2 \sim N/(gg_2) \\ (u_1, g_2u_2) = (u_2, g) = 1}} |\beta_{gu_1} \beta_{gg_2u_2}|
    \Bigg\vert
    \sum_{\substack{\ell_1, \ell_2 \in \mL_{\vartheta g} \\ (u_1\ell_1, u_2\ell_2) = 1 \\ a_1, a_2 \sim A \\ a_1\ell_1 \neq a_2\ell_2}} \nu_{a_1} \bar\nu_{a_2}
    \sum_{\substack{|d| \sim D \\ (u_1\ell_1-g_2u_2\ell_2, gd) = |d| \\ m := (u_1\ell_1-g_2u_2\ell_2)/d \\ (a_1\ell_1 - a_2\ell_2)\bar{m\ell_2} \equiv r \pmod{gg_2}}} 
    \\
    \times\ \gamma_{\ell_1,\ell_2,d,g}\
    e\left(\vartheta d (a_1\ell_1 - a_2\ell_2)\frac{\bar{u_1\ell_1gg_2}}{u_2\ell_2}\right)
    \Bigg\vert,
\end{aligned}
\end{equation}
where
\begin{equation} \label{eq:gamma-coefficients}
    \gamma_{\ell_1,\ell_2,d,g} := \frac{L|d|}{\ell_1 D} e\left(-\xi\frac{\ell_2}{\ell_1} \cdot \frac{NL}{gDM} - t\xi \frac{Ld}{\ell_1D}\right) \sum_{\substack{e \in \mE \\ e \mid dg}} 1 \ll (MN)^{o(1)}.
\end{equation}
\end{proposition}

\begin{proof}
Let $\eta > 0$.
We first separate variables in the smooth weight from \cref{eq:S3} by a double Fourier integral. We use linear changes of variables to write
\begin{equation} \label{eq:fourier-sep-vars-1}
\begin{aligned}
    \Phi\left(\frac{u_1\ell_1 - g_2u_2\ell_2}{dM}\right)
    &=
    \int \hat\Phi(\xi)\, e\left(\xi\frac{u_1\ell_1 - g_2u_2\ell_2}{dM}\right) d\xi
    \\
    &=
    \frac{NL|d|}{gg_2u_2\ell_1D}
    \int \hat\Phi\left(\xi \frac{NLd}{gg_2u_2\ell_1D}\right) e\left(\xi\left(\frac{u_1}{g_2u_2} - \frac{\ell_2}{\ell_1}\right) \frac{NL}{gDM}\right) d\xi,
\end{aligned}
\end{equation}
where $NL|d| \asymp gg_2u_2\ell_1D$ on the support of \cref{eq:S3}. By the Schwartz decay of $\hat\Phi$, we may restrict the last integral to $|\xi| \le (MN)^\eta$ with an error of $O_{\eta,B}((MN)^{-B})$ for any fixed $B > 0$. For these $\xi$, we write
\begin{equation} \label{eq:fourier-sep-vars-2}
\begin{aligned}
    \hat\Phi\left(\xi \frac{NLd}{gg_2u_2\ell_1D}\right)
    &=
    \int 
    \Phi(t)\,
    e\left(-t\xi \frac{NLd}{gg_2u_2\ell_1D}\right) dt
    \\
    &=
    \frac{gg_2u_2}{N}
    \int 
    \Phi\left(t\frac{gg_2u_2}{N}\right)
    e\left(-t\xi \frac{Ld}{\ell_1D}\right) dt.
\end{aligned}
\end{equation}
Since $\Phi$ is supported in $[0.1, 10]$ and $gg_2u_2 \sim N$ in \cref{eq:S3}, the integral on the second line in \cref{eq:fourier-sep-vars-2} is supported on $t \in [0.01, 100]$ for any $\xi$. We plug this into \cref{eq:S3} and apply the triangle inequality in $g, g_2, u_1, u_2, \xi, t$, which allows us to eliminate the factor
\[
    \Phi\left(t\frac{gg_2u_2}{N}\right) 
    e\left(\xi \frac{u_1}{g_2u_2} \cdot \frac{NL}{gDM}\right) 
    \ll 1.
\]
This leads to the bound
\begin{equation} \label{eq:S3-dec}
    |\msS_3| \ll \msS_3' + O_\eta\left((MN)^{-50} \|\beta\|^2 \|\nu\|^2\right),
\end{equation}
where, recalling \cref{eq:gamma-coefficients},
\[
\begin{aligned}
    \msS_3' := \sum_{\substack{g \sim G \\  g_2 \sim G_2 \\ g_2 \mid g^\infty }} \sum_{\substack{u_1 \sim N/g \\ u_2 \sim N/(gg_2) \\ (u_1, g_2u_2) = (u_2, g) = 1}} |\beta_{gu_1} \beta_{gg_2u_2}|
    \int_{0.01}^{100} \int_{-(MN)^\eta}^{(MN)^\eta}
    \Bigg\vert
    \sum_{\substack{\ell_1, \ell_2 \in \mL_{\vartheta g} \\ (u_1\ell_1, u_2\ell_2) = 1 \\ a_1, a_2 \sim A \\ a_1\ell_1 \neq a_2\ell_2}} \nu_{a_1} \bar\nu_{a_2}
    \sum_{\substack{|d| \sim D \\ (u_1\ell_1-g_2u_2\ell_2, gd) = |d| \\ m := (u_1\ell_1-g_2u_2\ell_2)/d}}
    \\
    \times\
    \gamma_{\ell_1,\ell_2,d,g}\
    e\left(\vartheta d (a_1\ell_1 - a_2\ell_2)\frac{\bar{u_1\ell_1gg_2}}{u_2\ell_2} + \vartheta(a_1\ell_1-a_2\ell_2)\frac{\bar{mu_2\ell_2}}{gg_2}\right)
    \Bigg\vert d\xi dt.
\end{aligned}
\]
We use an $L^\infty$-bound in $t,\xi$ to get
\begin{equation} \label{eq:S3-prime-dec}
    \msS_3' \ll (MN)^\eta \sup_{t, \xi \in \R} \msS_3'',
\end{equation}
where
\[
\begin{aligned}
    \msS_3'' := \sum_{\substack{g \sim G \\  g_2 \sim G_2 \\ g_2 \mid g^\infty }} \sum_{\substack{u_1 \sim N/g \\ u_2 \sim N/(gg_2) \\ (u_1, g_2u_2) = (u_2, g) = 1}} |\beta_{gu_1} \beta_{gg_2u_2}|
    \Bigg\vert
    \sum_{\substack{\ell_1, \ell_2 \in \mL_{\vartheta g} \\ (u_1\ell_1, u_2\ell_2) = 1 \\ a_1, a_2 \sim A \\ a_1\ell_1 \neq a_2\ell_2}} \nu_{a_1} \bar\nu_{a_2}
    \sum_{\substack{|d| \sim D \\ (u_1\ell_1-g_2u_2\ell_2, gd) = |d| \\ m := (u_1\ell_1-g_2u_2\ell_2)/d}} 
    \\
    \times\ \gamma_{\ell_1,\ell_2,d,g}\
    e\left(\vartheta d (a_1\ell_1 - a_2\ell_2)\frac{\bar{u_1\ell_1gg_2}}{u_2\ell_2} + \vartheta(a_1\ell_1-a_2\ell_2)\frac{\bar{mu_2\ell_2}}{gg_2}\right)
    \Bigg\vert.
\end{aligned}
\]
We also aim to eliminate the second term in the exponential phase, which has a small denominator $gg_2$. To this end, we apply the triangle inequality in the residue
\[
    r \equiv (a_1\ell_1 - a_2\ell_2) \bar{m\ell_2} \pmod{gg_2},
\]
keeping out $e(\tfrac{\vartheta r\bar{u_2}}{gg_2})$, to reach
\begin{equation} \label{eq:S3-second-dec}
    \msS_3'' \le \msS_4,
\end{equation}
where $\msS_4$ is as in \cref{eq:S4}.
Putting \cref{eq:S3-dec,eq:S3-prime-dec,eq:S3-second-dec} together and recalling that $\eta > 0$ was arbitrary completes our proof.
\end{proof}

%\fix{After applying CS, you should actually recover an archimedean condition by integrating over $\xi, t$ with a smooth weight. But you probably won't need this, so it might be best to just take a sup over $\xi, t$ already.}
%\fix{Alternatively: keep the smooth weight which comes with an archimedean condition, perhaps separate in one of the weights after CS, and then Poisson with what's left; the $u_2$ should simplify! See sketch}

\subsection{The second Cauchy--Schwarz}
We now apply Cauchy--Schwarz again, using \cref{lem:box-decomposition} to eliminate certain `cross-diagonal terms' (with $\ell_1' = \ell_2$ or $\ell_1 = \ell_2'$).

\begin{proposition}\label{prop:second-CS}
Let $\msS_4$ be as in \cref{eq:S4}. Then one has
\[
    \msS_4^2 \ll (MN)^{o(1)} \left(N^2(L+A)L \|\beta\|^4 \|\nu\|^4 + \sup_{D_0 \le D} GG_2 \|\beta\|^4 \msS_5\right),
\]
where
\begin{equation} \label{eq:S5}
\begin{aligned}
    \msS_5 := \sum_{\substack{g \sim G \\ (g, \vartheta) = 1 \\ g_0g_1 = g \\ g_2 \sim G_2 \\ g_2 \mid g^\infty }} 
    \sum_{\substack{\ell_1, \ell_1', \ell_2, \ell_2' \in \mL_{\vartheta g} \\ d_0 \sim D_0 \\ |\delta|, |\delta'| \sim D/d_0 \\ (\delta, \delta') = 1 \\ (\ell_1\ell_1',\ell_2\ell_2') = 1 \\ (d_0\delta\delta', g_2\ell_1\ell_1'\ell_2\ell_2') = 1}}
    \Bigg\vert
    \sum_{\substack{a_1, a_1', a_2, a_2' \sim A \\ a_1\ell_1 \neq a_2\ell_2,\ a_1'\ell_1' \neq a_2' \ell_2' \\ \Delta_1 := \frac{1}{g_1g_2} (\delta\ell_1'\ell_2'(a_1\ell_1 - a_2\ell_2) - \delta'\ell_1\ell_2(a_1'\ell_1'-a_2'\ell_2')) \in \Z \setminus \{0\} \\ \Delta_2 := \frac{1}{g_1}(\delta{\ell_2'}^2(a_1\ell_1-a_2\ell_2) - \delta'\ell_2^2(a_1'\ell_1' - a_2'\ell_2')) \in \Z \\ (\Delta_1, g_0) = 1 }} \nu_{a_1} \bar \nu_{a_1'} \bar\nu_{a_2} \nu_{a_2'}
    \\
    \times
    \sum_{\substack{u_1, u_2 \in \Z \\ (u_1, g_2) = (u_2, g) = 1 \\ (u_1 u_2, d_0\delta\delta') = (u_1\ell_1\ell_1', u_2\ell_2\ell_2') = 1 \\ (u_1\ell_1-g_2u_2\ell_2,gd_0\delta) = |d_0\delta| \\ 
    (u_1\ell_1' - g_2u_2\ell_2', gd_0\delta') = |d_0\delta'| \\ u_1\Delta_1 \equiv u_2 \Delta_2
    \pmod{d_0\delta\delta'g_0}}}
    \Phi\left(\frac{u_1 G}{N}\right) \Phi\left(\frac{u_2 GG_2}{N}\right)
    e\left(\vartheta d_0\Delta_1\frac{\bar{u_1\ell_1\ell_1'g_0}}{u_2\ell_2\ell_2'}\right)
    \Bigg\vert.
\end{aligned}
\end{equation} 
\end{proposition}

\begin{proof}
We apply \cref{lem:box-decomposition} with $X = L$ to partition the pairs $(\ell_1; \ell_2)$ from the sum in \cref{eq:S4}. Since there are no terms with $\ell_1 = \ell_2$ in \cref{eq:S4}, the diagonal term in \cref{eq:box-decomposition} can be ignored, and we obtain
\begin{equation} \label{eq:S4-dec}
    \msS_4 \ll L^{o(1)} \max_{1 \le j \le J_L} \msS_4',
\end{equation}
where (recalling that $(\beta_n)$ is supported on integers $n$ with $(n, \vartheta) = 1$)
\[
\begin{aligned}
    \msS_4' := \sum_{\substack{g \sim G \\ (g, \vartheta) = 1 \\ g_2 \sim G_2 \\ g_2 \mid g^\infty  \\ r \pmod{gg_2}}} \sum_{\substack{u_1 \sim N/g \\ u_2 \sim N/(gg_2) \\ (u_1, g_2u_2) = (u_2, g) = 1}} |\beta_{gu_1} \beta_{gg_2u_2}|
    \Bigg\vert
    \sum_{\substack{\ell_1, \ell_2 \in \mL_{\vartheta g} \\ (\ell_1; \ell_2) \in B_{L,j} \\ (u_1\ell_1, u_2\ell_2) = 1 \\ a_1, a_2 \sim A \\ a_1\ell_1 \neq a_2\ell_2}} \nu_{a_1} \bar\nu_{a_2}
    \sum_{\substack{|d| \sim D \\ (u_1\ell_1-g_2u_2\ell_2, gd) = |d| \\ m := (u_1\ell_1-g_2u_2\ell_2)/d \\ (a_1\ell_1 - a_2\ell_2)\bar{m\ell_2} \equiv r \pmod{gg_2}}}
    \\
    \times\
     \gamma_{\ell_1,\ell_2,d,g}\
    e\left(\vartheta d (a_1\ell_1 - a_2\ell_2)\frac{\bar{u_1\ell_1gg_2}}{u_2\ell_2}\right)
    \Bigg\vert.
\end{aligned}
\]
We now apply Cauchy--Schwarz in the variables $g, g_2, r, u_1, u_2$ and insert smooth majorants for the supports of $u_1, u_2$ to obtain
\begin{equation} \label{eq:S4-prime-dec}
    {\msS_4'}^2 \ll 
    \left(\sum_{n_1, n_2 \sim N} |\beta_{n_1} \beta_{n_2}|^2 GG_2 \right) \msS_4''
    = 
    GG_2 \|\beta\|^4 \msS_4'', 
\end{equation}
where
\[
\begin{aligned}
    \msS_4'' := \sum_{\substack{g \sim G \\ (g, \vartheta) = 1 \\  g_2 \sim G_2 \\ g_2 \mid g^\infty  \\ r \pmod{gg_2}}} \sum_{\substack{u_1, u_2 \in \Z \\ (u_1, g_2u_2) = (u_2, g) = 1}} \Phi\left(\frac{u_1 G}{N}\right) \Phi\left(\frac{u_2 GG_2}{N}\right)
    \Bigg\vert
    \sum_{\substack{\ell_1, \ell_2 \in \mL_{\vartheta g} \\ (\ell_1; \ell_2) \in B_{L,j} \\ (u_1\ell_1, u_2\ell_2) = 1 \\ a_1, a_2 \sim A \\ a_1\ell_1 \neq a_2\ell_2}} \nu_{a_1} \bar\nu_{a_2}
    \\
    \times
    \sum_{\substack{|d| \sim D \\ (u_1\ell_1-g_2u_2\ell_2, gd) = |d| \\ m := (u_1\ell_1-g_2u_2\ell_2)/d \\ (a_1\ell_1 - a_2\ell_2)\bar{m\ell_2} \equiv r \pmod{gg_2}}}
     \gamma_{\ell_1,\ell_2,d,g}\
    e\left(\vartheta d (a_1\ell_1 - a_2\ell_2)\frac{\bar{u_1\ell_1gg_2}}{u_2\ell_2}\right)
    \Bigg\vert^2.
\end{aligned}
\]
We expand the square, using $v'$ to denote a copy of each inner variable $v \in \{\ell_1,\ell_2,a_1,a_2,d,m\}$. Note that since $(\ell_1; \ell_2), (\ell_1', \ell_2') \in B_{L,j}$ simultaneously, we must have $\ell_1' \neq \ell_2$ and $\ell_1 \neq \ell_2'$ (due to the disjoint-coordinate-projections property of $B_{L,j}$ from \cref{eq:disjoint-projections}).  We also denote $d_0 := (d, d')$, $\delta := d/d_0$, $\delta' := d'/d_0$, and
\[
\begin{cases}
    \tilde \Delta_1 := \delta\ell_1'\ell_2'(a_1\ell_1 - a_2\ell_2) - \delta'\ell_1\ell_2(a_1'\ell_1'-a_2'\ell_2')
    \\
    \tilde \Delta_2 := \delta{\ell_2'}^2(a_1\ell_1-a_2\ell_2) - \delta'\ell_2^2(a_1'\ell_1' - a_2'\ell_2')
\end{cases}
\ll 
\frac{DAL^3}{d_0}.
\]
We obtain
\[
\begin{aligned}
    \msS_4'' = \sum_{\substack{g \sim G \\ (g, \vartheta) = 1 \\  g_2 \sim G_2 \\ g_2 \mid g^\infty  \\ r \pmod{gg_2}}} \sum_{\substack{u_1, u_2 \in \Z \\ (u_1, g_2u_2) = (u_2, g) = 1}} \Phi\left(\frac{u_1 G}{N}\right) \Phi\left(\frac{u_2 GG_2}{N}\right)
    \sum_{\substack{\ell_1, \ell_2, \ell_1', \ell_2' \in \mL_{\vartheta g} \\ (\ell_1; \ell_2), (\ell_1'; \ell_2') \in B_{L,j} \\ (u_1\ell_1, u_2\ell_2) = (u_1\ell_1', u_2\ell_2') = 1 \\ a_1, a_2, a_1', a_2' \sim A \\ a_1\ell_1 \neq a_2\ell_2,\ a_1'\ell_1' \neq a_2'\ell_2'}} \nu_{a_1} \bar\nu_{a_2} \bar\nu_{a_1'} \nu_{a_2'}
    \\
    \times
    \sum_{\substack{d_0 \le 2D \\ |\delta|, |\delta'| \sim D/d_0,\ (\delta,\delta') = 1 \\ (u_1\ell_1-g_2u_2\ell_2, gd_0\delta) = |d_0\delta| \\ (u_1\ell_1'-g_2u_2\ell_2', gd_0\delta') = |d_0\delta'| \\ m := (u_1\ell_1-g_2u_2\ell_2)/(d_0\delta) \\ m' := (u_1\ell_1'-g_2u_2\ell_2')/(d_0\delta') \\ (a_1\ell_1 - a_2\ell_2)\bar{m\ell_2} \equiv r \pmod{gg_2} \\ (a_1'\ell_1' - a_2'\ell_2')\bar{m'\ell_2'} \equiv r \pmod{gg_2}}}
     \gamma_{\ell_1,\ell_2,d_0\delta,g}\bar\gamma_{\ell_1',\ell_2',d_0\delta',g}\
    e\left(\vartheta d_0 \tilde\Delta_1\frac{\bar{u_1\ell_1\ell_1'gg_2}}{u_2\ell_2\ell_2'}\right).
\end{aligned}
\]

We then evaluate the sum over $r \pmod{gg_2}$, which leaves the congruence
\[
    (a_1\ell_1 - a_2\ell_2)m'\ell_2' \equiv (a_1'\ell_1' - a_2'\ell_2')m\ell_2 \pmod{gg_2},
\]
i.e., after multiplying everything by $d_0\delta\delta'$,
\[
    \delta(a_1\ell_1 - a_2\ell_2)(u_1\ell_1' - g_2u_2\ell_2')\ell_2' \equiv \delta'(a_1'\ell_1' - a_2'\ell_2')(u_1\ell_1 - g_2u_2\ell_2)\ell_2 \pmod{d_0\delta \delta'gg_2},
\]
which can be rephrased as
\[
    u_1\tilde \Delta_1 \equiv g_2u_2\tilde \Delta_2
    \pmod{d_0\delta\delta'gg_2}.
\]
This implies that $g_2 \mid \tilde\Delta_1$. Furthermore, letting $g_1 := (\tfrac{\tilde \Delta_1}{g_2}, g)$, we have $g_1 \mid \tilde\Delta_2$ since $(u_2, g) = 1$. We thus denote
\[
    \tilde\Delta_1 =: g_1g_2\Delta_1, 
    \qquad\qquad 
    \tilde\Delta_2 =: g_1\Delta_2,
    \qquad\qquad 
    g =: g_0 g_1.
\]
Overall, we obtain
\[
\begin{aligned}
    \msS_4'' 
    &= \sum_{\substack{g \sim G \\ (g, \vartheta) = 1 \\ g_0g_1 = g \\ g_2 \sim G_2 \\ g_2 \mid g^\infty }} 
    \sum_{\substack{\ell_1, \ell_1', \ell_2, \ell_2' \in \mL_{\vartheta g} \\ (\ell_1; \ell_2), (\ell_1'; \ell_2') \in B_{L,j} \\ d_0 \le 2D \\ |\delta|, |\delta'| \sim D/d_0 \\ (\delta, \delta') = 1}}
    \gamma_{\ell_1,\ell_2,d_0\delta,g} \bar\gamma_{\ell_1',\ell_2',d_0\delta',g}
    \sum_{\substack{a_1, a_1', a_2, a_2' \sim A \\ a_1\ell_1 \neq a_2\ell_2,\ a_1'\ell_1' \neq a_2' \ell_2' \\ \Delta_1 := \frac{1}{g_1g_2} (\delta\ell_1'\ell_2'(a_1\ell_1 - a_2\ell_2) - \delta'\ell_1\ell_2(a_1'\ell_1'-a_2'\ell_2')) \in \Z \\ \Delta_2 := \frac{1}{g_1}(\delta{\ell_2'}^2(a_1\ell_1-a_2\ell_2) - \delta'\ell_2^2(a_1'\ell_1' - a_2'\ell_2')) \in \Z \\ (\Delta_1, g_0) = 1}} 
    \\
    &\times
    \nu_{a_1} \bar \nu_{a_1'} \bar\nu_{a_2} \nu_{a_2'}
    \sum_{\substack{u_1, u_2 \in \Z \\ (u_1, g_2u_2) = (u_2, g) = 1 \\ (u_1\ell_1, u_2\ell_2) = (u_1\ell_1', u_2\ell_2') = 1 \\ (u_1\ell_1-g_2u_2\ell_2,gd_0\delta) = |d_0\delta| \\ 
    (u_1\ell_1' - g_2u_2\ell_2', gd_0\delta') = |d_0\delta'| \\ 
    u_1\Delta_1 \equiv u_2 \Delta_2
    \pmod{d_0 \delta \delta' g_0}}}
    \Phi\left(\frac{u_1 G}{N}\right) \Phi\left(\frac{u_2 GG_2}{N}\right)
    e\left(\vartheta d_0 \Delta_1\frac{\bar{u_1\ell_1\ell_1'g_0}}{u_2\ell_2\ell_2'}\right).
\end{aligned}
\]
We can now insert some coprimality constraints at no cost in this sum:
\begin{itemize}
    \item[$(i).$] The constraints $(u_1\ell_1, u_2\ell_2) = (u_1\ell_1', u_2\ell_2') = 1$ can be replaced with $(u_1\ell_1\ell_1', u_2\ell_2\ell_2') = 1$. This is possible due to the previous observation that $\ell_1' \neq \ell_2$ and $\ell_1 \neq \ell_2'$ when $(\ell_1; \ell_2), (\ell_1', \ell_2') \in B_{L,j}$ (and to the fact that $\ell_1, \ell_1', \ell_2, \ell_2'$ are primes).
    \item[$(ii).$] We can add the constraints $(g_2, d_0\delta\delta') = (u_1u_2, d_0\delta\delta') = 1$. Indeed, if a prime $p \mid (u_1, d_0\delta)$, then $p \mid u_1\ell_1$ and $p \mid d \mid u_1\ell_1 - g_2u_2\ell_2$, so $p \mid (u_1\ell_1, g_2u_2\ell_2)$. But this GCD equals $1$ since $(u_1\ell_1,g_2) = (u_1\ell_1,u_2\ell_2) = 1$. The constraints $(u_1, \delta') = 1$ and then $(g_2u_2, d_0\delta\delta') = 1$ follow similarly.
    \item[$(iii).$] We can add the constraint $(d_0\delta\delta', \ell_1\ell_1'\ell_2\ell_2') = 1$. Indeed, recalling from \cref{prop:switch-div} that $D \le 30\tfrac{NL}{GM}$, we have $|d_0\delta| \le 60\tfrac{NL}{GM}$. At the same time, the shape of $\gamma_{\ell_1,\ell_2,d_0\delta,g}$ from \cref{eq:gamma-coefficients} forces $d_0\delta g$ to have a divisor $e \in \mE$. So any prime factor of $d_0\delta$ must divide $e$ or $\tfrac{d_0\delta g}{e}$ for some $e \in \mE$. But by \cref{eq:E-defn}, $\ell_1, \ell_1', \ell_2, \ell_2'$ cannot divide $e$, and they cannot divide $\tfrac{d_0\delta g}{e}$ since $|\tfrac{d_0\delta g}{e}| \le 120\tfrac{NL}{ME} \le \frac{120}{1000}L$. This forces $(d_0\delta, \ell_1\ell_1'\ell_2\ell_2') = 1$, and $(\delta', \ell_1\ell_1'\ell_2\ell_2') = 1$ follows symmetrically.
\end{itemize}

We then apply the triangle inequality in $g, g_2, \ell_1, \ell_1', \ell_2, \ell_2', d_0, \delta, \delta'$ (recalling the bound $\gamma_{\ell_1,\ell_2,d_0\delta,g} \ll (MN)^{o(1)}$ from \cref{eq:gamma-coefficients}) to obtain
\[
\begin{aligned}
    \msS_4'' \ll (MN)^{o(1)} \sum_{\substack{g \sim G \\ (g, \vartheta) = 1 \\ g_0g_1 = g \\  g_2 \sim G_2 \\ g_2 \mid g^\infty }} 
    \sum_{\substack{\ell_1, \ell_1', \ell_2, \ell_2' \in \mL_{\vartheta g} \\ (\ell_1; \ell_2), (\ell_1'; \ell_2') \in B_{L,j} \\ d_0 \le 2D \\ |\delta|, |\delta'| \sim D/d_0 \\ (\delta, \delta') = 1 \\ (d_0\delta\delta', g_2\ell_1\ell_1'\ell_2\ell_2') = 1}}
    \Bigg\vert
    \sum_{\substack{a_1, a_1', a_2, a_2' \sim A \\ a_1\ell_1 \neq a_2\ell_2,\ a_1'\ell_1' \neq a_2' \ell_2' \\ \Delta_1 := \frac{1}{g_1g_2} (\delta\ell_1'\ell_2'(a_1\ell_1 - a_2\ell_2) - \delta'\ell_1\ell_2(a_1'\ell_1'-a_2'\ell_2')) \in \Z \\ \Delta_2 := \frac{1}{g_1}(\delta{\ell_2'}^2(a_1\ell_1-a_2\ell_2) - \delta'\ell_2^2(a_1'\ell_1' - a_2'\ell_2')) \in \Z \\ (\Delta_1, g_0) = 1}} \nu_{a_1} \bar\nu_{a_2} \bar\nu_{a_1'} \nu_{a_2'}
    \\
    \times
    \sum_{\substack{u_1, u_2 \in \Z \\ (u_1, g_2) = (u_2, g) = 1 \\ (u_1 u_2, d_0\delta\delta') = (u_1\ell_1\ell_1', u_2\ell_2\ell_2') = 1 \\ (u_1\ell_1-g_2u_2\ell_2,gd_0\delta) = |d_0\delta| \\ 
    (u_1\ell_1' - g_2u_2\ell_2', gd_0\delta') = |d_0\delta'| \\ u_1\Delta_1 \equiv u_2 \Delta_2
    \pmod{d_0\delta\delta'g_0}}}
    \Phi\left(\frac{u_1 G}{N}\right) \Phi\left(\frac{u_2 GG_2}{N}\right)
    e\left(\vartheta d_0\Delta_1\frac{\bar{u_1\ell_1\ell_1'g_0}}{u_2\ell_2\ell_2'}\right)
    \Bigg\vert.
\end{aligned}
\]
The contribution of the terms with $\Delta_1 \neq 0$ to the sum on the right-hand side is bounded by $\sup_{D_0 \le D} \msS_5$ as in \cref{eq:S5}, after dropping the constraint that $(\ell_1; \ell_2), (\ell_1'; \ell_2') \in B_{L,j}$ by positivity and putting $d_0$ in dyadic ranges. In the contribution of the terms with $\Delta_1 = 0$, we use the triangle inequality and drop the congruence modulo $d_0\delta\delta'g_0$. This gives
\begin{equation}\label{eq:S4-second-dec}
    \msS_4'' \ll (MN)^{o(1)}(\msD_4 + \sup_{D_0 \le D} \msS_5), 
\end{equation}
where, grouping back $d = d_0\delta$ and $d' = d_0\delta'$,
\[
\begin{aligned}
    \msD_4 :=
    \sum_{\substack{g \sim G \\ (g, \vartheta) = 1 \\  g_2 \sim G_2 \\ g_2 \mid g^\infty }} 
    \sum_{\substack{\ell_1, \ell_1', \ell_2, \ell_2' \in \mL_{\vartheta g} \\ (\ell_1;\ell_2), (\ell_1';\ell_2') \in B_{L,j} \\ |d|, |d'| \sim D \\  (dd', g_2\ell_1\ell_1'\ell_2\ell_2') = 1}}
    \sum_{\substack{a_1, a_1', a_2, a_2' \sim A \\ a_1\ell_1 \neq a_2\ell_2,\ a_1'\ell_1' \neq a_2' \ell_2' \\ d\ell_1'\ell_2'(a_1\ell_1 - a_2\ell_2) = d'\ell_1\ell_2(a_1'\ell_1'-a_2'\ell_2')}} |\nu_{a_1} \nu_{a_1'} \nu_{a_2} \nu_{a_2'}|
    \sum_{\substack{0.1\frac{N}{G} \le u_1 \le 10\frac{N}{G} \\ 0.1\frac{N}{GG_2} \le u_2 \le 10\frac{N}{GG_2} \\ (u_1, g_2) = (u_2, g) = 1 \\ (u_1 u_2, dd') = 1 \\  (u_1\ell_1\ell_1', u_2\ell_2\ell_2') = 1 \\ (u_1\ell_1-g_2u_2\ell_2,gd) = |d| \\ 
    (u_1\ell_1' - g_2u_2\ell_2', gd') = |d'|}}
    1.
\end{aligned}
\]
Note that the diagonal condition $d\ell_1'\ell_2'(a_1\ell_1 - a_2\ell_2) = d'\ell_1\ell_2(a_1'\ell_1'-a_2'\ell_2')$ can be rewritten as $\tfrac{d(a_1\ell_1 - a_2\ell_2)}{\ell_1\ell_2} = \tfrac{d'(a_1'\ell_1'-a_2'\ell_2')}{\ell_1'\ell_2'}$. Writing this nonzero fraction as $\tfrac{v}{w}$ for some $v \in \Z \setminus \{0\}$ and $w \in \Z_+$ with $(v, w) = 1$ (which forces $v \ll DAL$ and $w \ll L^2$), we see that $d \mid v\ell_1\ell_2$, so $d \mid v$.
We can therefore bound
\[
\begin{aligned}
    \msD_4 
    &\le 
    \sum_{\substack{v \in \Z \setminus \{0\} \\ w \in \Z_+ \\ (v, w) = 1}}
    \sum_{\substack{g \sim G \\ (g, \vartheta) = 1 \\  g_2 \sim G_2 \\ g_2 \mid g^\infty }}
    \sum_{\substack{0.1\frac{N}{G} \le u_1 \le 10\frac{N}{G} \\ 0.1\frac{N}{GG_2} \le u_2 \le 10\frac{N}{GG_2} \\ (u_1, g_2u_2) = (u_2, g) = 1}}
    \Bigg\vert 
    \sum_{\substack{|d| \sim D \\ d \mid v}}
    \sum_{\substack{\ell_1, \ell_2 \in \mL_{\vartheta g} \\ (\ell_1;\ell_2) \in B_{L,j} \\ (u_1\ell_1, u_2\ell_2) = 1 \\ u_1\ell_1 \equiv g_2u_2\ell_2 \pmod{d} \\ a_1,a_2 \sim A \\ d(a_1\ell_1-a_2\ell_2)w = v\ell_1\ell_2}} 
    |\nu_{a_1}\nu_{a_2}| 
    \Bigg\vert^2
    \\
    &\ll 
    \sum_{\substack{v \in \Z \setminus \{0\} \\ w \in \Z_+ \\ (v, w) = 1 }}
    \sum_{\substack{g \sim G \\ (g, \vartheta) = 1 \\  g_2 \sim G_2 \\ g_2 \mid g^\infty }}
    \sum_{\substack{0.1\frac{N}{G} \le u_1 \le 10\frac{N}{G} \\ 0.1\frac{N}{GG_2} \le u_2 \le 10\frac{N}{GG_2} \\ (u_1, g_2u_2) = (u_2, g) = 1}}
    (MN)^{o(1)}
    \sum_{\substack{|d| \sim D \\ d \mid v}}
    \Bigg\vert 
    \sum_{\substack{\ell_1, \ell_2 \in \mL_{\vartheta g} \\ (\ell_1;\ell_2) \in B_{L,j} \\ (u_1\ell_1, u_2\ell_2) = 1 \\ u_1\ell_1 \equiv g_2u_2\ell_2 \pmod{d} \\ a_1,a_2 \sim A \\ d(a_1\ell_1-a_2\ell_2)w = v\ell_1\ell_2}} 
    |\nu_{a_1}\nu_{a_2}| 
    \Bigg\vert^2.
\end{aligned}
\]
Expanding the square and evaluating the sum over $v, w$ now yields a simpler diagonal condition $\ell_1'\ell_2'(a_1\ell_1-a_2\ell_2) = \ell_1\ell_2(a_1'\ell_1'-a_2'\ell_2')$, which can be rewritten as $\ell_1\ell_1'(a_1\ell_2'-a_1'\ell_2) = \ell_2\ell_2'(a_2\ell_1' - a_2'\ell_1)$. The constraint $(\ell_1\ell_1',\ell_2\ell_2') = 1$ (which follows using $(\ell_1; \ell_2), (\ell_1', \ell_2') \in B_{L,j}$ as before) now forces 
\[
    \frac{a_1\ell_2' - a_1'\ell_2}{\ell_2\ell_2'} = \frac{a_2\ell_1' - a_2'\ell_1}{\ell_1\ell_1'} 
    =: h
    \in \Z,
\]
so we obtain
\[
    \msD_4
    \ll 
    (MN)^{o(1)}
    \sum_{\substack{g \sim G \\ (g, \vartheta) = 1 \\  g_2 \sim G_2 \\ g_2 \mid g^\infty }}
    \sum_{\substack{0.1\frac{N}{G} \le u_1 \le 10\frac{N}{G} \\ 0.1\frac{N}{GG_2} \le u_2 \le 10\frac{N}{GG_2} \\ (u_1, g_2u_2) = (u_2, g) = 1}}
    \sum_{|h| \le \frac{2A}{L}}
    \sum_{\substack{\ell_1,\ell_1',\ell_2,\ell_2' \in \mL_{\vartheta g} \\ (u_1\ell_1\ell_1', u_2\ell_2\ell_2') = 1 \\ a_1, a_1', a_2, a_2' \sim A \\ \frac{a_1\ell_2' - a_1'\ell_2}{\ell_2\ell_2'} = \frac{a_2\ell_1'-a_2'\ell_1}{\ell_1\ell_1'} = h}} |\nu_{a_1}\nu_{a_1'}\nu_{a_2}\nu_{a_2'}|
    \sum_{\substack{|d| \sim D \\ d \mid u_1\ell_1 - g_2u_2\ell_2}} 1.
\]
Since $(u_1\ell_1, g_2u_2\ell_2) = 1$, the sum over $d$ is $(MN)^{o(1)}$ by the divisor bound. We then drop the coprimality constraint $(u_1\ell_1\ell_1', u_2\ell_2\ell_2') = 1$ and sum trivially over $g, g_2, u_1, u_2$ (using \cref{lem:spec-prime-factors} for $g_2$) to reach
\begin{equation} \label{eq:D4-bound}
    \msD_4 \ll (MN)^{o(1)} G \frac{N}{G} \frac{N}{GG_2} \sum_{|h| \le \frac{2A}{L}}
    \Bigg\vert \sum_{\substack{\ell, \ell' \in \mL_\vartheta \\ a, a' \sim A \\ a\ell' - a'\ell = \ell\ell' h}} |\nu_a \nu_{a'}| \Bigg\vert^2.
\end{equation}
We split the inner sum into two parts, depending on whether the primes $\ell, \ell'$ are equal or not. If $\ell = \ell'$, we may write $a' = a-\ell h$. If $\ell \neq \ell'$, the equality $a\ell' - a'\ell = \ell\ell' h$ forces $\ell \mid a$ and $\ell' \mid a'$, and we may write $a = b\ell$, $a' = (b-h)\ell'$. Recalling that the sequence $(\nu_a)$ is extended with zeros outside of the range $a \sim A$, we obtain
\begin{equation} \label{eq:D4-bound-split}
\begin{aligned} 
    \sum_{\substack{\ell, \ell' \in \mL_\vartheta \\ a, a' \sim A \\ a\ell' - a'\ell = \ell\ell' h}} |\nu_a \nu_{a'}|
    &\le
    \sum_{\substack{\ell \in \mL_\vartheta \\ a \sim A}} |\nu_a \nu_{a-\ell h}|
    +
    \sum_{\substack{\ell, \ell' \in \mL_\vartheta \\ \ell \neq \ell' \\b \le \frac{2A}{L}}} 
    |\nu_{b\ell} \nu_{(b-h)\ell'}|
    \\
    &\le 
    \Bigg(
    \sum_{\substack{\ell \in \mL_\vartheta \\ a \sim A}} |\nu_{a}|^2
    \Bigg)^{1/2}
    \Bigg(
    \sum_{\substack{\ell \in \mL_\vartheta \\ a \sim A}} |\nu_{a-\ell h}|^2
    \Bigg)^{1/2}
    +
    \Bigg(\sum_{\substack{\ell, \ell' \in \mL_\vartheta \\ b \le \frac{2A}{L}}} 
    |\nu_{b\ell}|^2
    \Bigg)^{1/2}
    \Bigg(\sum_{\substack{\ell, \ell' \in \mL_\vartheta \\ b \le \frac{2A}{L}}} 
    |\nu_{(b-h)\ell'}|^2
    \Bigg)^{1/2}
    \\
    &\ll 
    L \|\nu\|^2 +
    A^{o(1)}
    L \|\nu\|^2
    ,
\end{aligned}
\end{equation}
by Cauchy--Schwarz. It follows that
\[
    \sum_{|h| \le \frac{2A}{L}}
    \Bigg\vert \sum_{\substack{\ell, \ell' \in \mL_\vartheta \\ a, a' \sim A \\ a\ell' - a'\ell = \ell\ell' h}} |\nu_a \nu_{a'}| \Bigg\vert^2
    \ll A^{o(1)} \left(1 + \frac{A}{L}\right) L^2 \|\nu\|^4
\]
and combining this with \cref{eq:D4-bound} gives
\[
    \msD_4 \ll (MN)^{o(1)} \frac{N^2}{GG_2} (L+A)L \|\nu\|^4.
\]
Together with \cref{eq:S4-dec,eq:S4-prime-dec,eq:S4-second-dec}, this completes our proof.
\end{proof}

\subsection{Completion of Kloosterman sums}
Here we finally reach sums of Kloosterman sums; the following counting lemma will be helpful.

\begin{lemma} \label{lem:y-count}
Let $d_0, \delta, \delta', g_0, g_1 \in \Z \setminus \{0\}$, $g = g_0g_1$, $x \in \Z/d_0\delta\Z$, $x' \in \Z/d_0\delta'\Z$, and $\Delta_1, \Delta_2 \in \Z$ such that $(\delta, \delta') = (\Delta_1, g_0) = 1$. Then one has
\[
    \sum_{\substack{z \pmod{d_0\delta\delta'g} \\ z \equiv x \pmod{d_0\delta} \\ z \equiv x' \pmod{d_0\delta'} \\  z\Delta_1 \equiv \Delta_2 \pmod{d_0\delta\delta'g_0}}} 1 \ll |g_1|.
\]
\end{lemma}
\begin{proof}
If the first two congruences are incompatible, the sum is empty. Otherwise, the congruences $z \equiv x \pmod{d_0\delta}$ and $z \equiv x' \pmod{d_0\delta'}$ fix the value of $z$ modulo $[d_0\delta, d_0\delta'] = |d_0\delta\delta'|$. Calling this residue $r \in \{0, 1, \ldots, |d_0\delta\delta'|-1\}$, we can write $z = d_0\delta\delta'k + r$ for some $k \in \Z/g\Z$. We are then left to bound the number of $k$'s such that
\[
    (d_0\delta\delta'k+r)\Delta_1 \equiv \Delta_2 \pmod{d_0\delta\delta'g_0}
    \qquad 
    \iff 
    \qquad 
    d_0\delta\delta'\Delta_1 k \equiv \Delta_2 - r\Delta_1 \pmod{d_0\delta\delta'g_0}.
\]
This forces $d_0\delta\delta' \mid \Delta_2 - r\Delta_1$, say $\Delta_2 - r\Delta_1 =: d_0\delta\delta'\Delta_3$. Then we obtain $\Delta_1 k \equiv \Delta_3 \pmod{g_0}$, which fixes $k$ modulo $g_0$, leaving $|g_1|$ possibilities for $k$ modulo $g$. 
\end{proof}

\begin{proposition}\label{prop:compl-Kl}
Let $\msS_5$ be as in \cref{eq:S5}. Then for any $\eps > 0$, one has
\[
    \msS_5 \ll_\eps (MN)^{o(1)} \left(\frac{N L^2 (L+A) A}{G} \|\nu\|^4 
    +
    \frac{ND_0}{G^2D^2L^2} \msS_6 \right),
\]
where
\begin{equation} \label{eq:S6}
\begin{aligned}
    \msS_6 := 
    \sum_{\substack{g \sim G \\ (g, \vartheta) = 1 \\ g_0g_1 = g \\ g_2 \sim G_2 \\ g_2 \mid g^\infty }} 
    \sum_{\substack{\ell_1, \ell_1', \ell_2, \ell_2' \in \mL_{\vartheta g} \\ d_0 \sim D_0 \\ |\delta|, |\delta'| \sim D/d_0 \\ (\delta, \delta') = 1 \\ (\ell_1\ell_1',\ell_2\ell_2') = 1 \\ (d_0\delta\delta', g_2\ell_1\ell_1'\ell_2\ell_2') = 1}}
    \sum_{\substack{a_1, a_1', a_2, a_2' \sim A \\ a_1\ell_1 \neq a_2\ell_2,\ a_1'\ell_1' \neq a_2' \ell_2' \\ \Delta_1 := \frac{1}{g_1g_2} (\delta\ell_1'\ell_2'(a_1\ell_1 - a_2\ell_2) - \delta'\ell_1\ell_2(a_1'\ell_1'-a_2'\ell_2')) \in \Z \setminus \{0\} \\ \Delta_2 := \frac{1}{g_1}(\delta{\ell_2'}^2(a_1\ell_1-a_2\ell_2) - \delta'\ell_2^2(a_1'\ell_1' - a_2'\ell_2')) \in \Z \\ (\Delta_1, g_0) = 1 }} |\nu_{a_1} \nu_{a_1'} \nu_{a_2} \nu_{a_2'}|
    \\
    \times
    \sum_{\substack{z \pmod{d_0\delta\delta'g} \\ z \equiv g_2\ell_2\bar{\ell_1} \pmod{d_0\delta} \\ z \equiv g_2\ell_2'\bar{\ell_1'} \pmod{d_0\delta'} \\  z\Delta_1 \equiv \Delta_2 \pmod{d_0\delta\delta'g_0}}}
    \Bigg\vert
    \sum_{1 \le |h| \le (MN)^\eps \frac{G}{G_2 D_0} D^2 L^2}
    e\left(h\frac{z \bar{\ell_2\ell_2'}}{d_0|\delta\delta'|g}\right)
    \sum_{\substack{u_2 \in \Z \\ (u_2, d_0\delta\delta'g \ell_1\ell_1') = 1}}
    \frac{1}{u_2}
    \\
    \times\
    \Phi\left(\frac{u_2 GG_2}{N}\right)
    \hat\Phi\left(\frac{hN}{G d_0|\delta\delta'|gu_2\ell_2\ell_2'}\right)
    S\left(\vartheta \Delta_1 \bar{\ell_1\ell_1'|\delta\delta'|gg_0}, h; u_2\ell_2\ell_2' \right)
    \Bigg\vert.
\end{aligned}
\end{equation}
\end{proposition}

\begin{proof}
We rewrite the sum from \cref{eq:S5} as
\begin{equation} \label{eq:S5-dec}
    \begin{aligned}
    \msS_5 = \sum_{\substack{g \sim G \\ (g, \vartheta) = 1 \\ g_0g_1 = g \\ g_2 \sim G_2 \\ g_2 \mid g^\infty }} 
    \sum_{\substack{\ell_1, \ell_1', \ell_2, \ell_2' \in \mL_{\vartheta g} \\ d_0 \sim D_0 \\ |\delta|, |\delta'| \sim D/d_0 \\ (\delta, \delta') = 1 \\ (\ell_1\ell_1',\ell_2\ell_2') = 1 \\ (d_0\delta\delta', g_2\ell_1\ell_1'\ell_2\ell_2') = 1}}
    \Bigg\vert
    \sum_{\substack{a_1, a_1', a_2, a_2' \sim A \\ a_1\ell_1 \neq a_2\ell_2,\ a_1'\ell_1' \neq a_2' \ell_2' \\ \Delta_1 := \frac{1}{g_1g_2} (\delta\ell_1'\ell_2'(a_1\ell_1 - a_2\ell_2) - \delta'\ell_1\ell_2(a_1'\ell_1'-a_2'\ell_2')) \in \Z \setminus \{0\} \\ \Delta_2 := \frac{1}{g_1}(\delta{\ell_2'}^2(a_1\ell_1-a_2\ell_2) - \delta'\ell_2^2(a_1'\ell_1' - a_2'\ell_2')) \in \Z \\ (\Delta_1, g_0) = 1 }} \nu_{a_1} \bar \nu_{a_1'} \bar\nu_{a_2} \nu_{a_2'}
    \\
    \times
    \sum_{\substack{u_2 \in \Z \\ (u_2, d_0\delta\delta'g \ell_1\ell_1') = 1}}
    \Phi\left(\frac{u_2 GG_2}{N}\right)
    \msS_5'
    \Bigg\vert,
\end{aligned}
\end{equation}
where
\[
    \msS_5' := \sum_{\substack{u_1 \in \Z \\ (u_1, d_0\delta\delta'g_2 u_2\ell_2\ell_2') = 1 \\ (u_1\ell_1-g_2u_2\ell_2,gd_0\delta) = |d_0\delta| \\ 
    (u_1\ell_1' - g_2u_2\ell_2', gd_0\delta') = |d_0\delta'| \\  u_1\Delta_1 \equiv u_2 \Delta_2
    \pmod{d_0\delta\delta'g_0}}}
    \Phi\left(\frac{u_1 G}{N}\right)
    e\left(\vartheta d_0 \Delta_1\frac{\bar{u_1\ell_1\ell_1'g_0}}{u_2\ell_2\ell_2'}\right).
\]
We then put $u_1$ into residue classes modulo $d_0\delta\delta'g$ (i.e., modulo $d_0|\delta\delta'|g$); since $(u_2, d_0\delta\delta'g) = 1$, we may write this residue class as $zu_2$. Since $g_2 \mid g^\infty$, the fact that $(u_1, g_2) = 1$ is really a condition on $z$. In fact, under $u_1 \equiv zu_2 \pmod{d_0\delta\delta'g}$, we have 
\[
\begin{cases}
    (u_1, d_0\delta\delta'g_2) = 1 \\ 
    (u_1\ell_1 - g_2u_2\ell_2, gd_0\delta) = |d_0\delta| \\ 
    (u_1\ell_1' - g_2u_2\ell_2', gd_0\delta') = |d_0\delta'| \\
    u_1\Delta_1 \equiv u_2\Delta_2 \pmod{d_0\delta\delta'g_0}
\end{cases}
    \quad 
    \iff 
    \quad
\begin{cases}
    (z, d_0\delta\delta'g_2) = 1 \\ 
    (z\ell_1 - g_2\ell_2, gd_0\delta) = |d_0\delta| \\ 
    (z\ell_1' - g_2\ell_2', gd_0\delta') = |d_0\delta'| \\ 
    z\Delta_1 \equiv \Delta_2 \pmod{d_0\delta\delta'g_0},
\end{cases}
\]
so we obtain
\begin{equation} \label{eq:S5-prime}
    \msS_5' = 
    \sum_{\substack{z \pmod{d_0\delta\delta'g} \\ (z, d_0\delta\delta'g_2) = 1 \\ (z\ell_1 - g_2\ell_2, gd_0\delta) = |d_0\delta| \\ (z\ell_1' - g_2\ell_2', gd_0\delta') = |d_0\delta'| \\ z\Delta_1 \equiv \Delta_2 \pmod{d_0\delta\delta'g_0}}}
    \sum_{\substack{u_1 \in \Z \\ (u_1, u_2\ell_2\ell_2') = 1 \\ u_1 \equiv zu_2 \pmod{d_0\delta\delta'g}}}
    \Phi\left(\frac{u_1 G}{N}\right)
    e\left(\vartheta d_0 \Delta_1\frac{\bar{u_1\ell_1\ell_1'g_0}}{u_2\ell_2\ell_2'}\right).
\end{equation}
The coprimality $(d_0\delta\delta'g, u_2\ell_2\ell_2') = 1$ can be deduced from the constraints in \cref{eq:S5-dec}, so we can apply \cref{lem:kloosterman-completion} with $q = d_0|\delta\delta'|g$, $c = u_2\ell_2\ell_2'$, and $H = (MN)^\eps \frac{G}{G_2 D_0} D^2 L^2$ to get, for any fixed $B > 0$,
\begin{equation} \label{eq:S5-prime-dec}
    \msS_5' = \frac{N}{G d_0|\delta\delta'|gu_2\ell_2\ell_2'} \left(\msE_5 + \msS_5''\right) + O_{\eps,B}\left((MN)^{-B}\right),
\end{equation}
where
\begin{equation} \label{eq:S5-second}
\begin{aligned}
    \msS_5'' = 
    \sum_{\substack{z \pmod{d_0\delta\delta'g} \\ (z, d_0\delta\delta'g_2) = 1 \\ (z\ell_1 - g_2\ell_2, gd_0\delta) = |d_0\delta| \\ (z\ell_1' - g_2\ell_2', gd_0\delta') = |d_0\delta'| \\  z\Delta_1 \equiv \Delta_2 \pmod{d_0\delta\delta'g_0}}}
    \sum_{1 \le |h| \le (MN)^\eps \frac{G}{G_2 D_0} D^2 L^2}
    \hat\Phi\left(\frac{hN}{G d_0|\delta\delta'|gu_2\ell_2\ell_2'}\right)
    e\left(h\frac{zu_2 \bar{u_2\ell_2\ell_2'}}{d_0|\delta\delta'|g}\right)
    \\
    \times
    S\left(\vartheta \Delta_1 \bar{\ell_1\ell_1'|\delta\delta'|gg_0}, h; u_2\ell_2\ell_2' \right),
\end{aligned}
\end{equation}
and
\[
    \msE_5 := 
    \sum_{\substack{z \pmod{d_0\delta\delta'g} \\ (z, d_0\delta\delta'g_2) = 1 \\ (z\ell_1 - g_2\ell_2, gd_0\delta) = |d_0\delta| \\ (z\ell_1' - g_2\ell_2', gd_0\delta') = |d_0\delta'| \\  z\Delta_1 \equiv \Delta_2 \pmod{d_0\delta\delta'g_0}}}
    \hat\Phi(0)\,
    S\left(\vartheta \Delta_1 \bar{\ell_1\ell_1'|\delta\delta'|gg_0}, 0; u_2\ell_2\ell_2' \right)
\]
is the contribution of $h = 0$.
Note that $u_2 \bar{u_2}$ simplifies in the exponential phase modulo $d_0|\delta\delta'|g$ from \cref{eq:S5-second}, which will be crucial for our use of spectral methods in the next section. 

By the Ramanujan bound from \cref{eq:ram-bound}, \cref{lem:y-count}, and the fact that $(\ell_1\ell_1', d_0\delta\delta') = (\ell_2\ell_2', \vartheta) = 1$, we have
\begin{equation} \label{eq:E5-bound}
\begin{aligned}
    \msE_5
    \ll 
    \sum_{\substack{z \pmod{d_0\delta\delta'g} \\ z \equiv g_2\ell_2\bar{\ell_1} \pmod{d_0\delta} \\ z \equiv g_2\ell_2'\bar{\ell_1'} \pmod{d_0\delta'} \\  z\Delta_1 \equiv \Delta_2 \pmod{d_0\delta\delta'g_0}}}
    (\vartheta \Delta_1, u_2\ell_2\ell_2')
    \le 
    g_1 (\vartheta \Delta_1, u_2\ell_2\ell_2')
    \le 
    g_1 (\vartheta\Delta_1, u_2) (\Delta_1, \ell_2\ell_2').
\end{aligned}
\end{equation}
We now plug \cref{eq:S5-second,eq:E5-bound} into \cref{eq:S5-prime-dec} and \cref{eq:S5-dec}; the contribution of $\msS_5''$ leads to the sum $\msS_6$ from \cref{eq:S6} after using the triangle inequality in $a_1, a_1', a_2, a_2', z$. We obtain
\begin{equation} \label{eq:S5-dec-2}
\begin{aligned}
    \msS_5 \ll 
    \frac{N D_0}{G^2D^2L^2} (\msF_5 + \msS_6)
    +
    O_\eps\left((MN)^{-50} \|\nu\|^4\right),
\end{aligned}
\end{equation}
where, after dropping some summation constraints and writing $\tilde\Delta_1 := g_1g_2\Delta_1$, $g_1 = g_1(\Delta_1, g_0) = (g_1\Delta_1, g) \le (\tilde\Delta_1, g)$, the contribution of $\msE_5$ leads to 
\[
\begin{aligned}
    \msF_5 := 
    \sum_{\substack{\ell_1, \ell_1', \ell_2, \ell_2' \in \mL_{\vartheta} \\ d_0 \sim D_0 \\ |\delta|, |\delta'| \sim D/d_0 \\ (\ell_1\ell_1',\ell_2\ell_2') = 1 \\ (d_0\delta\delta', \ell_1\ell_1'\ell_2\ell_2') = 1}}
    \sum_{\substack{a_1, a_1', a_2, a_2' \sim A \\ a_1\ell_1 \neq a_2\ell_2,\ a_1'\ell_1' \neq a_2'\ell_2' \\ \tilde \Delta_1 := \delta\ell_1'\ell_2'(a_1\ell_1 - a_2\ell_2) - \delta'\ell_1\ell_2(a_1'\ell_1'-a_2'\ell_2') \neq 0}} |\nu_{a_1} \nu_{a_1'} \nu_{a_2} \nu_{a_2'}|
    \\
    \times
    \sum_{u_2 \in \Z}
    \frac{1}{u_2}
    \Phi\left(\frac{u_2 GG_2}{N}\right)
    \sum_{\substack{g \sim G \\  g_2 \sim G_2 \\ g_2 \mid g^\infty }} 
    (\tilde\Delta_1, g) (\vartheta\tilde\Delta_1, u_2) (\tilde\Delta_1, \ell_2\ell_2').
\end{aligned}
\]
It remains to bound $\msF_5$, and the only nontrivial task is to absorb the GCD factors, by evaluating the sums in a careful order. Using \cref{lem:spec-prime-factors} and then \cref{lem:gcd-avg}, we obtain
\[
    \sum_{\substack{g \sim G \\ g_2 \sim G_2 \\ g_2 \mid g^\infty}}
    (\tilde\Delta_1, g)
    \ll (MN)^{o(1)}G,
    \qquad\quad 
    \sum_{u_2 \in \Z}
    \frac{1}{u_2}
    \Phi\left(\frac{u_2 GG_2}{N}\right) (\vartheta\tilde\Delta_1, u_2)
    \ll (MN)^{o(1)}.
\]
Using that $(\ell_1\ell_1', \ell_2\ell_2') = (\delta\delta',\ell_2\ell_2') = 1$, we also have
\[
\begin{aligned}
    (\tilde\Delta_1, \ell_2\ell_2')
    =
    (\delta a_1\ell_1\ell_1'\ell_2' - \delta'a_1'\ell_1\ell_1'\ell_2, \ell_2\ell_2')&
    \\
    =
    (\delta a_1\ell_2' - \delta'a_1'\ell_2, \ell_2\ell_2')&
    \le 
    \left(\one_{\ell_2 = \ell_2'} \ell_2 (1 +
    \one_{\ell_2 \mid \delta a_1 - \delta'a_1'} \ell_2)
    +
    \one_{\ell_2 \neq \ell_2'} (a_1,\ell_2)(a_1',\ell_2')
    \right).
\end{aligned}
\]
This leaves us with
\[
\begin{aligned}
    \msF_5 \ll (MN)^{o(1)} G
    &\sum_{\substack{a_1, a_1', a_2, a_2' \sim A }} |\nu_{a_1} \nu_{a_1'} \nu_{a_2} \nu_{a_2'}|
    \\
    &\times
    \sum_{\substack{\ell_1, \ell_1', \ell_2, \ell_2' \in \mL_{\vartheta} \\ d_0 \sim D_0 \\ |\delta|, |\delta'| \sim D/d_0 }} 
    \left(\one_{\ell_2 = \ell_2'} L (1 + \one_{\delta a_1 = \delta'a_1'} L + \one_{\substack{\ell_2 \mid \delta a_1 - \delta'a_1' \neq 0}} L)
    +
    (a_1,\ell_2)(a_1',\ell_2')\right).
\end{aligned}
\]
We treat the contributions to $\msF_5$ of the four terms in the last parenthesis separately:
\begin{itemize}
    \item[$(i).$] For the contribution of $\one_{\ell_2=\ell_2'}L$, we sum trivially over $\ell_1,\ell_1',\ell_2,d_0,\delta,\delta'$ and use Cauchy--Schwarz in $a_1, a_1', a_2, a_2' \sim A$ to get a bound of $\ll (MN)^{o(1)} G A^2 \|\nu\|^4 L^4 D^2 D_0^{-1}$.
    \item[$(ii).$] For the contribution of $\one_{\ell_2 = \ell_2'} \one_{\delta a_1=\delta'a_1'}L^2$, we write $d^* = |\delta a_1| = |\delta'a_1'|$, then use Cauchy--Schwarz and the divisor bound to obtain %\fix{you could improve this with the $d_0 \mid \ell_1\ell_2' - \ell_1'\ell_2$ condition. This could simplify your later bounds, the term wasn't covered by others!}
    \[
    \begin{aligned} 
        G \sum_{\substack{a_1, a_1', a_2, a_2' \sim A }} |\nu_{a_1} \nu_{a_1'} \nu_{a_2} \nu_{a_2'}|
        \sum_{\substack{\ell_1, \ell_1', \ell_2 \in \mL_{\vartheta} \\ d_0 \sim D_0 \\ |\delta|, |\delta'| \sim D/d_0 \\ \delta a_1 = \delta'a_1'}} L^2
        &\ll 
        G L^5 A \|\nu\|^2
        \sum_{\substack{d_0 \sim D_0 \\ d^* \le 4DA}}
        \Bigg( 
        \sum_{\substack{|\delta| \sim D/d_0 \\ a_1 \sim A \\ |\delta a_1| = d^*}}
        |\nu_{a_1}|
        \Bigg)^2
        \\
        &\ll 
        (MN)^{o(1)}
        GL^5 A \|\nu\|^2 
        \sum_{\substack{d_0 \sim D_0 \\ |\delta| \sim D/d_0 \\ a_1 \sim A}} |\nu_{a_1}|^2
        \\
        &\ll 
        (MN)^{o(1)} GL^5 A \|\nu\|^4 D.
    \end{aligned}
    \]
    \item[$(iii).$] For the contribution of $\one_{\ell_2=\ell_2'}\one_{\ell_2 \mid \delta a_1-\delta'a_1' \neq 0} L^2$, we sum over $\ell_2$ using the divisor bound, sum trivially over $\ell_1,\ell_1',d_0,\delta,\delta'$, and use Cauchy--Schwarz in $a_1, a_1', a_2, a_2' \sim A$ to get a bound of $\ll (MN)^{o(1)} G A^2 \|\nu\|^4 L^4 D^2 D_0^{-1}$.
    \item[$(iv).$] For the contribution of $(a_1,\ell_2)(a_1',\ell_2')$, we sum over $\ell_2, \ell_2'$ using \cref{lem:gcd-avg}, sum trivially over $\ell_1,\ell_1',d_0,\delta,\delta'$, and use Cauchy--Schwarz in $a_1, a_1', a_2, a_2' \sim A$ to reach a total bound of $\ll (MN)^{o(1)} G A^2 \|\nu\|^4 L^4 D^2 D_0^{-1}$.
\end{itemize}
Overall, since $D_0 \le D$, we conclude that
\[
    \msF_5 \ll (MN)^{o(1)} G D^2 D_0^{-1} L^4 (L + A) A \|\nu\|^4.
\]
Combining this with \cref{eq:S5-dec-2} completes our proof.
\end{proof}

It will be useful to state the combined result of our work in this section so far. 

\begin{proposition}\label{prop:kl-fr-to-kl-sums}
Let $\msS_1$ be as in \cref{eq:S1} with $\|\nu\| \neq 0$, and let $\eps > 0$. Recalling the notation in \cref{eq:smooth-maj,eq:L-defn}, let $\msS_6$ be as in \cref{eq:S6} (which depends on $\eps$). Then one has
\[
\begin{aligned}
    \msS_1^2 
    &\ll_\eps (MN)^{o(1)} \|\alpha\|^2 \|\beta\|^2 \|\nu\|^2 
    MNA
    \Bigg(
    \frac{(M+N)(N+A^2)}{N^2AL} + \frac{(M+N)}{M\sqrt{NL}} 
    \\ 
    &+
    \frac{|\vartheta|A}{MN} + \frac{1}{L^{3/2}\sqrt{A}}
    +
    \frac{1}{\sqrt{MA}}
    + 
    \frac{1}{\sqrt{N} A L^3\|\nu\|^2} 
    \sup_{\substack{G, G_2 \le 30NL/M \\ D_0 \le D \le 30\frac{NL}{GM}}} 
    \sqrt{\frac{G_2 D_0 \msS_6}{GD^2}}
    \Bigg).
\end{aligned}
\]
\end{proposition}
\begin{proof}
Let $C > 0$ be a parameter to be chosen shortly. First, we split
\[
    \msS_1 := \mathop{\sum\sum}_{\substack{m \sim M, n \sim N \\ (m, n) = 1 }} \alpha_m (\beta_n \one_{\frac{n}{\rad(n)} \ge C})
    \sum_{a \sim A} \nu_a\, 
    e\left(\vartheta a \frac{\bar{m}}{n} \right)
    +
    \tilde\msS_1,
\]
where
\[
    \tilde \msS_1 := \mathop{\sum\sum}_{\substack{m \sim M, n \sim N \\ (m, n) = 1 }} \alpha_m (\beta_n \one_{\frac{n}{\rad(n)} < C})
    \sum_{a \sim A} \nu_a\, 
    e\left(\vartheta a \frac{\bar{m}}{n} \right).
\]
By combining \cref{lem:large-sieve-large-powers} with \cref{prop:first-CS,prop:switch-div,prop:sep-vars,prop:second-CS,prop:compl-Kl}, we obtain
\[
\begin{aligned}
    \msS_1^2 
    \ll_\eps (MN)^{o(1)} \|\alpha\|^2 \|\beta\|^2 \|\nu\|^2 
    &\Bigg(\frac{(M+N)(N+A)}{C} +
    \frac{M(M+N)}{L}\Big(1 + \frac{A^2}{N} + \frac{\sqrt{NL}A}{M}\Big) 
    \\ 
    &+
    \frac{M}{L^2} \Big(\frac{|\vartheta|L^2A^2}{M} + \sqrt{N^2(L+A)L}\Big)
    \\ 
    + 
    \frac{M}{L^2} 
    &\sup_{\substack{G \le 30NL/M \\ G_2 \le C \\ D_0 \le D \le 30\frac{NL}{GM}}} 
    \sqrt{GG_2 \Big(\frac{NL^2(L+A)A}{G} + \frac{N D_0 \msS_6}{G^2D^2L^2\|\nu\|^4} \Big)}
    \Bigg).
\end{aligned}
\]
The first term inside the square-root does not depend on $D$ or $G$, and it attains its maximum when $G_2 = C$; this gives
\begin{equation} \label{eq:S1-bound-preliminary}
\begin{aligned}
    \msS_1^2 
    &\ll_\eps (MN)^{o(1)} \|\alpha\|^2 \|\beta\|^2 \|\nu\|^2 
    \Bigg( \frac{(M+N)(N+A)}{C} +
    \frac{M(M+N)}{L}\Big(1 + \frac{A^2}{N} + \frac{\sqrt{NL}A}{M}\Big) 
    \\ 
    &+
    \frac{M}{L^2} \Big(\frac{|\vartheta|L^2A^2}{M} + \sqrt{N^2(L+A)L} + \sqrt{CNL^2(L+A)A} \Big)
    + 
    \frac{M\sqrt{N}}{L^3\|\nu\|^2} 
    \sup_{\substack{G \le 30NL/M \\ G_2 \le C \\ D_0 \le D \le 30\frac{NL}{GM}}} 
    \sqrt{\frac{G_2 D_0 \msS_6}{GD^2}}
    \Bigg).
\end{aligned}
\end{equation}
We now pick $C := 30NL/M$, so that
\[
    \frac{(M+N)(N+A)}{C} \ll \frac{M(M+N)}{L} \left(1 + \frac{A}{N}\right),
\]
which gets absorbed by the second term in \cref{eq:S1-bound-preliminary}.
Our bound therefore simplifies to
\[
\begin{aligned}
    \msS_1^2 
    &\ll_\eps (MN)^{o(1)} \|\alpha\|^2 \|\beta\|^2 \|\nu\|^2 
    \Bigg(
    \frac{M(M+N)(N+A^2)}{NL} + \frac{(M+N)A\sqrt{N}}{\sqrt{L}} 
    \\ 
    &+
    |\vartheta|A^2 + \frac{MN}{L} + \frac{MN\sqrt{A}}{L^{3/2}}
    +
    N\sqrt{MA}
    +
    \frac{N\sqrt{M}A}{\sqrt{L}}
    + 
    \frac{M\sqrt{N}}{L^3\|\nu\|^2} 
    \sup_{\substack{G, G_2 \le 30NL/M \\ D_0 \le D \le 30\frac{NL}{GM}}} 
    \sqrt{\frac{G_2 D_0 \msS_6}{GD^2}}
    \Bigg).
\end{aligned}
\]
The term $MN/L$ on the second line is absorbed by the term $M(M+N)(N+A^2)/(NL)$ on the first line. Moreover, the term $N\sqrt{M}A/\sqrt{L}$ on the second line is absorbed by the term $(M+N)A\sqrt{N/L}$ on the first line, since $\sqrt{MN} \ll M+N$. The remaining terms give the desired bound.
\end{proof}

\subsection{Using the Weil bound} \label{subsec:weil-bound-argument} We take a moment to note the more direct consequence of the Weil bound for Kloosterman sums. This improves on the work of Bettin--Chandee \cite{bettin2018trilinear} in the $A$-aspect, and is preferable to using spectral methods if $\vartheta$ is large.

\begin{proposition}\label{prop:using-weil-bound}
Let $\eps > 0$ and $\msS_6$ be as in \cref{eq:S6}. Then one has
\[
    \msS_6 \ll_\eps
    (MN)^{3\eps}
    \frac{G^{3/2}}{G_2^{3/2} D_0^2}
    N^{1/2} L^7 D^4 A^2 \|\nu\|^4.
\]
\end{proposition}

\begin{proof}
Plugging the Weil bound from \cref{eq:weil-bound} into \cref{eq:S6}, we obtain
\[
\begin{aligned}
    \msS_6 \ll
    (MN)^{o(1)}
    \sum_{\substack{g \sim G \\ g_0g_1 = g \\ g_2 \sim G_2 \\ g_2 \mid g^\infty }} 
    \sum_{\substack{\ell_1, \ell_1', \ell_2, \ell_2' \in \mL_{\vartheta g} \\ d_0 \sim D_0 \\ |\delta|, |\delta'| \sim D/d_0 \\ (\ell_1\ell_1',\ell_2\ell_2') = 1 \\ (\delta, \delta') = 1 \\ (d_0\delta\delta', g_2\ell_1\ell_1'\ell_2\ell_2') = 1}}
    \sum_{\substack{a_1, a_1', a_2, a_2' \sim A \\ a_1\ell_1 \neq a_2\ell_2,\ a_1'\ell_1' \neq a_2' \ell_2' \\ \Delta_1 := \frac{1}{g_1g_2} (\delta\ell_1'\ell_2'(a_1\ell_1 - a_2\ell_2) - \delta'\ell_1\ell_2(a_1'\ell_1'-a_2'\ell_2')) \in \Z \setminus \{0\} \\ \Delta_2 := \frac{1}{g_1}(\delta{\ell_2'}^2(a_1\ell_1-a_2\ell_2) - \delta'\ell_2^2(a_1'\ell_1' - a_2'\ell_2')) \in \Z \\ (\Delta_1, g_0) = 1 }} |\nu_{a_1} \nu_{a_1'} \nu_{a_2} \nu_{a_2'}|
    \\
    \times
    \sum_{\substack{z \pmod{d_0\delta\delta'g} \\ z \equiv g_2\ell_2\bar{\ell_1} \pmod{d_0\delta} \\ z \equiv g_2\ell_2'\bar{\ell_1'} \pmod{d_0\delta'} \\  z\Delta_1 \equiv \Delta_2 \pmod{d_0\delta\delta'g_0}}}
    \sum_{\substack{1 \le |h| \le (MN)^\eps \frac{G}{G_2 D_0} D^2 L^2 \\ 0.1\frac{N}{GG_2} \le u_2 \le 10 \frac{N}{GG_2}}}
    \frac{1}{u_2}
    (h, u_2\ell_2\ell_2')^{1/2} \sqrt{u_2\ell_2\ell_2'}.
\end{aligned}
\]
We first sum over $h$, using \cref{lem:gcd-avg} to absorb $(h, u_2\ell_2\ell_2')$. We then sum trivially over $u_2$ and use \cref{lem:y-count} to obtain
\[
\begin{aligned}
    \msS_6 \ll_\eps
    (MN)^{2\eps}
    \sum_{\substack{g \sim G \\ g_0g_1 = g \\ g_2 \sim G_2 \\ g_2 \mid g^\infty }} 
    \sum_{\substack{\ell_1, \ell_1', \ell_2, \ell_2' \in \mL_{\vartheta g} \\ d_0 \sim D_0 \\ |\delta|, |\delta'| \sim D/d_0 \\ (\ell_1\ell_1',\ell_2\ell_2') = 1 \\ (d_0\delta\delta', g_2\ell_1\ell_1'\ell_2\ell_2') = 1}}
    \sum_{\substack{a_1, a_1', a_2, a_2' \sim A \\ a_1\ell_1 \neq a_2\ell_2,\ a_1'\ell_1' \neq a_2' \ell_2' \\ \Delta_1 := \frac{1}{g_1g_2} (\delta\ell_1'\ell_2'(a_1\ell_1 - a_2\ell_2) - \delta'\ell_1\ell_2(a_1'\ell_1'-a_2'\ell_2')) \in \Z \setminus \{0\} \\ (\Delta_1, g_0) = 1}} |\nu_{a_1} \nu_{a_1'} \nu_{a_2} \nu_{a_2'}|
    \\
    \times
    g_1 
    \frac{G}{G_2 D_0} D^2 L^2
    \sqrt{\frac{NL^2}{GG_2}}.
\end{aligned}
\]
Once again we write $\tilde\Delta_1 := g_1g_2\Delta_1$, $g_1 = (g_1\Delta_1, g) \le (\tilde\Delta_1, g)$, to bound
\[
\begin{aligned}
    \msS_6 \ll_\eps 
    (MN)^{2\eps}
    \sum_{\substack{\ell_1, \ell_1', \ell_2, \ell_2' \in \mL_{\vartheta} \\ d_0 \sim D_0 \\ |\delta|, |\delta'| \sim D/d_0}} 
    \sum_{\substack{a_1, a_1', a_2, a_2' \sim A \\ \tilde\Delta_1 := \delta\ell_1'\ell_2'(a_1\ell_1 - a_2\ell_2) - \delta'\ell_1\ell_2(a_1'\ell_1'-a_2'\ell_2') \neq 0}} |\nu_{a_1} \nu_{a_1'} \nu_{a_2} \nu_{a_2'}|
    \\
    \times \sum_{\substack{g \sim G \\ g_2 \sim G_2 \\ g_2 \mid g^\infty }} 
    (\tilde\Delta_1, g) 
    \frac{G}{G_2 D_0} D^2 L^2
    \sqrt{\frac{NL^2}{GG_2}}.
\end{aligned}
\]
We bound the sum over $g, g_2$ using \cref{lem:spec-prime-factors} and then \cref{lem:gcd-avg}. We then apply Cauchy--Schwarz in $a_1, a_1', a_2, a_2' \sim A$, and sum trivially over $\ell_1, \ell_1', \ell_2, \ell_2', d_0, \delta, \delta'$ to reach
\[
    \msS_6 \ll_\eps 
    (MN)^{3\eps}
    L^4 \frac{D^2}{D_0} A^2 \|\nu\|^4 G \frac{G}{G_2 D_0} D^2 L^2
    \sqrt{\frac{NL^2}{GG_2}},
\]
which rearranges to the desired bound.
\end{proof}

\begin{theorem}\label{cor:main-weil-bound-unbalanced}
Let $M, N, A \ge \tfrac{1}{2}$ and $(\alpha_m)_{m \sim M}, (\beta_n)_{n \sim N}, (\nu_a)_{a \sim A}$ be complex sequences. For any $\vartheta \in \Z \setminus \{0\}$, one has
\[
\begin{aligned}     
    \mathop{\sum\sum}_{\substack{m \sim M, n \sim N \\ (m, n) = 1 }} \sum_{a \sim A} \alpha_m \beta_n  \nu_a\, 
    e\left(\vartheta a \frac{\bar{m}}{n} \right)
    &\ll (MN)^{o(1)} \|\alpha\| \|\beta\| \|\nu\| 
    \sqrt{MNA}
    \\
    &\times
    \Bigg(
    \frac{|\vartheta|A}{MN} 
    +
    \frac{(M+N)^{3/5}(N+A^2)^{3/5}}{M^{2/5}N^{9/10}A^{3/5}} 
    +
    \frac{N^{3/8}}{M^{1/2}A^{1/4}} 
    + \frac{N^{3/4}}{M}
    \Bigg)^{1/2}.
\end{aligned}
\]
\end{theorem}

\begin{proof}
If $\|\nu\| = 0$, the claim is trivial, so let us assume that $\|\nu\| \neq 0$.
In light of the trivial bound $ \|\alpha\|\|\beta\|\|\nu\| \sqrt{MNA}$, the claim is trivial unless $A, \vartheta \ll (MN)^{O(1)}$, so let us assume this. 

Moreover, by writing $d=(n,\vartheta)$ and $n=dn'$, we obtain
\[
    \msS_1=
    \sum_{d\mid |\vartheta|}
    \mathop{\sum\sum}_{\substack{m\sim M,\ n'\sim N/d \\ (m,n')=1}}
    (\alpha_m\one_{(m,d)=1})
    (\beta_{dn'}\one_{(n',\vartheta/d)=1})
    \sum_{a\sim A}\nu_a\,
    e\left(\frac{\vartheta}{d}a\frac{\bar m}{n'}\right).
\]
The right-hand side of the claimed bound, including the factor $\sqrt{MNA}$ and apart from the norms, is non-decreasing in $N$ and $|\vartheta|$. Thus it suffices to prove the claim when $(\beta_n)$ is supported on integers with $(n,\vartheta)=1$, at a cost of $(MN)^{o(1)}$ by the divisor bound. This brings us to the setting considered in \cref{subsec:setup}. We adopt the notation for $\msS_1, \mL, \mE$ introduced in \cref{eq:S1,eq:L-defn,eq:E-defn}, where the parameter $L$ is available for optimization.

By combining \cref{prop:kl-fr-to-kl-sums,prop:using-weil-bound}, we obtain
\[
\begin{aligned}
    \msS_1^2 
    &\ll (MN)^{o(1)} \|\alpha\|^2 \|\beta\|^2 \|\nu\|^2 
    MNA
    \Bigg(
    \frac{(M+N)(N+A^2)}{N^2AL} + \frac{(M+N)}{M\sqrt{NL}} 
    \\ 
    &+
    \frac{|\vartheta|A}{MN} + \frac{1}{L^{3/2}\sqrt{A}}
    +
    \frac{1}{\sqrt{MA}}
    + 
    \frac{1}{\sqrt{N} A L^3} 
    \sup_{\substack{G, G_2 \le 30NL/M \\ D_0 \le D \le 30\frac{NL}{GM}}} 
    \sqrt{\frac{G_2 D_0}{GD^2} \frac{G^{3/2}}{G_2^{3/2} D_0^2} N^{1/2} L^7 D^4 A^2}
    \Bigg).
\end{aligned}
\]
The expression inside the supremum is non-increasing in $D_0$ and non-decreasing in $D$, so we may substitute $D_0 \asymp 1$ and $D \asymp \frac{NL}{GM}$. The supremum is therefore
\[
    \ll \sup_{G, G_2 \le 30NL/M} \sqrt{\frac{G_2}{G} \frac{G^{3/2}}{G_2^{3/2}} N^{1/2} L^7 \frac{N^2L^2}{G^2M^2} A^2}
    \ll 
    \sqrt{N^{1/2} L^7 \frac{N^2L^2}{M^2} A^2}.
\]
Plugging this into the previous bound gives
\[
\begin{aligned}
    \msS_1^2 
    \ll (MN)^{o(1)} \|\alpha\|^2 \|\beta\|^2 \|\nu\|^2 
    MNA
    \Bigg(
    \frac{(M+N)(N+A^2)}{N^2AL} + \frac{(M+N)}{M\sqrt{NL}} 
    \\ 
    +
    \frac{|\vartheta|A}{MN} + \frac{1}{L^{3/2}\sqrt{A}}
    +
    \frac{1}{\sqrt{MA}}
    + 
    \frac{N^{3/4} L^{3/2}}{M} 
    \Bigg).
\end{aligned}
\]
We may now pick the parameter $L$. Note that we have
\[
    \frac{(M+N)(N+A^2)}{N^2AL} + \frac{(M+N)}{M\sqrt{NL}} 
    + 
    \frac{1}{L^{3/2}\sqrt{A}}
    \ll
    \frac{N^{3/4} L^{3/2}}{M} 
\]
provided that
\[
    L \gg L_0 := \left(\frac{M(M+N)(N+A^2)}{N^{11/4}A}\right)^{2/5} + \frac{\sqrt{M+N}}{N^{5/8}}
    +
    \frac{M^{1/3}}{N^{1/4}A^{1/6}}.
\]
We pick $L := L_0+C_0+2\log(2|\vartheta|M)$, which satisfies \cref{eq:L-defn}. This gives
\[
\begin{aligned}
    \msS_1^2 
    &\ll (MN)^{o(1)} \|\alpha\|^2 \|\beta\|^2 \|\nu\|^2 
    MNA
    \Bigg(
    \frac{|\vartheta|A}{MN} + 
    \frac{1}{\sqrt{MA}}
    \\ 
    &+ 
    \frac{N^{3/4}}{M}
    \Bigg(
    \Big(\frac{M(M+N)(N+A^2)}{N^{11/4}A}\Big)^{3/5} + \Big(\frac{\sqrt{M+N}}{N^{5/8}}\Big)^{3/2}
    +
    \Big(\frac{M^{1/3}}{N^{1/4}A^{1/6}}\Big)^{3/2} + 1\Bigg)
    \Bigg).
\end{aligned}
\]
We can absorb the term $1/\sqrt{MA}$ since
\[
    \frac{1}{\sqrt{MA}}\ll \frac{N^{3/4}}{M}
    \left(\frac{M^{1/3}}{N^{1/4}A^{1/6}}\right)^{3/2}.
\]
If $M\le N$, then $\sqrt{M+N}/N^{5/8}\ll1$. If $M\ge N$, then $(N+A^2)/A \ge 2\sqrt N$ and $M\gg N^{11/12}$, so
\[
    \frac{\sqrt{M+N}}{N^{5/8}}\ll \frac{M^{4/5}}{N^{9/10}}
    \ll \left(\frac{M(M+N)(N+A^2)}{N^{11/4}A}\right)^{2/5}.
\]
Expanding the remaining terms gives the desired bound.
%\fix{After optimizing in $L$, simplify the final expression using \cref{lem:large-sieve-large-powers} with $C = 1$. And finally, write it as $M^2N^2A^2 \cdot $ something and simplify this in light of the trivial bound!}
\end{proof}

\begin{corollary} \label{cor:main-weil-bound}
Let $N, A \in \Z_+$ and $(\alpha_m)_{m = 1}^N, (\beta_n)_{n = 1}^N, (\nu_a)_{a = 1}^A$ be complex sequences. Then for any $\vartheta \in \Z \setminus \{0\}$, one has
\[
\begin{aligned}
    \mathop{\sum_{m=1}^N \sum_{n=1}^N}_{(m, n) = 1} \sum_{a = 1}^A \alpha_m \beta_n  \nu_a\,
    e\left(\vartheta a \frac{\bar{m}}{n} \right)
    &\ll 
    \|\alpha\|\|\beta\|\|\nu\| N^{1+o(1)}\sqrt{A} 
    \\
    &\times \left(\frac{\sqrt{|\vartheta|A}}{N} + \frac{1}{N^{1/20}A^{3/10}} + \frac{1}{N^{1/16} A^{1/8}}\right).
\end{aligned}
\]    
\end{corollary}

\begin{proof}[Proof of \cref{cor:main-weil-bound}]
As before, we may assume that $|\vartheta|A\le N^2$ in light of the trivial bound $\|\alpha\|\|\beta\|\|\nu\| N\sqrt A$. We put $m,n,a$ into dyadic ranges up to $N$ and $A$. To apply \cref{lem:large-sieve-large-powers} with $C=1$, we first split according to $(n,\vartheta)$ as in the proof of \cref{cor:main-weil-bound-unbalanced}. Combining the two estimates, and using that their right-hand sides including the normalizing factors are non-decreasing in the two lengths, we obtain
\[
\begin{aligned}     
   \sum_{\substack{m, n \le N \\ (m, n) = 1 }} \sum_{a = 1}^A \alpha_m \beta_n  \nu_a\, 
    e\left(\vartheta a \frac{\bar{m}}{n} \right)
    &\ll \|\alpha\| \|\beta\| \|\nu\| 
    N^{1+o(1)} \sqrt{A}
    \\
    &\times
    \min \Bigg(\frac{1}{A} + \frac{1}{N},
    \frac{|\vartheta|A}{N^2} 
    +
    \frac{(N+A^2)^{3/5}}{N^{7/10}A^{3/5}} 
    +
    \frac{1}{N^{1/8}A^{1/4}} 
    + \frac{1}{N^{1/4}}
    \Bigg)^{1/2}.
\end{aligned}
\]
The term $N^{-1/4}$ can be dropped: if $A \ll \sqrt{N}$ then $N^{-1/4} \ll N^{-1/8} A^{-1/4}$, and if $A \gg \sqrt{N}$ then $(N+A^2)^{3/5} N^{-7/10} A^{-3/5} \gg A^{-1} + N^{-1}$ (so the first entry in the minimum dominates with or without the $N^{-1/4}$ term). Expanding $(N+A^2)^{3/5}$, the minimum then becomes
\[
    \asymp \min \Bigg(\frac{1}{A} + \frac{1}{N},
    \frac{|\vartheta|A}{N^2} 
    +
    \frac{1}{N^{1/10}A^{3/5}}
    +
    \frac{A^{3/5}}{N^{7/10}} 
    +
    \frac{1}{N^{1/8}A^{1/4}} 
    \Bigg).
\]
We claim that the term $A^{3/5}N^{-7/10}$ can also be dropped. Firstly, if $A \ll N^{1/6}$ then $A^{3/5}N^{-7/10} \le N^{-1/8} A^{-1/4}$. Secondly, if $N^{1/6} \ll A \ll N$, then $N^{-1/8}A^{-1/4} \gg A^{-1} + N^{-1}$ (so the first entry in the minimum dominates with or without the $A^{3/5} N^{-7/10}$ term).
Thirdly, if $A \gg N$, then we similarly have $N^{-1/8}A^{-1/4} \gg A^{-1}$ and $|\vartheta|AN^{-2} \gg N^{-1}$.

After ignoring the term $A^{3/5}N^{-7/10}$ and then the first entry in the minimum for an upper bound, we reach the desired bound by expanding the square-root.
\end{proof}

\section{Using the spectral theory of automorphic forms} \label{sec:kuznetsov-exc-spectrum}

Here we assume all of the notation from \cref{subsec:aut-forms-kuznetsov}.

\subsection{A large sieve with level averaging} \label{subsec:large-sieve}

As was briefly mentioned in \cref{subsec:aut-forms-kuznetsov},
Deshouillers--Iwaniec \cite{deshouillers1982kloosterman} proved several on-average bounds for the spectral parameters and Fourier coefficients of automorphic forms of level $q$, known as spectral large sieve inequalities.
The author's previous work in \cite{pascadi2026large} built on the framework of Deshouillers--Iwaniec in the exceptional spectrum, by exploiting the special additive structure of the sequences $(\alpha_n)$ that appear inside the spectral large sieve in some applications of Linnik's dispersion method \cite{linnik1963dispersion}. Here, we take this a step further by allowing averaging over the levels $q$, where the sequence $(\alpha_n)$ may depend on $q$ in a very specific way -- more precisely, we have $\alpha_n \approx e(\tfrac{tn}{q})$ for some $t \in \Z/q\Z$. The resulting large sieve inequality for exceptional Maass forms, given below, will serve as a key technical ingredient to our main results.

\begin{theorem}[Exceptional-spectrum large sieve with level-dependent phases] \label{thm:exc-ls-level-avg}
Let $X, Q, N, Z \gg 1$ and $\eps \in (0, 1)$.
Let $(w_{q,t})$ be complex weights supported on $q \asymp Q$ and $t \in \Z/q\Z$, and let $\|w\|_1 := \sum_{q,t} |w_{q,t}|$. Let $\Psi : \R \to \C$ be a smooth function with $\Psi(x)$ supported in $x \asymp 1$ and $\Psi^{(j)} \ll_j Z^j$ for all $j \ge 0$. Let 
\begin{equation} \label{eq:exc-ls-W}
    W :=
    \max_{\substack{|x|, |y| \le ZX^{1+\eps} \\ C_1 \frac{NX}{Q} \le c \le C_2 \frac{NX}{Q} \\ c \nmid x-1,\ c \nmid x+1}}
    \sum_{\substack{q \in \Z_+ \\ t \pmod{q} \\ (x-ct)(y+ct) \equiv 1 \pmod{cq}}}
    |w_{q,t}|,
\end{equation}
where $C_1, C_2$ are sufficiently small, respectively large, constants.
Then with the notation from \cref{eq:on-basis}, one has
\begin{equation} \label{eq:exc-ls-level-avg}
    \sum_{\substack{q \in \Z_+ \\ t \pmod{q}}}
    w_{q,t}
    \sum_{\substack{f \\ \lambda_f < \frac{1}{4}}}^{\Gamma_0(q)}
    X^{2\theta_f}
    \left\vert
    \sum_n \Psi\left(\frac{n}{N}\right) e\left(\frac{tn}{q}\right) \rho_{f\infty}(n)
    \right\vert^2
    \ll_\eps
    Z^{O(1)} (QNX)^{3\eps} 
    \|w\|_1 \left(1 + \frac{N}{Q}\right) N
\end{equation}
in the range
\begin{equation} \label{eq:exc-ls-range}
    X W \ll (QNX)^{o(1)} \|w\|_1 \left(1 + \frac{N}{Q}\right).
\end{equation}
\end{theorem}

\begin{remark}
Except for the $X^{3\eps}$ factor, the right-hand side of \cref{eq:exc-ls-level-avg} does not depend on $X$, and matches the contribution of large sieve inequalities for the regular spectrum \cite[Theorem 2]{deshouillers1982kloosterman}. The goal of exceptional-spectrum results is to keep the regular spectrum dominant for $X$ as large as possible, so that the factor $X^{2\theta_f}$ gets absorbed on average; this was illustrated in a simpler case by \cref{lem:density-thm}.

To the author's knowledge, \cref{thm:exc-ls-level-avg} is the first exceptional-spectrum large sieve inequality in which averaging over $q$ improves the admissible range of $X$ for level-dependent sequences. Although we stated this result for general weights $w_{q,t}$, a motivating case of interest is when $w_{q,t} = \one_{q \sim Q} \one_{t \equiv t_q \pmod{q}}$, where $t_q \in \Z/q\Z$ is such that the congruence $(x-ct_q)(y+ct_q) \equiv 1 \pmod{cq}$ allows us to absorb the sum over $q$ from \cref{eq:exc-ls-W} by the divisor bound. In this model case we have $W \ll (QNX)^{o(1)}$ and $\|w\|_1 \asymp Q$, so \cref{eq:exc-ls-range} holds provided that
\begin{equation}\label{eq:model-range}
    X \ll \max(N, Q).
\end{equation}
Let us compare this with two relevant previous results.
\begin{itemize}
    \item[$(i)$.] Applying \cite[Theorem~1.5]{pascadi2026large} at each level and summing trivially over $q$ gives the uniform range $X \ll \max(\sqrt{N}, Q/\sqrt{N}, N/Q)$; so \cref{eq:model-range} allows $X$ larger by a factor of $\sqrt{N}$ or $Q$.
    \item[$(ii)$.] For $w_{q,t} = \one_{q \sim Q} \one_{t = 0}$, Deshouillers--Iwaniec \cite[Theorem~7]{deshouillers1982kloosterman} obtained the range $X \ll \max(N,Q^2/N)$ by averaging over $q$. Their proof applies more generally to a common phase $e(n\omega)$ with $\omega \in \R/\Z$ independent of $q$ (see, e.g., \cite[Theorem~3.11]{pascadi2026large}).
\end{itemize}
\end{remark}

\begin{example} \label{ex:fixed-t0}
Let $t_0$ be an integer independent of $q$ with $|t_0| \ll (QNX)^{O(1)}$, and suppose that $w_{q,t}$ is supported on $t \equiv t_0 \pmod q$. Then every admissible $q$ divides $(x-ct_0)(y+ct_0)-1$, which is nonzero since $c \nmid x-1$ and $c \nmid x+1$. So $W \ll (QNX)^{o(1)}$, and we attain the range in \cref{eq:model-range}.
\end{example}

\begin{example} \label{eq:n2+1}
Take $w_{q,t}$ to be supported (and of size $1$) on $t^2+1 \equiv 0 \pmod q$. Consider $x,y,c$ occurring in \cref{eq:exc-ls-W} (for which the inner sum is nonempty). Then $c \mid xy-1$, and on expanding the congruence modulo $cq$, dividing by $c$, and using $t^2\equiv-1\pmod q$, we obtain
\[
    \frac{xy-1}{c}+c\equiv -(x-y)t\pmod q.
\]
Squaring then gives
\[
    q\mid \left(\frac{xy-1}{c}+c\right)^2+(x-y)^2.
\]
The integer on the right is positive, since otherwise $x = y$ and $x^2+c^2 = 1$, which is impossible for $c \ge 2$ (and $c = 1$ is excluded by $c \nmid x-1$). Thus the divisor bound controls the number of choices for $q$ in \cref{eq:exc-ls-W}, and for each $q$ there are at most $2^{\omega(q)}$ roots of $t^2+1 \equiv 0\pmod q$ by the Chinese remainder theorem. So $W$ is once again divisor-bounded, giving the range in \cref{eq:model-range}.

Such weights $w_{q,t}$ occur in the study of the greatest prime factor of $n^2+1$. In \cite[\S6]{pascadi2026large}, the exponential sums arising in this problem were bounded using a large sieve for exponential phases without level averaging. Grimmelt--Merikoski's improvement \cite[\S1.2]{grimmelt2025greatest} exploited this average on the geometric side of an automorphic kernel, also absorbing the sum over $q$ by the divisor bound. One can view this case of \cref{thm:exc-ls-level-avg} as a spectral-side counterpart to Grimmelt--Merikoski's argument.
\end{example}

\begin{example} \label{ex:tq}
In our current application, we will need large sieve inequalities roughly of the shape in \cref{eq:sketch-large-sieve}, with $q = \ell_1\ell_1'\ell_2\ell_2'dd'$ and residues $t_q$ as in \cref{eq:sketch-phase}.
If all $6$ factors of $q$ are pairwise coprime, these residues satisfy
$t_q \equiv 0 \pmod{\ell_1\ell_1'\ell_2\ell_2'}$, $t_q \equiv \ell_1' \ell_2 \pmod{d}$, and $t_q \equiv \ell_1 \ell_2' \pmod{d'}$. So the congruence $(x-ct_q)(y+ct_q) \equiv 1 \pmod{cq}$ from \cref{eq:exc-ls-W} implies that
\[
    \ell_1\ell_1'\ell_2\ell_2'\mid xy-1, \quad d \mid (x-c\ell_1'\ell_2)(y+c\ell_1'\ell_2)-1
    \quad 
    d' \mid (x-c\ell_1\ell_2')(y+c\ell_1\ell_2')-1.
\]
This allows us to absorb the sums over $d, d', \ell_1, \ell_1', \ell_2, \ell_2'$ by the divisor bound in generic cases, reaching once again the range $X \ll \max(N, Q)$ from \cref{eq:model-range}; in \cref{eq:sketch-large-sieve}, this means $X_2 \ll L^6$, as we claimed.
We give the details in the proof of \cref{prop:using-spectral-bounds}.
\end{example}

\begin{proof}[Proof of \cref{thm:exc-ls-level-avg}]
By \cref{lem:to-exc-large-sieve}, it suffices to show that
\[
    \mL := \sum_{\substack{q \in \Z_+ \\ t \pmod{q}}} |w_{q,t}|
    \left\vert
    \sum_{c \in \Z_+}
    \frac{1}{cq}
    \sum_{m, n} \Psi\left(\frac{m}{N}\right) \bar\Psi\left(\frac{n}{N}\right) e\left(\frac{tm - tn}{q}\right) S(m, n; cq)\, \Phi\left(\frac{\sqrt{mn}}{cq}X\right)
    \right\vert
\]
is bounded by the right-hand side of \cref{eq:exc-ls-level-avg}, where $\Phi$ is a smooth function as in \cref{lem:to-exc-large-sieve}. We may eliminate the absolute values around the sum over $c$ by incorporating unit factors into the weights, which we again denote by $w_{q,t}$.
We then expand the Kloosterman sum and write
\[
    F_z(u, v) := \Psi(u) \bar\Psi(v) \Phi\left(z\sqrt{uv}\right)
\]
to reach
\[
    \mL
    =
    \sum_{\substack{q \in \Z_+ \\ t \pmod{q}}} w_{q,t}
    \sum_{c \in \Z_+}
    \frac{1}{cq}
    \sum_{\substack{r, s \pmod{cq} \\ rs \equiv 1 \pmod{cq}}}
    \sum_{m, n} e\left(\frac{m(r+ct) + n(s-ct)}{cq}\right) F_{\frac{NX}{cq}}\left(\frac{m}{N}, \frac{n}{N}\right),
\]
where the supports of $w_{q,t}$, $\Psi$ and $\Phi$ force $q \asymp Q$ and $cq \asymp NX$, so that $c \asymp NX/Q$.

Applying Poisson summation in $m,n$ and collecting the frequencies $x \equiv r+ct \pmod{cq}$ and $y \equiv s-ct \pmod{cq}$, we obtain
\[
    \mL
    =
    \sum_{\substack{q \in \Z_+ \\ t \pmod q}} w_{q,t}
    \sum_{c \asymp NX/Q} \frac{N^2}{cq}
    \sum_{\substack{x,y \in \Z \\ (x-ct)(y+ct) \equiv 1 \pmod{cq}}}
    \widehat F_{\frac{NX}{cq}}\left(-\frac{Nx}{cq},-\frac{Ny}{cq}\right),
\]
where $\widehat F_z(\xi,\eta) := \int_{\R^2} F_z(u,v)e(-u\xi-v\eta)\,du\,dv$. Since $NX/(cq) \asymp 1$, repeated integration by parts gives, for any $B \ge 2$,
\[
    \widehat F_{\frac{NX}{cq}}\left(-\frac{Nx}{cq},-\frac{Ny}{cq}\right)
    \ll_B Z^2
    \left(1+\frac{|x|}{ZX}\right)^{-B}
    \left(1+\frac{|y|}{ZX}\right)^{-B}.
\]
Dropping the congruence condition, we see that the contribution of $\max(|x|,|y|)>ZX^{1+\eps}$ is
\[
    \ll_B Z^2 \frac{N}{X} \|w\|_1 \frac{NX}{Q}
    (ZX)^2 \left(1+X^\eps\right)^{1-B}
    \ll_\eps Z^4 \|w\|_1 \frac{N^2}{Q},
\]
on taking $B>1+2/\eps$.

We next consider the terms with $x \equiv \pm 1 \pmod c$. For fixed $q,t,c$ and either choice of sign, $x$ lies in one residue class modulo $c$. Once $x$ is fixed, the congruence $(x-ct)(y+ct)\equiv 1 \pmod{cq}$ determines at most one residue class for $y$ modulo $cq$. Since the sum of $(1+|x|/(ZX))^{-2}$ over any residue class modulo $d$ is $\ll 1+ZX/d$, the contribution of these terms is
\[
    \ll Z^2 \frac{N}{X} \|w\|_1 \frac{NX}{Q}
    \left(1+\frac{ZQ}{N}\right)
    \left(1+\frac{Z}{N}\right)
    \ll Z^4 \|w\|_1 \left(1+\frac{N}{Q}\right)N.
\]
This is absorbed by the right-hand side of \cref{eq:exc-ls-level-avg}.

For the remaining terms, reducing the congruence modulo $c$ gives $c \mid xy-1$. Moreover, $xy=1$ would force $x=y=\pm1$, which we have already excluded. Thus the remaining contribution is
\[
    \ll Z^2 \frac{N}{X}
    \sum_{\substack{|x|,|y|\le ZX^{1+\eps} \\ xy\ne 1}}
    \sum_{\substack{c\asymp NX/Q \\ c\mid xy-1 \\ c\nmid x-1,\ c\nmid x+1}}
    \sum_{\substack{q \in \Z_+ \\ t \pmod q \\ (x-ct)(y+ct)\equiv 1 \pmod{cq}}}
    |w_{q,t}|.
\]
For each $x,y$ in this sum, there are at most $\tau(|xy-1|) \ll (ZQNX)^{o(1)}$ choices of $c$, by the divisor bound. With $W$ as in \cref{eq:exc-ls-W}, we can therefore bound the remaining contribution by
\[
    \ll 
    (QNX)^{o(1)} Z^{O(1)}
    \frac{N}{X} X^{2+2\eps} W,
\]
which is dominated by the right-hand side of \cref{eq:exc-ls-level-avg} in the range \cref{eq:exc-ls-range}.
\end{proof}

\subsection{Spectral bounds for sums of Kloosterman sums}

We now combine \cref{thm:exc-ls-level-avg} with the Kuznetsov formula, following \cite[Section 5]{pascadi2026large}.

\begin{corollary}[Kloosterman sums averaged over levels] \label{cor:kl-sums-level-avg}
Let $R,S,M\ge\tfrac12$, $N,C\gg1$, $X\ge1$, and $0<\eps<1$. Let $(w_{r,s,t})$ be nonnegative weights supported on $r\sim R$, $s\sim S$, $(r,s)=1$, and $t\in\Z/rs\Z$. For each triple in the support let $m=m_{r,s,t}\sim M$ be a positive integer, and let
\[
    \|w\|_1 := \sum_{r,s,t}w_{r,s,t} > 0,\qquad
    \kappa := \max_{w_{r,s,t}>0}(m_{r,s,t},rs),
\]
Let $F_{r,s,t}(x,y)$ be smooth functions supported in a fixed box $x,y\asymp1$ in $(0,\infty)^2$, with $\partial_x^i\partial_y^jF_{r,s,t}\ll_{i,j}1$ uniformly in $r,s,t$. With the constants $C_1, C_2$ from \cref{eq:exc-ls-W}, put
\begin{equation}\label{eq:kl-level-multiplicity}
    W = W(X) :=
    \max_{\substack{|x|,|y|\le 2X^{1+\eps} \\
        C_1NX/(RS)\le j\le C_2NX/(RS) \\ j\nmid x-1,\ j\nmid x+1}}
    \sum_{\substack{r\sim R,\, s\sim S \\ t\pmod{rs} \\
        (x-jt)(y+jt)\equiv1\pmod{jrs}}}w_{r,s,t},
\end{equation}
where $x,y\in\Z$, $j\in\Z_+$, and an empty maximum is taken to be zero. Then for either choice of sign one has
\begin{equation}\label{eq:kl-sums-level-avg}
\begin{aligned}
    \sum_{\substack{r\sim R,\, s\sim S \\ t\pmod{rs}}}w_{r,s,t}
    \bigg|\sum_{n\in\Z_+}e\left(\frac{tn}{rs}\right)
    \sum_{\substack{c\in\Z_+ \\ (c,r)=1}}\frac1c
    F_{r,s,t}\left(\frac{n}{N},\frac{c}{C}\right)
    S(m_{r,s,t}\bar r,\pm n;sc)\bigg| \\
    \ll_\eps (2+RSMNCX)^\eps\|w\|_1 S\sqrt{RN}
    \left(1+\frac{MN}{RS^2C^2}+\frac{\sqrt{\kappa M}}{RS}\right)^{1/2} \\
    \times
    \left(1+\min\left(\frac{N}{RS},\frac{SC^2}{M}\right)\right)^{1/2}
    \left(1+\frac{C^2}{RNX}\min\left(\kappa,\frac{R^2S^2}{M}\right)\right)^{\theta_{\max}}
\end{aligned}
\end{equation}
in the range
\[
    XW \ll (2+RSNX)^{o(1)}\|w\|_1\left(1+\frac{N}{RS}\right).
\]
Here $\bar r$ denotes the inverse of $r$ modulo $sc$.
\end{corollary}

\begin{proof}
We follow the proof of \cite[Corollary~5.5]{pascadi2026large}, keeping the average over levels when applying the large sieve to the sum over $n$. For fixed $r,s,t$ write $q=rs$ and $m=m_{r,s,t}$. With the scaling matrices from \cite[(3.9)]{pascadi2026large}, the Kloosterman sum in \cref{eq:kl-sums-level-avg} is the one between the cusps $\infty$ and $1/s$ of $\Gamma_0(q)$, of generalized modulus $\gamma=s\sqrt r\,c$. As $c$ runs over the positive integers coprime to $r$, $\gamma$ runs over all admissible moduli, by \cite[(3.12)]{pascadi2026large}. We have $1/c=s\sqrt r/\gamma$, which gives a factor $O(S\sqrt R)$ when passing to the Kuznetsov formula.

Let $G_{r,s,t}(\xi,y)$ be the Fourier transform in $x$ of $F_{r,s,t}(x,\sqrt{x}/y)$, extended by zero outside $x,y>0$. It is supported in $y\asymp1$, and for all $i,B\ge0$ we have
\[
    \partial_y^iG_{r,s,t}(\xi,y)\ll_{i,B}(1+|\xi|)^{-B}.
\]
Take a fixed smooth function $\Psi$ supported in $x\asymp1$ and equal to $1$ on the common $x$-support of the $F_{r,s,t}$. By Fourier inversion,
\[
    F_{r,s,t}\left(\frac{n}{N},\frac{c}{C}\right)
    =\Psi\left(\frac{n}{N}\right)\int_{\R}e\left(\frac{\xi n}{N}\right)
    G_{r,s,t}\left(\xi,\frac Cc\sqrt{\frac{n}{N}}\right)\,d\xi .
\]
Thus the sequence in $n$ and the test function of $u=4\pi\sqrt{mn}/\gamma$ in the Kuznetsov formula are respectively
\[
    n\mapsto\Psi(n/N)e\left(\frac{\xi n}{N}+\frac{tn}{q}\right),
    \qquad
    u\mapsto G_{r,s,t}\left(\xi,\frac{s\sqrt r\, C}{4\pi\sqrt{mN}}u\right).
\]
The test function is supported in $u\asymp\sqrt{MN}/(S\sqrt R\,C)$, with the same decay in $\xi$ after any fixed number of applications of $u\partial_u$.

Let $0<\eta<\eps/2$ be sufficiently small in terms of $\eps$. Both cusps have $\mu=q^{-1}$, and the squared norm of the sequence in $n$ is $O(N)$. For the regular Maass spectrum, the Eisenstein series and the holomorphic spectrum, apply the Bessel bounds and the large sieve estimates as in \cite[Lemmas~3.8--3.9 and (5.23)]{pascadi2026large}. Summing with the weights $w_{r,s,t}$ and using $q\asymp RS$, $m\asymp M$ and $(q,m)\le\kappa$, the regular contribution to the integrand is
\[
\begin{aligned}
    &\ll_{\eta,B}(2+RSMNC)^\eta(1+|\xi|)^{-B}
      \frac{\|w\|_1 S\sqrt{RN}}{1+\sqrt{MN}/(S\sqrt R\,C)} \\
    &\qquad\times\left(1+\frac{MN}{RS^2C^2}+\frac{\sqrt{\kappa M}}{RS}\right)^{1/2}
      \left(1+\frac{MN}{RS^2C^2}+\frac{N}{RS}\right)^{1/2}.
\end{aligned}
\]
This gives the first two factors in parentheses in \cref{eq:kl-sums-level-avg}, since
\[
    \frac{1+MN/(RS^2C^2)+N/(RS)}{(1+\sqrt{MN}/(S\sqrt R\,C))^2}
    \asymp 1+\min\left(\frac{N}{RS},\frac{SC^2}{M}\right).
\]

For the exceptional spectrum, put
\[
    X_1=\max\left(1,\frac{R^2S^2}{\kappa M}\right).
\]
By \cite[(3.19) and (5.24)]{pascadi2026large}, the Bessel transform is bounded, up to a factor of $O_{\eta,B}((2+RSMNC)^\eta(1+|\xi|)^{-B})$, by
\[
    \frac{1+(RS^2C^2/(MN))^{\theta_f}}{1+\sqrt{MN}/(S\sqrt R\,C)} .
\]
In the numerator we use
\[
    1+\left(\frac{RS^2C^2}{MN}\right)^{\theta_f}
    \ll X_1^{\theta_f}X^{\theta_f}
    \left(1+\frac{C^2}{RNX}\min\left(\kappa,\frac{R^2S^2}{M}\right)\right)^{\theta_{\max}}.
\]
We keep the last factor outside Cauchy--Schwarz and apply the latter over $r,s,t$ and $f$. For the sum of squares at $m$, we have $X_1\ll\max(1,q^2/((q,m)m))$ on the whole support of $w$, so \cref{lem:density-thm} gives
\[
    \sum_{r,s,t} w_{r,s,t}
    \sum_{\substack{f \\ \lambda_f<1/4}}^{\Gamma_0(q)}
    X_1^{2\theta_f}|\rho_{f\infty}(m)|^2
    \ll_\eta(2+RSM)^\eta\|w\|_1\left(1+\frac{\sqrt{\kappa M}}{RS}\right).
\]
For the sum of squares in $n$, the scaling matrix at $1/s$ in \cite[(3.9)]{pascadi2026large} is
\[
    \frac1{\sqrt r}\begin{pmatrix}r&-b\\q&ra\end{pmatrix},
    \qquad a,b\in\Z,\quad ra+sb=1.
\]
This is an Atkin--Lehner matrix for the exact divisor $r$ of $q$. Composition with it changes the orthonormal basis unitarily on each Laplace eigenspace, giving at $\infty$ the Fourier coefficients of the original basis at $1/s$. For the negative sign we also compose with the reflection $z\mapsto-\bar z$. Thus we may replace $\rho_{f,1/s}(\pm n)$ by $\rho_{f\infty}(n)$ in the sum of squares and collect the weights as
\[
    w_{q,t}=\sum_{rs=q}w_{r,s,t},\qquad
    \sum_{q,t}w_{q,t}=\|w\|_1.
\]

We apply \cref{thm:exc-ls-level-avg} with $Q=RS$, test function $\Psi(x)e(\xi x)$, $Z=1+|\xi|$, and $\eps$ there equal to $\eta$. If $|\xi|\le X^{\eps/2}$, then $ZX^{1+\eta}\le 2X^{1+\eps}$, and the maximum in \cref{eq:exc-ls-W} for $w_{q,t}$ is at most $W$. Hence the range condition in the corollary gives \cref{eq:exc-ls-range}. If $|\xi|>X^{\eps/2}$, we apply the theorem with $X$ replaced by $1$, when its range condition holds since the maximum is at most $\|w\|_1$. In this case $X^{2\theta_f}\le (1+|\xi|)^{4\theta_{\max}/\eps}$. In either case we obtain
\[
\begin{aligned}
    &\sum_{q,t}w_{q,t}\sum_{\substack{f \\ \lambda_f<1/4}}^{\Gamma_0(q)}
      X^{2\theta_f}\bigg|\sum_n\Psi\left(\frac{n}{N}\right)
      e\left(\frac{\xi n}{N}+\frac{tn}{q}\right)\rho_{f\infty}(n)\bigg|^2 \\
    &\qquad
    \ll_{\eps,\eta}(2+RSNX)^{3\eta}(1+|\xi|)^{O_\eps(1)}
      \|w\|_1\left(1+\frac{N}{RS}\right)N .
\end{aligned}
\]
Taking the square roots of the two bounds, the factors $1+\sqrt{\kappa M}/(RS)$ and $1+N/(RS)$ are bounded by the corresponding factors in the regular contribution before simplification. We recover the same bound with the additional $\theta_{\max}$-factor in \cref{eq:kl-sums-level-avg}. Integrating in $\xi$ with $B$ sufficiently large in terms of $\eps$, and recalling the choice of $\eta$, proves the result.
\end{proof}

\subsection{Completing the proof}
We now resume the notation from \cref{subsec:setup}. Throughout this subsection, let $\theta \in [\theta_{\max}, 7/64]$ for $\theta_{\max}$ as in \cref{eq:theta-max}; in particular, the choice $\theta = 7/64$ is admissible due to Kim--Sarnak \cite{kim2003functoriality}, and $\theta = 0$ would be admissible under Selberg's eigenvalue conjecture.

\begin{proposition}\label{prop:using-spectral-bounds}
Let $\eps>0$ and let $\msS_6$ be as in \cref{eq:S6}, in the ranges of \cref{prop:kl-fr-to-kl-sums}. Suppose that $M\ge N$ and $|\vartheta|A\le MN$. Then one has
\[
    \msS_6\ll_\eps (MN)^{O(\eps)}
    \frac{G^3D^4L^8}{G_2D_0}A^2\|\nu\|^4
    \left(1+\frac{|\vartheta|A}{(GD)^3L^5}\right)^{1/4-\theta}
    \left(1+\frac{|\vartheta|A}{(GD)^3L^5}
             +\frac{N^2}{(GD)^6L^8}\right)^\theta.
\]
\end{proposition}

\begin{proof}
For each tuple up to and including $z$ in \cref{eq:S6}, write
\[
    b=d_0|\delta\delta'|g,\qquad P=\ell_1\ell_1'\ell_2\ell_2',\qquad
    s=\ell_2\ell_2',\qquad r=\ell_1\ell_1' b g_0.
\]
The level, the fixed argument of the Kloosterman sum, and the numerator of the additive phase are respectively
\[
    q=rs=Pbg_0,\qquad k=d_0|\vartheta\Delta_1|,\qquad
    t=Pg_0 z\bar s\pmod q,
\]
where $\bar s$ is taken modulo $b$. Indeed, $(r,s)=1$, the condition on $u_2$ is $(u_2,r)=1$, and since $(d_0,su_2)=1$ we have
\[
    d_0\vartheta\Delta_1\bar r\equiv
    \vartheta\Delta_1\overline{\ell_1\ell_1'|\delta\delta'|gg_0}
    \pmod{su_2}.
\]
We can change the signs of both arguments of a Kloosterman sum simultaneously, so we may take the first argument to be $k\bar r$.

Split according to the sign of $\vartheta\Delta_1$ and put $g_0,g_1,k,(q,k)/d_0$ into dyadic ranges of sizes $G_0,G_1,M_1,K$, respectively. Here $G_0G_1\asymp G$ and $d_0\mid(q,k)$. We also split $r,s$ into dyadic ranges $r\sim R$, $s\sim S$. There are $O(1)$ such pairs of ranges on each of the preceding blocks, since
\[
    R\asymp L^2\frac{GD^2G_0}{D_0},\qquad S\asymp L^2.
\]
We take $C=N/(GG_2)$, so that the cutoff $\Phi(u_2/C)$ restricts $u_2\asymp C$.
Split $h$ according to its sign and use a smooth partition of unity in $|h|\asymp H$. To do this at the upper endpoint of \cref{eq:S6}, we may replace the sharp cutoff by a smooth one equal to $1$ below that endpoint and vanishing above twice the endpoint. Indeed, where these cutoffs differ the argument of $\hat\Phi$ has absolute value $\gg (MN)^\eps$, since $u_2\asymp N/(GG_2)$ and $b\asymp GD^2/D_0$. By the rapid decay of $\hat\Phi$ and the trivial bound for the Kloosterman sums, this changes $\msS_6$ by $O_{\eps,B}((MN)^{-B}\|\nu\|^4)$ for any fixed $B$. The same decay allows us to use $1/2$ in the cutoff for $h$ if $\eps\ge1$, so we may assume $0<\eps<1$. The derivatives of the resulting weights in $|h|/H$ and $u_2/C$ are uniformly bounded, using $y^j\hat\Phi^{(j)}(y)\ll_j 1$. Finally,
\[
    H\ll (MN)^\eps\frac{GD^2L^2}{D_0G_2},
    \qquad k\ll\frac{|\vartheta|ADL^3}{G_1G_2}.
\]
In particular $H/(RS)\ll (MN)^\eps$ and
\[
    \frac{M_1H}{RS^2C^2}\ll (MN)^\eps\frac{|\vartheta|A GD}{N^2L}
    \ll (MN)^\eps,
\]
since $GD\ll NL/M$.

For each $r,s,t$ in the block, let $w_{r,s,t}$ be the sum of $|\nu_{a_1}\nu_{a_1'}\nu_{a_2}\nu_{a_2'}|$ over the tuples giving rise to $r,s,t$. The integer $k$ and the smooth weight may still depend on the tuple. For each $r,s,t$ we take them from a tuple maximizing the absolute value of the inner sum on this block. Multiplying this maximum by $w_{r,s,t}$ bounds the contribution of all tuples giving rise to $r,s,t$. We apply \cref{cor:kl-sums-level-avg} with $m=k$, $M,N$ there equal to $M_1,H$ here, and $c=u_2$. The smooth weights are functions of $|h|/H$ and $u_2/C$, and the factor $1/u_2$ is the factor $1/c$ in the corollary.

Put
\[
    X_0=\max\left(1,\frac{L^4D^2}{D_0G_2}\right)\asymp\frac{L^4D^2}{D_0G_2},
\]
where we used $G_2\ll NL/M\ll L$ and $D\ge D_0\gg1$. In order to choose $X\le X_0$ separately on each block, put $W_0 =\max_{1\le X\le X_0}W(X)$, with $W(X)$ as in \cref{eq:kl-level-multiplicity} for this block. The integers in these congruences are bounded by a fixed power of $2+MN$, uniformly in $X$. Either sign of the phase is allowed, since negating $t$ is equivalent to negating both $x$ and $y$. We claim that
\begin{equation}\label{eq:spectral-phase-multiplicity}
    W_0 \ll (MN)^{o(1)}GG_1 A^2\|\nu\|^4.
\end{equation}
To see this, fix a triple $x,y,j$ occurring in this maximum, with $j$ the auxiliary modulus in the congruence modulo $jq$. Fix also $g,g_0,g_1,g_2$. Since $P$ divides both $q$ and $t$, we have $P\mid xy-1$. Also $xy\ne1$, since $j\nmid x-1$ and $j\nmid x+1$. Thus the four primes have $(MN)^{o(1)}$ possibilities by the divisor bound. Next, the congruences on $z$ in \cref{eq:S6} give
\[
    t\equiv g_0g_2\ell_1'\ell_2\pmod{d_0|\delta|},
    \qquad
    t\equiv g_0g_2\ell_1\ell_2'\pmod{d_0|\delta'|}.
\]
Consequently $d_0|\delta|$ divides
\[
    (x-jg_0g_2\ell_1'\ell_2)(y+jg_0g_2\ell_1'\ell_2)-1,
\]
and $d_0|\delta'|$ divides the same expression with $\ell_1'\ell_2$ replaced by $\ell_1\ell_2'$. Both integers are nonzero, since $(x-jv)(y+jv)=1$ would imply $x\equiv\pm1\pmod j$. The divisor bound therefore gives $(MN)^{o(1)}$ choices for the two divisors. Their GCD is $d_0$, so this also bounds the number of choices for $d_0,\delta,\delta'$. All integers in these divisor bounds have polynomial size in $MN$. We sum over the four $a$'s by Cauchy--Schwarz at a cost of $O(A^2\|\nu\|^4)$, and over $z$ at a cost of $O(G_1)$ by \cref{lem:y-count}. Since each $g$ admits $(MN)^{o(1)}$ values of $g_0,g_1,g_2$ by the divisor bound and \cref{lem:spec-prime-factors}, this proves \cref{eq:spectral-phase-multiplicity}.

We will use the following two bounds for the mass of the block.
\begin{align}
    \|w\|_1+X_0 W_0 &\ll (MN)^{o(1)}
        GG_1 A^2\|\nu\|^4\frac{D^2L^4}{D_0G_2},
    \label{eq:spectral-block-count}\\
    K^{1/2}\|w\|_1&\ll (MN)^{o(1)}
        GG_1^{1/2} A^2\|\nu\|^4\frac{D^2L^4}{D_0}.
    \label{eq:spectral-block-gcd}
\end{align}
In view of \cref{eq:spectral-phase-multiplicity}, the first inequality leaves us to count the tuples in $\|w\|_1$, saving a factor of $G_2$. As in the proof of \cref{prop:using-weil-bound}, let
\[
   \tilde\Delta_1=g_1g_2\Delta_1
   =\delta\ell_1'\ell_2'(a_1\ell_1-a_2\ell_2)
      -\delta'\ell_1\ell_2(a_1'\ell_1'-a_2'\ell_2')\ne0.
\]
The congruence $\tilde\Delta_1\equiv0\pmod{g_2}$ is linear in each of $\ell_1,\ell_2,\ell_1',\ell_2'$ when the other three primes are fixed. Its constant term is respectively $a_2,a_1,a_2',a_1'$ times a unit modulo $g_2$, since $(\delta\delta'P,g_2)=1$. For example, as a congruence in $\ell_1'$ its constant term is $\delta'\ell_1\ell_2a_2'\ell_2'$. If the congruence is soluble, the GCD of its linear coefficient with $g_2$ divides its constant term. Hence, writing $a_j$ for the corresponding one of the four $a$'s, the number of choices for this prime is
\[
    \ll 1+\frac{L(a_j,g_2)}{g_2}\ll\frac{L(a_j,g_2)}{G_2},
\]
where we used $G_2\ll NL/M\ll L$. We sum trivially over the other three primes and take the least of these four bounds, which is at most their geometric mean. The four GCD factors thus occur to the power $1/4$, and for each one we have
\[
    \sum_{a\sim A}|\nu_a|(a,g_2)^{1/4}
    \le\|\nu\|\left(\sum_{a\sim A}(a,g_2)^{1/2}\right)^{1/2}
    \ll (MN)^{o(1)}\sqrt A\,\|\nu\|
\]
by \cref{lem:gcd-avg}. Summing over $d_0$ and the two switched divisors costs $O(D^2/D_0)$. We sum over the remaining variables as in \cref{eq:spectral-phase-multiplicity}, proving \cref{eq:spectral-block-count}.

For \cref{eq:spectral-block-gcd}, we use $(g,\vartheta)=(g_0,\Delta_1)=1$ to write
\[
\begin{aligned}
    \frac{(q,k)}{d_0}
    &=\big(P|\delta\delta'|gg_0,\vartheta\Delta_1\big)\\
    &\le g_1\,
    \big(\delta,\vartheta(a_1'\ell_1'-a_2'\ell_2')\big)
    \big(\delta',\vartheta(a_1\ell_1-a_2\ell_2)\big)
    (P,\tilde\Delta_1).
\end{aligned}
\]
Indeed the part from $gg_0=g_0^2g_1$ is at most $g_1$, while the other factors follow by reducing $\tilde\Delta_1$ modulo $\delta,\delta'$ and using $(\delta,\delta')=(\delta\delta',P)=(\vartheta,P)=1$. A prime occurring just once in $P$ contributes $(\ell_1,a_2)$, $(\ell_2,a_1)$, $(\ell_1',a_2')$ or $(\ell_2',a_1')$ to $(P,\tilde\Delta_1)$, respectively. A repeated prime contributes at most its square, and occurs at most twice since $(\ell_1\ell_1',\ell_2\ell_2')=1$.

Fix first the four $a$'s, the primes and $d_0,\delta,\delta'$. The displayed expression for $\tilde\Delta_1$ is independent of $g_0,g_1,g_2$, so $g_1\mid\tilde\Delta_1$ leaves $(MN)^{o(1)}$ choices for $g_1$. We sum trivially over $g_0\asymp G_0$, use \cref{lem:spec-prime-factors} for $g_2$, and use \cref{lem:y-count} for $z$. Including the square root of the factor $g_1$ in the GCD bound, these sums cost $(MN)^{o(1)}G_0G_1^{3/2}\asymp (MN)^{o(1)}GG_1^{1/2}$. Next, for each $d_0$ we sum over $\delta,\delta'$ using \cref{lem:gcd-avg}, since $a_1\ell_1-a_2\ell_2$ and $a_1'\ell_1'-a_2'\ell_2'$ are both nonzero. Together with $d_0$, these variables cost $(MN)^{o(1)}D^2/D_0$. For the square root of $(P,\tilde\Delta_1)$, each prime occurring once costs $(MN)^{o(1)}L$ by \cref{lem:gcd-avg}, and a repeated pair costs at most $\sum_{\ell\sim L}\ell\ll L^2$. Thus all four primes cost $(MN)^{o(1)}L^4$, and summing the weights over the four $a$'s finishes the proof of \cref{eq:spectral-block-gcd}.

We may assume $\|w\|_1>0$. In \cref{cor:kl-sums-level-avg} take
\[
    X=\min\left(X_0,\frac{\|w\|_1}{W_0}\right),
\]
with $X=X_0$ if $W_0 = 0$. Since $W_0 \le \|w\|_1$, we have $1\le X\le X_0$ and $XW(X)\le XW_0 \le\|w\|_1$, as required. Recalling that $q\asymp RS$, $k\sim M_1$ and $(q,k)\asymp D_0K$ on the block, we have $\kappa\ll D_0K$ in \cref{eq:kl-sums-level-avg}. Using the bounds above for $H/(RS)$ and $M_1H/(RS^2C^2)$, and writing $\min(a,b)=a/\max(1,a/b)$ for $a,b>0$, we obtain the bound
\[
\begin{aligned}
    &\ll_\eps (MN)^{O(\eps)}S\sqrt{RH}\|w\|_1
    \left(1+\frac{D_0KM_1}{R^2S^2}\right)^{1/4} \\
    &\qquad\times
    \left(1+\frac{C^2D_0K}
     {RHX\max(1, \frac{D_0 K M_1}{R^2S^2})}\right)^\theta.
\end{aligned}
\]
We may replace $X$ by $X_0$ in the last factor at the cost of $(X_0/X)^\theta$. By the choice of $X$ and the fact that $\theta\le1/2$, we have
\[
    \|w\|_1\left(\frac{X_0}{X}\right)^\theta
    \le\sqrt{\|w\|_1\big(\|w\|_1+X_0 W_0\big)}.
\]
For $y\ge0$, the function $\sqrt H(1+y/H)^\theta$ is increasing in $H$. Thus we may use $H=GD^2L^2/(D_0G_2)$ in this expression, at a cost of $(MN)^{O(\eps)}$ by the upper bound for $H$ above. The number of dyadic intervals contributes a further logarithm. At this value we have
\[
\begin{aligned}
    S\sqrt{RH}&\asymp \frac{L^4GD^2}{D_0}\sqrt{\frac{G_0}{G_2}},
    \qquad
    \frac{D_0KM_1}{R^2S^2}\ll KD_0^3\frac{|\vartheta|A}{(GD)^3L^5},\\
    \frac{C^2D_0K}{RHX_0}&\ll KD_0^4\frac{G^2}{G_0}\frac{N^2}{(GD)^6L^8}.
\end{aligned}
\]
We use $(1+p)^{1/4}(1+y/\max(1,p))^\theta\ll
(1+p)^{1/4-\theta}(1+p+y)^\theta$ for $p,y\ge0$. Taking $p=D_0KM_1/(R^2S^2)$ and $y=C^2D_0K/(RHX_0)$, we may extract $KD_0^3$ from $1+p$ and $KD_0^4G^2/G_0$ from $1+p+y$, up to absolute factors. Thus the last two factors in the corollary cost at most
\[
    \ll K^{1/4}D_0^{3/4+\theta}
    \left(\frac{G^2}{G_0}\right)^\theta
    \left(1+\frac{|\vartheta|A}{(GD)^3L^5}\right)^{1/4-\theta}
    \left(1+\frac{|\vartheta|A}{(GD)^3L^5}
             +\frac{N^2}{(GD)^6L^8}\right)^\theta.
\]
By \cref{eq:spectral-block-count,eq:spectral-block-gcd},
\[
    K^{1/4}\sqrt{\|w\|_1\big(\|w\|_1+X_0 W_0\big)}\ll (MN)^{o(1)}
    \frac{GG_1^{3/4}D^2L^4}{D_0\sqrt{G_2}}A^2\|\nu\|^4.
\]
Multiplying these estimates and summing the remaining dyadic blocks gives the claimed bound with an extra factor of
\[
    D_0^{\theta-1/4}G_0^{\theta-1/2}G_1^{2\theta-1/4}\ll1,
\]
since $D_0,G_0,G_1\gg1$ on every nonempty block and $\theta\le7/64<1/8$.
\end{proof}

\begin{remark}
In the proof of \cref{prop:using-spectral-bounds}, we have kept $\vartheta$ inside the Fourier coefficient when applying the density theorem. For large $|\vartheta|$, one could also use Hecke relations at primes away from the level and bounds towards the Ramanujan conjecture to treat factors of $\vartheta$ separately, as in \cite{assing2021uniform}.
\end{remark}

\begin{theorem}\label{thm:unbalanced-spectral-methods}
Let $M,N,A \ge \tfrac{1}{2}$ with $M\ge N$, and let $(\alpha_m)_{m\sim M}$, $(\beta_n)_{n\sim N}$, $(\nu_a)_{a\sim A}$ be complex sequences. For any $\vartheta\in\Z\setminus\{0\}$, one has
\begin{equation}\label{eq:spectral-unbalanced}
\begin{aligned}
    &\bigg|\mathop{\sum\sum}_{\substack{m\sim M, n\sim N\\(m,n)=1}} \sum_{a\sim A}
        \alpha_m\beta_n\nu_a\,
        e\left(\vartheta a\frac{\bar m}{n}\right)\bigg|^2
    \ll (MN)^{o(1)}\|\alpha\|^2\|\beta\|^2\|\nu\|^2\\
    &\quad\times\Bigg(|\vartheta|A^2+
    \sum_{(K,b)}\Big(
       K+\big(KM^{2b}\big)^{\frac1{b+1}}
       +\big(K(M^3A^2)^b\big)^{\frac1{2b+1}}
       +\big(K^3(M^2N^2A)^b\big)^{\frac1{2b+3}}
    \Big)\Bigg),
\end{aligned}
\end{equation}
where $(K,b)$ runs through the four pairs
\[
\begin{gathered}
    \big(AN^{3/2},2\big),\qquad
    \big(|\vartheta|^{1/8}A^{9/8}M^{3/8}N^{9/8},1\big),\\
    \big(AM^{3\theta}N^{3/2-2\theta},2-7\theta\big),\qquad
    \big(|\vartheta|^{1/8}A^{9/8}M^{3/8+3\theta}N^{9/8-2\theta},1-7\theta\big).
\end{gathered}
\]
For $M<N$, the same bound holds with $M,N$ interchanged on the right-hand side.
\end{theorem}

\begin{proof}
Recalling $\msS_1$ from \cref{eq:S1}, we first show that for any $L\ge1$ and $M\ge N$,
\begin{equation}\label{eq:spectral-unbalanced-amplifier}
\begin{aligned}
    |\msS_1|^2
    &\ll (MN)^{o(1)}\|\alpha\|^2\|\beta\|^2\|\nu\|^2
    \Bigg(|\vartheta|A^2+\frac{M^2}{L}+\frac{M^{3/2}A}{\sqrt L}
        +\frac{MN\sqrt A}{L^{3/2}}\\
    &\qquad+AN^{3/2}L^2
        \left(1+\frac{|\vartheta|AM^3}{N^3L^8}\right)^{1/8}
        \left(1+\frac{M^3}{N^2L^7}\right)^\theta\Bigg).
\end{aligned}
\end{equation}
We will also use the estimate obtained by replacing the last term in parentheses by
\begin{equation}\label{eq:spectral-unbalanced-coupled}
    AN^{3/2}L^2
    \left(1+\frac{|\vartheta|AM^3}{N^3L^8}
             +\frac{M^6}{N^4L^{14}}\right)^{1/8}.
\end{equation}
The claim is immediate if any of the three norms vanishes. We may assume $|\vartheta|A\le MN$, since otherwise the term $|\vartheta|A^2$ suffices by the trivial bound. We first treat $M\ge N$ and sequences $(\beta_n)$ supported on integers coprime to $\vartheta$. Take $L$ in the range \cref{eq:L-defn}. Using $1+p+y\le(1+p)(1+y)$ for $p,y\ge0$ in \cref{prop:using-spectral-bounds}, the last term in parentheses in \cref{prop:kl-fr-to-kl-sums} is
\[
    \ll_\eps (MN)^{O(\eps)}
    \sup_{GD\ll NL/M}
    \frac{GD L}{\sqrt N}
    \left(1+\frac{|\vartheta|A}{(GD)^3L^5}\right)^{1/8}
    \left(1+\frac{N^2}{(GD)^6 L^8}\right)^{\theta/2}.
\]
The expression is non-decreasing in $GD$, since $1-3/8-3\theta>0$. Substituting $GD\asymp NL/M$ gives an upper bound of
\[
    \ll_\eps (MN)^{O(\eps)}
    \frac{\sqrt N L^2}{M}
    \left(1+\frac{|\vartheta|AM^3}{N^3L^8}\right)^{1/8}
    \left(1+\frac{M^3}{N^2L^7}\right)^\theta.
\]
After multiplication by $MNA$, this is the last term of \cref{eq:spectral-unbalanced-amplifier}. To obtain \cref{eq:spectral-unbalanced-coupled}, we instead bound the two factors in \cref{prop:using-spectral-bounds} by
\[
    \left(1+\frac{|\vartheta|A}{(GD)^3L^5}
             +\frac{N^2}{(GD)^6L^8}\right)^{1/4},
\]
since $0\le\theta\le1/4$. In the supremum we now have $GD$ times the eighth root of the expression in parentheses, which is again non-decreasing in $GD$, since $1-6/8>0$. Substituting $GD\asymp NL/M$ and multiplying by $MNA$ gives \cref{eq:spectral-unbalanced-coupled}. Here and below we may absorb $(MN)^{O(\eps)}$ into $(MN)^{o(1)}$ in the final bound, since $\eps>0$ in \cref{prop:kl-fr-to-kl-sums} was arbitrary.

After the same multiplication, the other terms of \cref{prop:kl-fr-to-kl-sums} become
\[
    (M+N)\left(\frac{M}{L}+\frac{MA^2}{NL}
                     + A\sqrt{\frac NL}\right)
    +|\vartheta|A^2+\frac{MN\sqrt A}{L^{3/2}}+N\sqrt{MA}.
\]
To simplify the term involving $A^2/L$, we also have \cref{lem:large-sieve-large-powers} with $C=1$, which gives $(M+N)(N+A)$ with the same factors of $(MN)^{o(1)}$ and the norms squared. Thus we may replace $MA^2/(NL)$ in the last display by
\[
    \min\left(\frac{MA^2}{NL},N+A\right)
    \le A+\min\left(\frac{MA^2}{NL},N\right)
    \le A+A\sqrt{\frac ML}.
\]
Since $M\ge N$, this gives the first four terms of \cref{eq:spectral-unbalanced-amplifier}, together with $O(MA+N\sqrt{MA})$. If $L\le M$, then $MA\le M^{3/2}A/\sqrt L$, and otherwise $MA\ll AN^{3/2}L^2$. Similarly, if $L^3\le M$, then $N\sqrt{MA}\le MN\sqrt A/L^{3/2}$, and otherwise $N\sqrt{MA}\ll AN^{3/2}L^2$.

For a general $L\ge1$, we may assume $L\le (2+MN)^2$, since otherwise the term $AN^{3/2}L^2$ already exceeds the trivial bound squared, apart from the norms. Replacing $L$ in the argument by $L+C_0+2\log(2+2|\vartheta|M)$ enforces \cref{eq:L-defn} and changes the resulting bound by at most a power of $\log(2+MN)$. To remove the condition $(n,\vartheta)=1$, we split according to $d=(n,\vartheta)$ as in the proof of \cref{cor:main-weil-bound-unbalanced}. This replaces $N$ by $N/d$ and $\vartheta$ by $\vartheta/d$ and restricts the two sequences. For fixed $L$, the right-hand side of \cref{eq:spectral-unbalanced-amplifier}, apart from the norms, is non-decreasing in $M,N,A,|\vartheta|$. For the last term, monotonicity in $N$ follows from $3/2-3/8-2\theta>0$. The same holds for \cref{eq:spectral-unbalanced-coupled}, since $3/2-4/8>0$. We may therefore sum over $d\mid |\vartheta|$ at a cost of $(MN)^{o(1)}$.

Finally suppose $M<N$. By additive reciprocity,
\[
    e\left(\vartheta a\frac{\bar m}{n}\right)
    = e\left(-\vartheta a\frac{\bar n}{m}\right)
      e\left(\frac{\vartheta a}{mn}\right).
\]
Expanding the last exponential in a Taylor series separates $m,n,a$. For the term of degree $j$ we apply the case already proved with $M,N$ interchanged, absorbing $(M/m)^j,(N/n)^j,(a/(2A))^j$ into the respective sequences. These factors are all $1$-bounded, and the sum of the absolute values of the scalar coefficients is $\exp(4\pi|\vartheta|A/(MN))\ll1$. This gives the same estimate with $M,N$ interchanged on the right-hand side.

It remains to choose $L$ in \cref{eq:spectral-unbalanced-amplifier}. Up to an absolute factor, its last term is the maximum of $KL^b$ over the four pairs in the statement. All four terms are increasing in $L$, since $b\ge1-7\theta>0$. For each of the three decreasing terms $DL^{-p}$, with
\[
    (D,p)=(M^2,1),\quad (M^{3/2}A, \tfrac{1}{2}),\quad
          (MN\sqrt A,\tfrac32),
\]
the first value of $L>0$ at which some $KL^b$ equals $DL^{-p}$ is $\min_{(K,b)}(D/K)^{1/(b+p)}$. We therefore take
\[
    L=\max\left(1,\ \max_{(D,p)}\min_{(K,b)}
                       (D/K)^{1/(b+p)}\right).
\]
If $L>1$, the largest decreasing term and the largest increasing term are equal at this choice. If $L=1$, every decreasing term is bounded by $\max_{(K,b)}K$. In the former case the common value for the two terms attaining the maxima is $(K^pD^b)^{1/(p+b)}$. For the three choices of $(D,p)$ these are precisely the last three expressions in the summand in \cref{eq:spectral-unbalanced}. This proves the bound.

Conversely, for any $L\ge1$ we have $K\le KL^b$ and
\[
    (K^pD^b)^{1/(p+b)}
    =(KL^b)^{p/(p+b)}(DL^{-p})^{b/(p+b)}
    \le \max(KL^b,DL^{-p}).
\]
Thus taking the infimum over $L\ge1$ in \cref{eq:spectral-unbalanced-amplifier} gives \cref{eq:spectral-unbalanced} up to an absolute factor.
\end{proof}

Specializing \cref{thm:unbalanced-spectral-methods} to dyadic ranges up to $N$, we deduce the following result.

\begin{corollary} \label{cor:balanced-spectral-methods}
Let $N,A\in\Z_+$ and let $(\alpha_m)_{m=1}^N$, $(\beta_n)_{n=1}^N$, and $(\nu_a)_{a=1}^A$ be complex sequences.
For any $\vartheta\in\Z\setminus\{0\}$, one has
\[
\begin{aligned}
    \mathop{\sum_{m=1}^N\sum_{n=1}^N}_{(m,n)=1} \sum_{a=1}^A
    \alpha_m\beta_n\nu_a\,
    e\left(\vartheta a\frac{\bar m}{n}\right)
    &\ll \|\alpha\|\|\beta\|\|\nu\| N^{1+o(1)}\sqrt A\\
    &\times\Bigg(
    \frac{\sqrt{|\vartheta|A}}{N}+\frac{(|\vartheta|A)^{1/16}}{N^{(1-2\theta)/4}}
    +\frac{1}{N^{1/12}A^{1/3}}
    +\left(\frac{|\vartheta|}{N^4A^7}\right)^{1/32}\\
    &
    +\frac{1}{N^{\frac{1-2\theta}{12-28\theta}}
               A^{\frac{2-7\theta}{6-14\theta}}}
    +\left(\frac{|\vartheta|}{N^{4(1-2\theta)}A^{7-56\theta}}\right)^{\frac{1}{16(2-7\theta)}}\\
    &
    +\frac{1}{N^{3/28}A^{\frac{2-7\theta}{14-28\theta}}}
    +\left(\frac{|\vartheta|^3 A^{56\theta-5}}{N^{12(1-2\theta)}}\right)^{\frac{1}{16(5-14\theta)}}
    \Bigg).
\end{aligned}
\]
In particular, if $A\ll N^{1/14}$ and $|\vartheta|A^{11/3}\ll N^{4/3}$, the factor in parentheses is $\ll N^{-1/12}A^{-1/3}$.
\end{corollary}

\begin{proof}
By Cauchy--Schwarz, the sum is bounded by $\|\alpha\|\|\beta\|\|\nu\|N\sqrt A$, so the first term suffices if $|\vartheta|A>N^2$. We may therefore assume $A\le N^2$. Extending the three sequences by zero beyond their upper endpoints, we split the sum into dyadic ranges $m\sim M_1$, $n\sim M_2$, $a\sim B$, with $M_1,M_2\le N$ and $B\le A$. We start each decomposition with $(\tfrac12,1]$, so that it includes the index $1$.

Fix $L\ge1$. On each block we apply \cref{eq:spectral-unbalanced-amplifier}, interchanging $M_1,M_2$ on the right-hand side if necessary. As observed in the proof of \cref{thm:unbalanced-spectral-methods}, the expression in parentheses there is non-decreasing in each of the three lengths. We may thus replace $M_1,M_2,B$ by $N,N,A$, respectively, before taking square roots. The norms on each block are at most the norms of the full sequences, and there are $O(\log^2(2N)\log(2A))=N^{o(1)}$ blocks. After summing and extracting the factor $\|\alpha\|\|\beta\|\|\nu\|N^{1+o(1)}\sqrt A$, the term $N^{3/2}A/\sqrt L$ contributes $(NL)^{-1/4}\le L/N^{1/4}$. It therefore suffices to bound
\begin{equation}\label{eq:spectral-balanced-amplifier}
    \frac{\sqrt{|\vartheta|A}}{N}+\frac{1}{\sqrt{AL}}+\frac{1}{A^{1/4}L^{3/4}}
    +\frac{L}{N^{1/4}}\left(1+\frac{|\vartheta|A}{L^8}\right)^{1/16}
                       \left(1+\frac{N}{L^7}\right)^{\theta/2}.
\end{equation}
Up to absolute factors, the last term is the maximum over $i,j\in\{0,1\}$ of
\[
    N^{-1/4+j\theta/2}(|\vartheta|A)^{i/16}L^{1-i/2-7j\theta/2}.
\]
These four terms are increasing, since $1/2-7\theta/2>0$. At $L=1$ their maximum is $\asymp (|\vartheta|A)^{1/16}N^{-(1-2\theta)/4}$, giving the second term in the statement.

For the contributions from balancing these terms with the two decreasing terms, we use the calculation in the proof of \cref{thm:unbalanced-spectral-methods}. Balancing $(AL)^{-1/2}$ with the terms indexed by $(i,j)=(0,0),(1,0),(0,1),(1,1)$ gives, respectively,
\[
\begin{gathered}
    \frac{1}{N^{1/12}A^{1/3}},\qquad
    \left(\frac{|\vartheta|}{N^4 A^7}\right)^{1/32},\\
    \frac{1}{N^{\frac{1-2\theta}{12-28\theta}}
               A^{\frac{2-7\theta}{6-14\theta}}},\qquad
    \left(\frac{|\vartheta|}{N^{4(1-2\theta)}A^{7-56\theta}}\right)^{\frac{1}{16(2-7\theta)}}.
\end{gathered}
\]
For $A^{-1/4}L^{-3/4}$ it suffices to use the terms with $j=1$. Indeed, at $L=N^{1/7}$ we have
\[
    A^{-1/4}L^{-3/4}=A^{-1/4}N^{-3/28}\ll N^{-3/28}
    =N^{-1/4+\theta/2}L^{1-7\theta/2}.
\]
Hence the crossing of $A^{-1/4}L^{-3/4}$ with the maximum of the two terms with $j=1$ occurs at $L\ll N^{1/7}$. There the terms with $j=0$ are bounded by constant multiples of the corresponding terms with $j=1$. The two balances with $j=1$ give
\[
    \frac{1}{N^{3/28}A^{\frac{2-7\theta}{14-28\theta}}},
    \qquad
    \left(\frac{|\vartheta|^3 A^{56\theta-5}}{N^{12(1-2\theta)}}\right)^{\frac{1}{16(5-14\theta)}}.
\]
For each of the two decreasing terms, consider the value of $L>0$ at which it equals the maximum of the four increasing terms. As in the unbalanced case, we take the larger of these two values, increasing it to $1$ if necessary. If $L>1$, the maximum at this choice is bounded up to an absolute factor by the largest of the six quantities just displayed. If $L=1$, it is bounded by the contribution at $L=1$ computed above. Conversely, each of the six quantities and the contribution from the increasing terms at $L=1$ are bounded up to an absolute factor by the sum of the $L$-dependent terms in \cref{eq:spectral-balanced-amplifier} for every $L\ge1$. This proves the stated bound.

For the last assertion take $L=N^{1/6}A^{-1/3}$ in \cref{eq:spectral-balanced-amplifier}, increasing it to $1$ if necessary. In the stated range we have $L\gg A$, $N/L^7\ll1$ and $|\vartheta|A/L^8\ll1$. Thus the $L$-dependent terms are $\ll (AL)^{-1/2}+LN^{-1/4}\ll N^{-1/12}A^{-1/3}$. Also $\sqrt{|\vartheta|A}/N\ll N^{-1/3}A^{-4/3}\ll N^{-1/12}A^{-1/3}$.
\end{proof}

\begin{corollary}\label{cor:balanced-bilinear}
Let $N\in\Z_+$ and let $(\alpha_m)_{m=1}^N$ and $(\beta_n)_{n=1}^N$ be complex sequences. For any $\vartheta\in\Z\setminus\{0\}$, one has
\[
\begin{aligned}
    \mathop{\sum_{m=1}^N\sum_{n=1}^N}_{(m,n)=1}
       \alpha_m\beta_n\,e\left(\vartheta\frac{\bar m}{n}\right)
    \ll \|\alpha\|\|\beta\|N^{1+o(1)}
    \left(
       \frac{1}{N^{1/12}}
       + \frac{|\vartheta|^{1/32}}{N^{1/8}}
       + \frac{\sqrt{|\vartheta|}}{N}
       \right).
\end{aligned}
\]
\end{corollary}

\begin{proof}
By the trivial bound we may assume $|\vartheta|\le N^2$. We split $m,n$ into dyadic ranges as in the proof of \cref{cor:balanced-spectral-methods}. On each block we take the singleton $a=1$ (of dyadic length $\tfrac12$), with $\nu_1=1$, in \cref{eq:spectral-unbalanced-amplifier}, using \cref{eq:spectral-unbalanced-coupled} for the last term. This expression is again non-decreasing in both dyadic lengths. Bounding them by $N$, taking square roots and summing over the $O(\log^2(2N))$ blocks, we find that after extracting $\|\alpha\|\|\beta\|N^{1+o(1)}$ it suffices to bound
\[
    \frac{\sqrt{|\vartheta|}}{N}+\frac{1}{\sqrt L}
    +\frac{L}{N^{1/4}}
       \left(1+\frac{|\vartheta|}{L^8}
                +\frac{N^2}{L^{14}}\right)^{1/16}.
\]
We take
\[
    L=\min\left(N^{1/6},\frac{N^{1/4}}{|\vartheta|^{1/16}}\right).
\]
Since $|\vartheta|\le N^2$, we have $L\ge1$. The last term is
\[
    \ll \frac{L}{N^{1/4}}
       +\frac{|\vartheta|^{1/16}L^{1/2}}{N^{1/4}}
       +\frac{L^{1/8}}{N^{1/8}}.
\]
The two upper bounds on $L$ give $\ll L^{-1/2}$ for the first and second terms, respectively. For the third, $L^5\le N$ since $L\le N^{1/6}$. Thus all three terms are $\ll L^{-1/2}$, and
\[
    L^{-1/2}=\max\left(N^{-1/12},
                \frac{|\vartheta|^{1/32}}{N^{1/8}}\right),
\]
as claimed.
\end{proof}

\begin{proof}[Proof of \cref{thm:main-bilinear}]
Apply \cref{cor:balanced-bilinear} with $\vartheta=a$. For $|a|\ll N^{4/3}$ we have
\[
    \frac{|a|^{1/32}}{N^{1/8}}\ll N^{-1/12},
    \qquad \frac{\sqrt{|a|}}{N}\ll N^{-1/3},
\]
which gives the claimed bound.
\end{proof}

\begin{corollary} \label{cor:balanced-trilinear}
Let $N,A\in\Z_+$ and let $(\alpha_m)_{m=1}^N$, $(\beta_n)_{n=1}^N$, and $(\nu_a)_{a=1}^A$ be complex sequences. Then one has
\[
\begin{aligned}
    \mathop{\sum_{m=1}^N\sum_{n=1}^N}_{(m,n)=1} \sum_{a=1}^A
    \alpha_m \beta_n \nu_a\,
    e\left(a\frac{\bar m}{n}\right)
    &\ll \|\alpha\|\|\beta\|\|\nu\| N^{1+o(1)}\sqrt A
    \\
    &\times
    \Bigg(
        \frac{1}{N^{1/12}A^{1/3}}
        +\frac{1}{N^{\frac{1-2\theta}{12-28\theta}}
                   A^{\frac{2-7\theta}{6-14\theta}}}
        +\frac{1}{N^{3/28}
                   A^{\frac{2-7\theta}{14-28\theta}}}
        +\frac{\sqrt A}{N}
    \Bigg).
\end{aligned}
\]
\end{corollary}

\begin{proof}
For $A>N^2$ the last term suffices by the trivial bound, so we may assume $A\le N^2$. Take $\vartheta=1$ in \cref{eq:spectral-balanced-amplifier}, and choose $L>0$ so that
\[
    \max\left(\frac{1}{\sqrt{AL}},\frac{1}{A^{1/4}L^{3/4}}\right)
    =\max\left(\frac{L}{N^{1/4}},
                \frac{L^{1-7\theta/2}}{N^{1/4-\theta/2}}\right).
\]
As in the proof of \cref{thm:unbalanced-spectral-methods}, the common value is the largest of the four quantities obtained by equating a term on the left with a term on the right. These are the first three terms in parentheses in the statement and
\[
    \frac{1}{N^{3/28}A^{1/7}}
    \ll \frac{1}{N^{3/28}A^{\frac{2-7\theta}{14-28\theta}}},
\]
where we used $\theta\ge0$ and $A\ge1$.

If $L\ge A^{1/3}$, then $L\ge1$ and $A/L^8\le1$, so the sum of the $L$-dependent terms in \cref{eq:spectral-balanced-amplifier} is bounded by a constant times this common value. This gives the desired bound. Otherwise the common value is at least $A^{-1/4}L^{-3/4}\ge A^{-1/2}$. We split into dyadic ranges $m\sim M_1$, $n\sim M_2$, $a\sim B$ as in the proof of \cref{cor:balanced-spectral-methods}. On each block we use the second bound in \cref{lem:large-sieve-large-powers} with $C=1$. Since $\vartheta=1$, the coprimality condition with $\vartheta$ is automatic, and the proof of the lemma applies to every $B$. Apart from the squared norms and a factor of $N^{o(1)}$, its bound for the square of the sum is
\[
    (M_1+M_2)(M_2+B)\ll N(N+A).
\]
Taking square roots and summing over the $O(\log^2(2N)\log(2A))$ blocks gives
\[
    \ll \|\alpha\|\|\beta\|\|\nu\|N^{1+o(1)}\sqrt A
       \left(A^{-1/2}+N^{-1/2}\right).
\]
Since $N^{-1/2}\le A^{-1/2}$ for $A\le N$ and $N^{-1/2}\le\sqrt A/N$ for $A \ge N$, this also gives the stated bound.
\end{proof}

\begin{proof}[Proof of \cref{thm:main-trilinear}]
Take $\vartheta=1$ in the last assertion of \cref{cor:balanced-spectral-methods}. The hypothesis $A\ll N^{1/14}$ gives $A^{11/3}\ll N^{11/42}\ll N^{4/3}$, so the asserted bound is
\[
    \ll \|\alpha\|\|\beta\|\|\nu\|N^{1+o(1)}\sqrt A\,
          N^{-1/12}A^{-1/3}
    =\|\alpha\|\|\beta\|\|\nu\|N^{11/12+o(1)}A^{1/6},
\]
which completes our proof.
\end{proof}

\section{Applications} \label{sec:applications}

\subsection{Twisted moments of zeta}
With $P_Y$ and $w$ as in \cref{subsec:applications}, we allow complex shifts $\alpha,\beta\ll1/\log T$ in the two zeta factors and write
\[
    I(\alpha,\beta):=\int_{\R}
        \zeta(\tfrac12+\alpha+it)\zeta(\tfrac12+\beta-it)
        |P_Y(\tfrac12+it)|^2\,w(t/T)\,dt.
\]
The implied constants may also depend on the constant in the bounds for the shifts.

\begin{theorem} \label{thm:appl-twisted-zeta}
Uniformly for $1\le Y\le T^{207/389}$ and $\alpha,\beta\ll1/\log T$, one has
\[
\begin{aligned}
    I(\alpha,\beta)
    &=\sum_{m,n\le Y}\frac{c_m\overline{c_n}}{[m,n]}
        \frac{(m,n)^{\alpha+\beta}}{m^\beta n^\alpha}\\
    &\quad\times\int_{\R}\left(
             \zeta(1+\alpha+\beta)
             +\zeta(1-\alpha-\beta)
              \left(\frac{2\pi mn}{t(m,n)^2}\right)^{\alpha+\beta}
        \right)w(t/T)\,dt\\
    &\quad +O\left(T^{o(1)}\left(
        T^{1/2}Y^{1/2}+T^{1/3}Y^{5/4}
        +T^{79/286}Y^{389/286}\right)\right).
\end{aligned}
\]
At $\alpha+\beta=0$ the expression in parentheses in the integral is interpreted by continuity.
\end{theorem}

\begin{proof}[Proof of \cref{thm:appl-twisted-zeta}]
We use the reduction of Bettin--Chandee--Radziwi{\l}{\l} \cite[\S3]{bettin2017mean} in the shifted form of Pratt--Robles \cite[\S3]{pratt2018perturbed}, referring to their proofs for most of the details. 

We first suppose that $|\alpha+\beta|\asymp1/\log T$. The diagonal and the zero frequency of the off-diagonal together give the stated main term. To see the placement of the shifts, the diagonal for the coefficient $c_m\overline{c_n}$ requires the two zeta indices to be $\ell n/(m,n)$ and $\ell m/(m,n)$, respectively. The other terms, apart from the nonzero frequencies in \cite[\S3.2.3]{pratt2018perturbed} and the corresponding terms from the dual sum, contribute $O(T^{1/2+o(1)}Y^{1/2})$.

We recall the ranges in the remaining sum, see also \cite[(5.1)]{dong2026bilinear}. After extracting a common divisor $d$, the two indices $n_1,n_2$ are coprime and have respective sizes $N_1/d,N_2/d$, where $N_1,N_2\ll Y$. If $M$ denotes the dyadic length of the zeta variable in the Poisson summation, the numerator of the Kloosterman fraction is $h\Delta$. For any fixed $\eps>0$, the frequency sums may be truncated to
\[
    0<|h|\ll\frac{N_2 T^\eps}{dM},
    \qquad 0<|\Delta|\ll\frac{MN_1}{dT^{1-\eps}},
\]
with negligible error. We group the pairs with the same product and split it according to its sign and its dyadic size $A$. Thus $A\ll T^{-1+2\eps}N_1N_2/d^2$, and the coefficient of each product is bounded by $O(\tau(|h\Delta|))$. For each choice of the continuous variables, we have a sum of the form
\[
    \sum_{a\sim A}\nu_a
    \sum_{\substack{n_1\asymp N_1/d,\ n_2\asymp N_2/d\\(n_1,n_2)=1}}
        b_{n_1}b'_{n_2}\,e\left(\pm a\frac{\bar n_1}{n_2}\right).
\]
We take the Mellin inversion of $w$ on $\Re s=0$, which is permissible since $w$ is smooth and supported in $[1,2]$. The resulting imaginary powers of $n_1,n_2$ have absolute value $1$, so the individual $b_{n_1},b'_{n_2}$ are $\ll T^{o(1)}$ by the coefficient hypothesis and the bounds on the shifts. Hence
\[
    \|\nu\|\|b\|\|b'\|\ll T^{o(1)}\sqrt A\,\frac{\sqrt{N_1N_2}}{d}.
\]
The factor outside the integral is $1/d$. Before Mellin inversion of $w$, the two real integration variables are restricted to intervals of lengths $O(dM/N_2)$ and $O(Td/(MN_1))$, respectively. We retain this rectangle when taking absolute values. The integral of the absolute value of the Mellin transform of $w$ on $\Re s=0$ is $\ll(\log T)^2$ by two integrations by parts. Thus the cost of the integrations, including the outside factor, is $\ll T^{1+o(1)}d/(N_1N_2)$. These bounds hold for both sums in the approximate functional equation.

Put $X=\max(N_1,N_2)/d$. Extending the three sequences by zero, we apply \cref{cor:balanced-trilinear} with upper bounds $\asymp X$ for the two indices and $\asymp A$ for the numerator. Apart from a factor of $T^{o(1)}$, a term $X^{-s}A^{-t}$ in the parentheses of that corollary bounds the contribution for fixed $d,N_1,N_2,M,A$ by
\[
\begin{aligned}
    \frac{Td}{N_1N_2}\,
    \sqrt A\,\frac{\sqrt{N_1N_2}}{d}\,
    X^{1-s}A^{1/2-t}
    &\ll T^{t+2\eps(1-t)}
        \frac{(N_1N_2)^{1/2-t}\max(N_1,N_2)^{1-s}}
             {d^{3-s-2t}}\\
    &\ll T^{t+2\eps(1-t)}Y^{2-s-2t}d^{-3+s+2t}.
\end{aligned}
\]
Here the four pairs $(s,t)$ are
\[
    (\tfrac{1}{12},\tfrac13),\qquad
    (\tfrac{25}{286},\tfrac{79}{286}),\qquad
    (\tfrac{3}{28},\tfrac{79}{700}),\qquad
    (1,-\tfrac12).
\]
In particular $1/2-t\ge0$, $1-s\ge0$, and the sum over $d$ converges in all four cases. Summing also over the logarithmically many dyadic choices and using the fact that the truncations hold for every fixed $\eps>0$, we obtain the error
\[
    \ll T^{o(1)}\left(
        T^{1/3}Y^{5/4}+T^{79/286}Y^{389/286}
        +T^{79/700}Y^{1167/700}+T^{-1/2}Y^2\right)
\]
from the nonzero frequencies under consideration. Since $Y^{389}\le T^{207}$, we have
\[
    \frac{T^{79/700}Y^{1167/700}}{T}
    =\left(\frac{Y^{389}}{T^{207}}\right)^{3/700}
    \le\left(\frac{Y^{389}}{T^{207}}\right)^{1/286}
    =\frac{T^{79/286}Y^{389/286}}{T}.
\]
Also $T^{-1/2}Y^2\le T^{1/2}Y^{1/2}$, proving the required error estimate. More generally, before substituting $\theta=7/64$, the first three terms of \cref{cor:balanced-trilinear} require, for $Y=T^\sigma$, that $\sigma<8/15$ and $\sigma<(8-14\theta)/(15-26\theta)$, the latter restriction occurring twice. Thus Selberg's eigenvalue conjecture would give every $\sigma<8/15$. At the limiting unconditional length $Y=T^{207/389}$, the natural length of the sum over $h\Delta$ is $Y^2/T=Y^{25/207}$. This is where we use the range of \cref{cor:balanced-trilinear} beyond \cref{thm:main-trilinear}.

Finally fix $|\beta|\le C/\log T$. The difference between the moment and its main term is holomorphic in $\alpha$ on the disc $|\alpha+\beta|\le(2C+1)/\log T$, since the two zeta poles in the integrand of the main term cancel. On the boundary the proof above applies uniformly, with $|\alpha|\le(3C+1)/\log T$. The maximum modulus principle gives the same estimate for every $|\alpha|\le C/\log T$.
\end{proof}

\begin{proof}[Proof of \cref{cor:appl-twisted-zeta}]
Take $\alpha=\beta=0$ and $Y=T^\sigma$ in \cref{thm:appl-twisted-zeta}. Since $\zeta(1+z)=z^{-1}+\gamma+O(z)$ as $z\to0$, for $t$ in the support of $w(t/T)$ we have
\[
    \lim_{z\to0}\left(\zeta(1+z)+\zeta(1-z)
        \left(\frac{2\pi mn}{t(m,n)^2}\right)^z\right)
    =\log\frac{t(m,n)^2}{2\pi mn}+2\gamma.
\]
The factor $(m,n)^{\alpha+\beta}/(m^\beta n^\alpha)$ is $1$ at these shifts, giving the main term in \cref{cor:appl-twisted-zeta}. The three exponents of $T$ in the error term are, up to $o(1)$,
\[
    \frac{1+\sigma}{2},\qquad \frac13+\frac{5\sigma}{4},\qquad
    \frac{79+389\sigma}{286}.
\]
Since $\sigma<207/389<8/15<1$, all three displayed exponents are less than $1$. Taking $\delta$ to be half the smallest of their distances from $1$ gives the corollary.
\end{proof}

\subsection{Convolutions in APs} \label{subsec:appl-convolutions}

Following Fouvry--Radziwi{\l}{\l} \cite{fouvry2022level}, we say that $(\beta_n)$ satisfies the Siegel--Walfisz condition if there is an integer $k_0\ge1$ such that, for every $B>0$,
\[
    \sum_{\substack{n\sim y\\n\equiv b\pmod r\\(n,s)=1}}\beta_n
    -\frac{1}{\varphi(r)}\sum_{\substack{n\sim y\\(n,rs)=1}}\beta_n
    \ll_B \tau_{k_0}(s)y(\log y)^{-B}
\]
uniformly for $y\ge2$, $s\ge1$ and $1\le |b|<r$ with $(b,r)=1$. 

\begin{theorem} \label{thm:appl-convolutions}
Let $k\ge1$ be an integer, $\eps>0$, $M,N\ge2$, $Q\ge1$ and $MN/2\le x\le4MN$. Let $(\alpha_m)_{m\sim M}$ and $(\beta_n)_{n\sim N}$ be complex sequences with $|\alpha_m|\le\tau_k(m)$ and $|\beta_n|\le\tau_k(n)$, and suppose that $(\beta_n)$ satisfies the Siegel--Walfisz condition. Let $a$ be a nonzero integer with $|a|\le x/12$. Then, for every $B>0$, the bound \cref{eq:appl-convolutions} holds provided that $\exp((\log x)^\eps)\le N$ and
\begin{equation} \label{eq:appl-convolution-range}
    Q\le x^{-\eps}\min\left\{
        \left(\frac{x^{207}}{N^{468}}\right)^{1/389},
        \left(\frac{x^{151}}{|a|^8N^{292}}\right)^{1/277},
        \frac{x^{3/4}}{|a|^{1/4}N}
    \right\}.
\end{equation}
%The implied constant may also depend on the constants in the Siegel--Walfisz condition.
\end{theorem}

\begin{proof}[Proof of \cref{thm:appl-convolutions}]
Write $X=MN\asymp x$. We may use $X$ in place of $x$ in the power bounds of \cref{eq:appl-convolution-range}, decreasing the margin $\eps$ to $\eps/2$ if necessary. First suppose that $N\ge X^{1/53}$. Then
\[
    X^{207/389}N^{-468/389}
    \le \min\left(\sqrt{XN},X^{4/7}N^{-6/7}\right).
\]
Indeed the first inequality is equivalent to $X\le N^{53}$, and the second follows from $207/389<4/7$ and $468/389>6/7$. In this range the result follows from Fouvry \cite[Corollaire~1]{fouvry1987autour}, in the form used in the proof of \cite[Proposition~8.1(ii)]{fouvry2022level}, with the restriction $x<mn\le2x$ inserted as in \cite[\S8]{fouvry2022level}.

We may now suppose that $N\le X^{1/53}$. We follow the proof of \cite[Theorem~1.1 and Corollary~1.1(i)]{fouvry2022level}, using the numbering of the arXiv version. We take $\eta>0$ sufficiently small in terms of $\eps$ and set
\[
    D=\exp\big((\log X)^{\min(\eps,1)/2}/20\big).
\]
We have $M>N>D^{10}$ and $Q\le MX^{-\eta}$ by the first bound in \cref{eq:appl-convolution-range}. In the dispersion estimate it suffices to prove $W^{\mathrm{Err1}}\ll MN^2 X^{-c}$ for some $c>0$. Indeed, the other terms under the square root in \cite[Theorem~1.1]{fouvry2022level} are
\[
    \ll \frac{M}{Q}\mathcal E^\star(\beta,N,Q)
       +(\log X)^{O_k(1)}\left(N^2Q+D^{-1/2}MN^2\right).
\]
By the Siegel--Walfisz condition and \cite[Lemma~8.1]{fouvry2022level}, $\mathcal E^\star(\beta,N,Q)\ll_B N^2Q(\log X)^{-B}$ for every $B$, since $\log N\gg(\log X)^\eps$. Thus these terms are $\ll_B MN^2(\log X)^{-B}$. Taking square roots and multiplying by $\|\alpha\|\ll\sqrt M(\log X)^{O_k(1)}$ gives the desired logarithmic saving.

We use the trilinear sum for $W^{\mathrm{Err1}}$ in \cite[(48)--(49)]{fouvry2022level}, also recorded in \cite[(5.10)]{dong2026bilinear}. After fixing $d,d_1\le D$, $\gamma\le D^3$ and summing over $\nu_1',\nu_2\le2N$, the phase is
\[
    e\left(\vartheta h\frac{\bar m}{n}\right),\qquad
    \vartheta=a(d_1\nu_1'-\nu_2)\ne0,\quad
    m=\gamma dd_1\nu_2 k_1',\quad n=\nu_1'k_2'.
\]
Here $|\vartheta|\ll D|a|N$, $k_1',k_2'\ll Q$, $(m,n)=1$, and $1\le h\ll(\log X)^4 Q^2/M$ after splitting according to sign. The remaining restrictions are incorporated in the three coefficient sequences as in \cite[(48)]{fouvry2022level}. For $h\sim A$ their norms are $\ll\sqrt A,\sqrt Q,\sqrt Q$, respectively. In particular, the multiples of $\gamma dd_1\nu_2$ and $\nu_1'$ supporting the last two sequences cost nothing in their norms. The two lengths are $\ll D^5NQ$. There are $O(N^2)$ choices for $\nu_1',\nu_2$, and the factor in front of these sums is $\ll D^{O(1)}X^{o(1)}M/Q^2$, including the remaining choices and dyadic subdivisions.

Put $H=Q^2/M$. If there are nonzero frequencies then $H\gg(\log X)^{-4}$. We split $m,n$ into dyadic ranges and apply \cref{eq:spectral-balanced-amplifier} on each pair, keeping $\vartheta$ as the multiplier. By the monotonicity in the two lengths observed in the proof of \cref{thm:unbalanced-spectral-methods}, we may use $D^5NQ$ as their common upper bound up to an absolute factor. After including the norms and the preceding factors and dividing by $MN^2$, we obtain, for any $L\ge1$,
\begin{equation} \label{eq:appl-dispersion}
\begin{aligned}
    \frac{|W^{\mathrm{Err1}}|}{MN^2}
    \ll X^{o(1)} NH\Bigg(
        &\frac{\sqrt{|a|NH}}{NQ}
         +\frac{1}{\sqrt{HL}}+\frac{1}{H^{1/4}L^{3/4}}\\
        &+\frac{L}{(NQ)^{1/4}}
          \left(1+\frac{|a|NH}{L^8}\right)^{1/16}
          \left(1+\frac{NQ}{L^7}\right)^{7/128}
    \Bigg).
\end{aligned}
\end{equation}
Indeed, for each block the two factors of $\sqrt A$ combine to give $A$, and $MQ^{-2}\cdot N^2\cdot Q\cdot NQ\cdot A=MN^2\cdot NA$. After multiplication by $A$ all terms are increasing in $A$, so we can use $A\ll H(\log X)^4$ throughout. All powers of $D$ and $\log X$ are absorbed by $X^{o(1)}$.

Take $L=X^\eta N^2 H$ in \cref{eq:appl-dispersion}. The two decreasing terms, including the factor $NH$, give $X^{-\eta/2}$ and $N^{-1/2}X^{-3\eta/4}$. The first term gives
\[
    \frac{\sqrt{|a|}(NQ)^2}{X^{3/2}},
\]
which has a power saving by the last bound in \cref{eq:appl-convolution-range}. The remaining term, again including $NH$, is at most
\[
    X^\eta\frac{N^{19/4}Q^{15/4}}{X^2}
      \left(1+\frac{|a|X^7}{N^{22}Q^{14}}\right)^{1/16}
      \left(1+\frac{X^7}{N^{20}Q^{13}}\right)^{7/128}.
\]
By $Q\ll X^{207/389}N^{-468/389}$ and $N\le X^{1/53}$, we have
\[
    \frac{X^7}{N^{20}Q^{13}}\gg
    \left(\frac{X}{N^{53}}\right)^{32/389}\ge1.
\]
Expanding the factor involving $|a|$, we can therefore bound the remaining term by
\[
    \ll X^\eta\left\{
        \left(\frac{Q^{389}N^{468}}{X^{207}}\right)^{1/128}
       +\left(\frac{|a|^8Q^{277}N^{292}}{X^{151}}\right)^{1/128}
    \right\}.
\]
Both terms have a power saving by the first two bounds in \cref{eq:appl-convolution-range}, on choosing $\eta$ small enough. This proves the required estimate for $W^{\mathrm{Err1}}$. Finally, the passage from the dyadic rectangle to $x<mn\le2x$ in \cite[\S8]{fouvry2022level} only involves smoothing the endpoints and twisting $\alpha_m,\beta_n$ by $m^{it},n^{it}$ with $|t|\le(\log X)^{O_B(1)}$. The norms are preserved, and partial summation preserves the Siegel--Walfisz bounds with an arbitrarily large logarithmic saving. This proves \cref{eq:appl-convolutions}.
\end{proof}

\begin{proof}[Proof of \cref{cor:appl-convolutions}]
For $x$ sufficiently large we have $x^{175/389}\le x/12$. To apply \cref{thm:appl-convolutions}, evaluate the second and third restrictions of \cref{eq:appl-convolution-range} at $Q=x^{207/389}N^{-468/389}$ before inserting the factor $x^{-\eps}$. They become $|a|\le x^{175/389}N^{2006/389}$ and $|a|\le x^{339/389}N^{316/389}$, respectively, both of which hold if $|a|\le x^{175/389}$.
\end{proof}

\subsection{A determinant equation}
Here we prove the following generalization of \cref{cor:appl-determinant}.

\begin{theorem} \label{thm:appl-determinant}
Let $N\ge2$ and $\Delta\in\Z\setminus\{0\}$. With $S_\Delta(N)$, $f,g$, and the two sequences as in \cref{subsec:applications}, one has
\[
\begin{aligned}
    \Bigg|S_\Delta(N)-
    \sum_{\substack{n_1,n_2\sim N\\(n_1,n_2)\mid\Delta}}
       \frac{(n_1,n_2)}{n_1n_2}\alpha_{n_1}\beta_{n_2}
       \int_{\R} f\left(\frac{x+\Delta}{Nn_2}\right)
                 g\left(\frac{x}{Nn_1}\right)\,dx\Bigg|
    \ll_{f,g}\|\alpha\|\|\beta\|N^{o(1)}\\
    \times\min\left\{N^{19/20},\
        N^{11/12}+N^{7/8}|\Delta|^{1/32}
        +|\Delta|^{1/2}\right\}.
\end{aligned}
\]
\end{theorem}

\begin{proof}[Proof of \cref{thm:appl-determinant}]
We may assume $|\Delta|\le3N^2$, since otherwise the supports of $f,g$ make both the sum and the main term vanish. For fixed $n_1,n_2$, write $d=(n_1,n_2)$, $n_1=du$ and $n_2=dv$. The determinant equation requires $d\mid\Delta$, and then it determines $m_1=(um_2+\Delta/d)/v$ with $m_2\equiv-(\Delta/d)\bar u\pmod v$. Poisson summation in $m_2$ writes $S_\Delta(N)$ as
\[
    \sum_{\substack{d\mid\Delta\\d>0}}\sum_{h\in\Z}
    \sum_{\substack{u,v\sim N/d\\(u,v)=1}}
        \alpha_{du}\beta_{dv}\,
        e\left(-\frac{h\Delta}{d}\frac{\bar u}{v}\right)
        \frac1v\int_{\R}
            f\left(\frac{uy+\Delta/d}{vN}\right)
            g(y/N)e(-hy/v)\,dy.
\]
For $h=0$ the change of variables $x=n_1y$ gives exactly the main term, since $1/(vn_1)=d/(n_1n_2)$.

For $h\ne0$, the change of variables $y=Nz$ writes the smooth weight as $K_{d,h}(du/N,dv/N)$, where
\[
    K_{d,h}(s,t)=\frac{d}{t}\int_{\R}
        f\left(\frac{sz+\Delta/N^2}{t}\right)g(z)
        e\left(-\frac{hdz}{t}\right)\,dz .
\]
Uniformly for $1\le s,t\le2$, $0\le i,j\le1$ and any $B\ge0$, we have
\[
    \partial_s^i\partial_t^j K_{d,h}(s,t)
    \ll_{f,g,B} d(1+|h|d)^{-B}.
\]
Indeed $|\Delta|/N^2\le3$, differentiation costs at most one power of $|h|d$ here, and we may integrate by parts arbitrarily often in $z$. By partial summation in $u,v$, it therefore suffices to bound the bilinear form with $u,v$ restricted to arbitrary subintervals of $(N/d,2N/d]$, at a cost of $O_{f,g,B}(d(1+|h|d)^{-B})$.

We apply \cref{cor:balanced-bilinear} with upper endpoint $\lfloor 2N/d\rfloor$ and $\vartheta=-h\Delta/d$. We extend the two sequences by zero outside their subintervals, so their norms are at most $\|\alpha\|,\|\beta\|$. Since $2N/d\le2N$ and $|\vartheta|\le |h\Delta|$, the bilinear form is $\ll\|\alpha\|\|\beta\|N^{o(1)}$ times
\[
    (1+|h|)^{1/2}
        \left(N^{11/12}+N^{7/8}|\Delta|^{1/32}+|\Delta|^{1/2}\right).
\]
Taking $B=4$, the sum of $d(1+|h|d)^{-4}(1+|h|)^{1/2}$ over all $d\ge1$ and $h\ne0$ converges. This gives the second entry in the minimum. The first follows from \cite[Corollary~1]{bettin2018trilinear} with all four lengths equal to $2N$ in its convention and with smooth weights $f(\cdot/N),g(\cdot/N)$.
\end{proof}

\begin{proof}[Proof of \cref{cor:appl-determinant}]
If $|\Delta|\ll N^{4/3}$, then $N^{7/8}|\Delta|^{1/32}\ll N^{7/8+1/24}=N^{11/12}$ and $|\Delta|^{1/2}\ll N^{2/3}$. The result then follows from \cref{thm:appl-determinant}.
\end{proof}

\let\oldaddcontentsline\addcontentsline
\renewcommand{\addcontentsline}[3]{}
\bibliographystyle{alpha}
\bibliography{main}
\let\addcontentsline\oldaddcontentsline

\end{document}